\documentclass[11pt]{article}%
\usepackage{amssymb}
\usepackage{amsfonts}
\usepackage{amsmath}
\usepackage{graphicx}%
\providecommand{\U}[1]{\protect\rule{.1in}{.1in}}
\newtheorem{theorem}{Theorem}

\newtheorem{corollary}{Corollary}

\newtheorem{property}{Property}

\begin{document}

\title{Binet's factorial series and extensions to Laplace transforms}
\author{Piet Van Mieghem\thanks{Faculty of Electrical Engineering, Mathematics and
Computer Science, P.O Box 5031, 2600 GA Delft, The Netherlands; \emph{email}:
P.F.A.VanMieghem@tudelft.nl. }}
\date{Delft University of Technology\\
11 February 2023 - for Nathan (5/12/1994-11/2/2012)}
\maketitle

\begin{abstract}
We investigate a generalization of Binet's factorial series in the parameter
$\alpha$
\[
\mu\left(  z\right)  =\sum_{m=1}^{\infty}\frac{b_{m}\left(  \alpha\right)
}{\prod_{k=0}^{m-1}(z+\alpha+k)}%
\]
due to Gilbert, for the Binet function
\[
\mu\left(  z\right)  =\log\Gamma\left(  z\right)  -\left(  z-\frac{1}%
{2}\right)  \log z+z-\frac{1}{2}\log\left(  2\pi\right)
\]
After a review of the Binet function $\mu\left(  z\right)  $ and Gilbert's
investigations of $\mu\left(  z\right)  $, several properties of the Binet
polynomials $b_{m}\left(  \alpha\right)  $ are presented. We compare Gilbert's
generalized factorial series with Stirling's \emph{asymptotic} expansion and
demonstrate by a numerical example that, with a same number of terms
evaluated, the Gilbert generalized factorial series with an optimized value of
$\alpha$ can beat the best possible accuracy of Stirling's expansion. Finally,
we extend Binet's method to factorial series of Laplace transforms.

\end{abstract}

\section{Introduction}

In his truly comprehensive\footnote{The reference in Whittaker and Watson
\cite[footnote, p. 248]{Whittaker_Watson} with page numbers from 123 to 143
suggests an ordinary article, but the correct page numbers pp. 123-343 point
to a book-size treatise. Apart from adding his own beautiful contributions,
Binet has reviewed the knowledge about the Beta-function before 1839. The
focal point at that time were the many properties of Euler's Beta integral
from which properties of the Gamma function were derived. Today, on the other
hand and perhaps after Weierstrass's product and Hankel's contour integral,
the theory concentrates on the Gamma function, with applications to the Beta
function.} \textquotedblleft memoir\textquotedblright\ \cite[Section 3, p.
223]{Binet1839} in 1839, Jacques Binet defines his function $\mu\left(
z\right)  $ in relation to the Gamma function $\Gamma\left(  z\right)  $ as
\begin{equation}
\mu\left(  z\right)  =\log\Gamma\left(  z\right)  -\left(  z-\frac{1}%
{2}\right)  \log z+z-\frac{1}{2}\log\left(  2\pi\right)
\label{def_Binet_mu(z)}%
\end{equation}
Binet has derived two integral representations of his function \cite[p.
248-251]{Whittaker_Watson}, \cite[p. 211-217]{Sansone}, for $\operatorname{Re}%
\left(  z\right)  >0$,
\begin{equation}
\mu\left(  z\right)  =2\int_{0}^{\infty}\frac{\arctan\left(  \frac{t}%
{z}\right)  }{e^{2\pi t}-1}dt \label{Binet_mu_integral_arctan}%
\end{equation}
and%
\begin{equation}
\mu\left(  z\right)  =\frac{1}{2}\int_{0}^{\infty}\frac{e^{-zt}}{t}\left(
\frac{1+e^{-t}}{1-e^{-t}}-\frac{2}{t}\right)  dt
\label{Binet_mu_integral_Bernoulli_generating_function}%
\end{equation}

There exist more representations\footnote{Blagouchine \cite{Blagouchine_2016}
lists 7 different formulae for $\log\Gamma\left(  z\right)  $.} of $\mu\left(
z\right)  $, but here we merely concentrate on elegant, converging factorial
series for $\operatorname{Re}\left(  z\right)  >0$ due to Binet\footnote{In
Binet's notation \cite[p. 234]{Binet1839},
\begin{equation}
2\mu\left(  z\right)  =\frac{I\left(  1\right)  }{z+1}+\frac{I\left(
2\right)  }{2\left(  z+1\right)  \left(  z+2\right)  }+\frac{I\left(
3\right)  }{3\left(  z+1\right)  \left(  z+2\right)  \left(  z+3\right)
}+\cdots\label{Binet_wrong}%
\end{equation}
where the integral form of the coefficients (derived at \cite[p.
238]{Binet1839}) is%
\begin{equation}
I\left(  m\right)  =\int_{0}^{1}x\left(  x+1\right)  \left(  x+2\right)
\cdots\left(  x+m-1\right)  \left(  2x-1\right)  dx
\label{Binet_coefficient_integral_wrong}%
\end{equation}
Binet \cite[p. 234]{Binet1839} lists the first few coefficients, $I\left(
1\right)  =\frac{1}{6}$, $I\left(  2\right)  =\frac{1}{3}$, $I\left(
3\right)  =\frac{59}{60}$ and $I\left(  4\right)  =\frac{227}{60}$.} in
\cite[p. 234]{Binet1839},
\begin{equation}
\mu\left(  z\right)  =\sum_{m=1}^{\infty}\frac{\beta_{m}}{\left(  z+1\right)
\left(  z+2\right)  \cdots\left(  z+m\right)  } \label{Binet_original}%
\end{equation}
where the coefficients are%
\begin{equation}
\beta_{m}=\frac{1}{m}\int_{0}^{1}\left(  u-\frac{1}{2}\right)  u\left(
u+1\right)  \cdots\left(  u+m-1\right)  du \label{Binet_coefficient_original}%
\end{equation}
Explicitly, $\beta_{1}=\beta_{2}=\frac{1}{12}$, $\beta_{3}=\frac{59}{360}$,
$\beta_{4}=\frac{29}{60}$, $\beta_{5}=\frac{533}{280}$, $\beta_{6}=\frac
{1577}{168}$, $\beta_{7}=\frac{280361}{5040}$, $\beta_{8}=\frac{69311}{180}$,
and $\beta_{m}$ is positive and rapidly increasing in $m$. Nearly at the end
of his memoir and somewhat hidden, Binet \cite[p. 342]{Binet1839} has given a
second factorial series (\ref{Binet_function_faculteits_expansie}), that is
rederived slightly differently in Theorem
\ref{theorem_Binet_convergent_expansion} in Section \ref{sec_Binet_expansion}
and generalized to Laplace transforms in Section \ref{sec_1_factorial_series}.

Recently, Nemes \cite[Theorem 2.1]{Nemes_2013} has generalized Binet's
expansion (\ref{Binet_original}), for $0\leq a\leq1$ and $\operatorname{Re}%
\left(  z\right)  >0$,%
\begin{equation}
\log\Gamma\left(  z+a\right)  =\left(  z+a-\frac{1}{2}\right)  \log
z-z+\frac{1}{2}\log\left(  2\pi\right)  +\sum_{m=1}^{\infty}\frac{c_{m}\left(
a\right)  }{\left(  z+1\right)  \left(  z+2\right)  \cdots\left(  z+m\right)
} \label{LogGamma_factorial_expansion_Nemes}%
\end{equation}
where%
\begin{equation}
c_{m}\left(  a\right)  =\frac{1}{m}\int_{0}^{1}\left(  u+a-\frac{1}{2}\right)
u\left(  u+1\right)  \cdots\left(  u+m-1\right)  du-\frac{1}{m}\int_{0}%
^{a}\prod_{j=0}^{m-1}\left(  u-a+1+j\right)  du
\label{LogGamma_factorial_expansion_coeff_Nemes}%
\end{equation}
Clearly, if $a=0$, we retrieve Binet's first factorial expansion
(\ref{Binet_original}) and coefficients $\beta_{m}$ in
(\ref{Binet_coefficient_original}). Nemes' expansion
(\ref{LogGamma_factorial_expansion_Nemes}), expressed in terms of Binet's
function $\mu\left(  z\right)  $ with the definition (\ref{def_Binet_mu(z)}),%
\[
\mu\left(  z\right)  =\left(  z-\frac{1}{2}\right)  \log\frac{z-a}{z}%
+a+\sum_{m=1}^{\infty}\frac{c_{m}\left(  a\right)  }{\prod_{k=1}^{m}(z-a+k)}%
\]
bears resemblance to
\begin{equation}
\mu\left(  z\right)  =\sum_{m=1}^{\infty}\frac{b_{m}\left(  \alpha\right)
}{\prod_{k=0}^{m-1}(z+\alpha+k)} \label{mu_factorial_alpha_expansions}%
\end{equation}
in our main Theorem \ref{theorem_mu_infinite_convergent_factorial_expansions},
where the coefficients $b_{m}\left(  \alpha\right)  $ are called Binet
polynomials, to honour Jacques Binet. Much earlier, Hermite \cite{Hermite1895}
has deduced the corresponding generalized Stirling asymptotic expansion
\begin{equation}
\log\Gamma\left(  z+a\right)  \simeq\left(  z+a-\frac{1}{2}\right)  \log
z-z+\frac{1}{2}\log\left(  2\pi\right)  +\sum_{k=1}^{K}\frac{\left(
-1\right)  ^{k-1}B_{k+1}\left(  a\right)  }{k\left(  k+1\right)  z^{k}}
\label{LogGamma_generalized_Stirling_Hermite}%
\end{equation}
in terms of the Bernoulli polynomials $B_{n}\left(  a\right)  $ that reduces
for $a=0$ to Stirling's original asymptotic series
(\ref{Binet_asymptotic_Stirling_expansion}). Starting from a complex integral
(\ref{Binet_complex_integral_Zeta}) for Binet's function $\mu\left(  z\right)
$, Stirling's asymptotic series (\ref{Binet_asymptotic_Stirling_expansion}) is
derived in Section \ref{sec_Stirling_asymptotic_series}, where also the
meaning of the upper bound $K$ is explained. Section
\ref{sec_Stirling_convergent_companion} presents the \emph{convergent}
companion (\ref{mu_series_TaylorZeta_1_op_z}) of Stirling's \emph{divergent}
series (if $K\rightarrow\infty$ in (\ref{Binet_asymptotic_Stirling_expansion})).

The main motivation that led to this article is twofold. Originally, I was
confused about Binet's achievements: he derived \emph{two} \emph{different}
factorial expansions (\ref{Binet_original}) and
(\ref{Binet_function_faculteits_expansie}) for the \emph{same} function
$\mu\left(  z\right)  $. While I thought initially that one of them must be
wrong, I discovered later with (\ref{mu_factorial_alpha_expansions}) that
infinitely many different factorial series exist. However, Gilbert
\cite{Gilbert1873} has anticipated me about 150 years earlier. The second
motivation was my unbelief that Stirling's asymptotic \emph{but divergent}
expansion seems unbeatable in performance for some optimal, finite $K$.

Rediscoveries seem to appear frequently in mathematics. In earlier versions, I
wrote that \textquotedblleft our main result is the demonstration that there
exist infinitely many different factorial expansions in a complex parameter
$\alpha$ for Binet's function $\mu\left(  z\right)  $\textquotedblright. After
Gerg\H{o} Nemes informed me in January 2023 about Gilbert's investigations in
\cite{Gilbert1873}, I have reoriented the article by integrating the beautiful
discoveries of predessors with my own findings in Theorem
\ref{theorem_mu_infinite_convergent_factorial_expansions} (Section
\ref{sec_infinitely_many_factorial series}) and for Laplace transforms in
Theorem
\ref{theorem_Laplace_transform_infinite_convergent_factorial_expansions}
(Section \ref{sec_infinitely_many_factorial_series_Laplance_transforms}). My
approach is different from Gilbert's and an integrated presentation provides a
broader view on Binet's function $\mu\left(  z\right)  $. For $\alpha=1$ in
(\ref{mu_factorial_alpha_expansions}), we recover Binet's first factorial
series (\ref{Binet_original}) with $b_{m}\left(  1\right)  =\beta_{m}$ in
(\ref{Binet_coefficient_original}), while $\alpha=0$ in
(\ref{mu_factorial_alpha_expansions}) corresponds to Binet's second factorial
series (\ref{Binet_function_faculteits_expansie}) with $b_{m}\left(  0\right)
=c_{m}$ in (\ref{def_Binet_coefficients}). Interestingly, the $\alpha=1$
factorial series has all positive coefficients and any truncation thus lower
bounds Binet's function $\mu\left(  z\right)  $, while the $\alpha=0$
factorial series is shown in Theorem \ref{theorem_Binet_coefficient_integral}
to possess coefficients $c_{m}<0$ for $m>2$. Thus, a truncation of
(\ref{Binet_function_faculteits_expansie}) upper bounds $\mu\left(  z\right)
$. For $\alpha$ around $\frac{1}{10}$, which is close to the largest zero of
the Binet polynomial $b_{m}\left(  \alpha\right)  $, numerical computations
exhibit the fastest convergence. With the same number $K$ of terms evaluated,
the variant $\alpha=0$ is more accurate than the variant $\alpha=1$. Perhaps,
the slower convergence of Binet's expansion (\ref{Binet_original}) has led to
its omission in handbooks of functions, like Abramowitz \& Stegun
\cite{Abramowitz} nor in its successor by Olver \emph{et al}. \cite{Olver}.

We will first discuss the main properties of Binet's function $\mu\left(
z\right)  $ in Section \ref{sec_Binet_function} and the deductions from the
complex integral (\ref{Binet_complex_integral_Zeta}) in Section
\ref{sec_complex_integral_Binet_function_mu}, before we review parts of
Binet's great treatise. In section \ref{sec_Binet_first_factorial_expansion},
we sketch Binet's route towards his first factorial series
(\ref{Binet_original}), that is covered in the literature (see e.g. \cite[p.
253]{Whittaker_Watson}, \cite[p. 30]{Paris_Kaminski}). Binet's second
factorial series, that I have not found in later works, is derived in more
detail in subsection \ref{sec_Binet_expansion}, also because I believe that,
being the case $\alpha=0$ in (\ref{mu_factorial_alpha_expansions}), it is
slightly more important than his first series (\ref{Binet_original})
corresponding to $\alpha=1$. Moreover, Binet's method towards his second
factorial series, which is a recipe in five steps, enables a far reaching
generalization to Laplace transforms as explained in Section
\ref{sec_factorial_series_Laplace_transforms}. Section \ref{sec_Gilbert}
reviews Gilbert's remarkable investigations in \cite{Gilbert1873}. I have
slightly generalized in Section \ref{sec_Gilbert_generalized_factorial_series}
his derivation of a generalized factorial series in
(\ref{Gilbert_factorial_series_Binet_functions}) for Binet's function
$\mu\left(  z\right)  $. Section \ref{sec_infinitely_many_factorial series}
presents my derivation of Gilbert's generalized factorial series for $p=1$ in
(\ref{Gilbert_factorial_series_Binet_functions}) and properties of its
coefficients that I have called the Binet polynomial $b_{m}\left(
\alpha\right)  $. Factorial expansions for the derivatives $\frac{d^{n}%
\mu\left(  z\right)  }{dz^{n}}$ are derived in Section
\ref{sec_derivatives_mu_digamma_polygamma} and applied to the digamma and
polygamma function. In particular for the digamma function $\psi\left(
z\right)  =\frac{\Gamma^{\prime}\left(  z\right)  }{\Gamma\left(  z\right)  }%
$, we thus add a convergent series (\ref{digamma_factorial_expansion}) to its
asymptotic counterpart (\ref{digamma_asymptotic}). Section
\ref{sec_Stirling_Binet} discusses and compares, with a same number of terms,
the accuracy of Stirling's \emph{asymptotic} expansion and the best possible
that can be attained by the generalized Binet factorial expansion. The
commonly accepted belief about the superiority of Stirling's asymptotic
expansion is demonstrated and plotted in Fig.
\ref{Fig_Stirlingoptimalaccuracy}. However, at the zeros of the Binet
polynomial $b_{m}\left(  \alpha\right)  $, the accuracy of the generalized
Binet factorial expansion (\ref{mu_factorial_alpha_expansions}) improves
considerably as drawn in Fig. \ref{Fig_Binetalpha}. By a numerical example, we
show that Stirling's series accuracy is not always better than a factorial
series (with a same number of terms)! In other words, the generalized Binet
expansion (\ref{mu_factorial_alpha_expansions}) can be optimized with respect
to the \textquotedblleft free\textquotedblright\ parameter $\alpha$ to
achieve, at least, a comparable accuracy with a same computational effort.
Perhaps, this observation deserves to list the generalized Binet factorial
expansion (\ref{mu_factorial_alpha_expansions}) in handbooks of functions.

Computations are deferred to the appendices in order to enhance the
readability and focus on the essential parts.

\section{Binet's function $\mu\left(  z\right)  $}

\subsection{Properties deduced from the definition (\ref{def_Binet_mu(z)})}

\label{sec_Binet_function}The definition (\ref{def_Binet_mu(z)}) of Binet's
function $\mu\left(  z\right)  $ directly shows, for a positive integer $z=n$,
that%
\[
\mu\left(  n\right)  =n+\log\left(  n-1\right)  !-\left(  n-\frac{1}%
{2}\right)  \log\left(  n\right)  -\frac{1}{2}\log\left(  2\pi\right)
\]
The sequence $\mu\left(  1\right)  \simeq0.0811$, $\mu\left(  2\right)
\simeq0.0413$, $\mu\left(  3\right)  \simeq0.0277$, $\mu\left(  4\right)
\simeq0.0208$, $\mu\left(  5\right)  \simeq0.0166$,\ldots, $\mu\left(
10\right)  \simeq0.0083$ demonstrates the slow decay roughly as $\mu\left(
n\right)  \approx\frac{1}{12n}$. The precise decay is given in
(\ref{Binet_explicit}) below.

\medskip\emph{Maximum at real, positive values of }$z$. Binet's integral
(\ref{Binet_mu_integral_Bernoulli_generating_function}) can be rewritten as a
Laplace transform%
\begin{equation}
\mu\left(  z\right)  =\int_{0}^{\infty}\frac{e^{-zt}}{t}\left(  \frac{1}%
{e^{t}-1}-\frac{1}{t}+\frac{1}{2}\right)  dt
\label{Binet_mu_integral_Bernoulli_generating_function_2}%
\end{equation}
where $0\leq\frac{1}{e^{t}-1}-\frac{1}{t}+\frac{1}{2}\leq\frac{1}{2}$ for
real, non-negative $t$. Since the integrand is positive for real $z=x>0$, we
observe that $\mu\left(  x\right)  >0$ for real, positive $x$. In addition,
for a complex number $z=x+iy$, the above integral shows that $\mu\left(
z\right)  $ is analytic for $\operatorname{Re}\left(  z\right)  >0$ and that
\[
\left\vert \mu\left(  z\right)  \right\vert \leq\int_{0}^{\infty}\frac
{e^{-xt}}{t}\left(  \frac{1}{e^{t}-1}-\frac{1}{t}+\frac{1}{2}\right)
dt=\mu\left(  x\right)
\]
In other words, the maximum absolute value of Binet's function $\mu\left(
z\right)  $ for $\operatorname{Re}\left(  z\right)  >0$ is attained at the
positive real axis. Moreover, $\mu\left(  x\right)  $ strictly decreases with
$x=\operatorname{Re}\left(  z\right)  $. Another rather straightforward bound
is%
\[
\left\vert \mu\left(  z\right)  \right\vert \leq\max_{0\leq t}\frac{1}%
{t}\left(  \frac{1}{e^{t}-1}-\frac{1}{t}+\frac{1}{2}\right)  \int_{0}^{\infty
}e^{-xt}dt=\frac{1}{x}\max_{0\leq t}\frac{1}{t}\left(  \frac{1}{e^{t}-1}%
-\frac{1}{t}+\frac{1}{2}\right)
\]
Since the maximum occurs at $t=0$, the generating function of the Bernoulli
numbers $B_{n}$
\begin{equation}
\frac{1\,}{e^{t}-1}=\frac{1}{t}-\frac{1}{2}+\sum_{n=1}^{\infty}B_{2n}%
\,\frac{t^{2n-1}}{(2n)!}\hspace{1cm}\text{for }\left\vert t\right\vert
\leq2\pi\label{gf_Bernoullinumbers}%
\end{equation}
illustrates that $\max_{0\leq t}\frac{1}{t}\left(  \frac{1}{e^{t}-1}-\frac
{1}{t}+\frac{1}{2}\right)  =\frac{B_{2}}{2}=\frac{1}{12}$. Hence, for
$x=\operatorname{Re}\left(  z\right)  >0$, we find \cite[p. 249]%
{Whittaker_Watson}%
\begin{equation}
\left\vert \mu\left(  z\right)  \right\vert \leq\mu\left(  x\right)  \leq
\frac{1}{12x} \label{bound_mu_Re(z)>0}%
\end{equation}

\medskip\emph{Difference }$\mu\left(  z+1\right)  -\mu\left(  z\right)  $.
Binet's definition (\ref{def_Binet_mu(z)})%
\[
\mu\left(  z\right)  =\log\Gamma\left(  z\right)  -\left(  z-\frac{1}%
{2}\right)  \log z+z-\frac{1}{2}\log\left(  2\pi\right)
\]
illustrates that the complex conjugate $\mu^{\ast}\left(  z\right)
=\mu\left(  z^{\ast}\right)  $, by the reflection principle \cite[p.
155]{Titchmarshfunctions}. The functional equation of the Gamma function
$\Gamma\left(  z+1\right)  =z\Gamma\left(  z\right)  $ leads to%
\[
\mu\left(  z+1\right)  =\log\Gamma\left(  z\right)  +\log z-\left(  z+\frac
{1}{2}\right)  \log\left(  z+1\right)  +z+1-\frac{1}{2}\log\left(
2\pi\right)
\]
The forward difference $\Delta\mu(z)=\mu\left(  z+1\right)  -\mu\left(
z\right)  $ equals%
\begin{equation}
\mu\left(  z+1\right)  -\mu\left(  z\right)  =\left(  z+\frac{1}{2}\right)
\log\frac{z}{z+1}+1 \label{forwards_difference_mu}%
\end{equation}
which is valid for any complex $z$, except at the negative real axis that is a
branch cut for $\mu\left(  z\right)  $. In particular, around $z=0$, the
forward difference (\ref{forwards_difference_mu}) shows, with $\mu\left(
1\right)  =1-\frac{1}{2}\log\left(  2\pi\right)  $, that%
\begin{equation}
\mu\left(  z\right)  \sim-\frac{1}{2}\log z-\frac{1}{2}\log\left(
2\pi\right)  \label{mu_around_z=0}%
\end{equation}
illustrating that Binet's function $\mu\left(  z\right)  $ possesses a
logarithmic singularity at $z=0$.

Erdelyi et al. \cite[p. 24]{Erdelyi_v1} deduce from Gauss's multiplication
formula that%
\[
\int_{z}^{z+1}\log\Gamma\left(  t\right)  dt=z\log z-z+\frac{1}{2}\log\left(
2\pi\right)
\]
which we rewrite with the definition (\ref{def_Binet_mu(z)}) as%
\begin{equation}
\int_{z}^{z+1}\mu\left(  t\right)  dt=\frac{1}{2}\left(  z\left(  z+1\right)
\log\left(  \frac{z}{z+1}\right)  +z+\frac{1}{2}\right)
\label{forwards_difference_mu_integrated}%
\end{equation}
Differentiation of (\ref{forwards_difference_mu_integrated}) with respect to
$z$ again leads to the forward difference (\ref{forwards_difference_mu}). For
$z=0$ in (\ref{forwards_difference_mu_integrated}), we find $\int_{0}^{1}%
\mu\left(  t\right)  dt=\frac{1}{4}$.

\emph{Gudderman's series}. If we replace $z$ by $z+k$ in the forward
difference (\ref{forwards_difference_mu}) and change the sign, then
(\ref{forwards_difference_mu}) becomes $\mu\left(  z+k\right)  -\mu\left(
z+k+1\right)  =\left(  z+k+\frac{1}{2}\right)  \log\left(  1+\frac{1}%
{z+k}\right)  -1$. Summing over integer $k$ results in a telescoping series
$\sum_{k=0}^{K}\left\{  \mu\left(  z+k\right)  -\mu\left(  z+k+1\right)
\right\}  =\mu\left(  z\right)  -\mu\left(  z+K+1\right)  $ leading to%
\[
\mu\left(  z\right)  -\mu\left(  z+K+1\right)  =\sum_{k=0}^{K}\left\{  \left(
z+k+\frac{1}{2}\right)  \log\left(  1+\frac{1}{z+k}\right)  -1\right\}
\]
When $K$ tends to infinity, the bound (\ref{bound_mu_Re(z)>0}) shows that
$\lim_{K\rightarrow\infty}\mu\left(  z+K+1\right)  =0$ and we obtain
Gudermann's series \cite{Gudermann1845}%
\begin{equation}
\mu\left(  z\right)  =\sum_{k=0}^{\infty}\left\{  \left(  z+k+\frac{1}%
{2}\right)  \log\left(  1+\frac{1}{z+k}\right)  -1\right\}
\label{Gudermann_series_mu}%
\end{equation}

\medskip\emph{Reflection formula of Binet's function }$\mu\left(  z\right)  $.
We replace $z\rightarrow1-z$ in Binet's definition (\ref{def_Binet_mu(z)})%
\[
\mu\left(  1-z\right)  =\log\Gamma\left(  1-z\right)  -\left(  \frac{1}%
{2}-z\right)  \log\left(  1-z\right)  -z+1-\frac{1}{2}\log\left(  2\pi\right)
\]
which added to $\mu\left(  z\right)  $ in (\ref{def_Binet_mu(z)}), yields
\[
\mu\left(  z\right)  +\mu\left(  1-z\right)  =\log\left(  \Gamma\left(
z\right)  \Gamma\left(  1-z\right)  \right)  -\left(  z-\frac{1}{2}\right)
\log\frac{z}{1-z}+1-\log\left(  2\pi\right)
\]
After invoking the reflection formula of the Gamma function $\Gamma\left(
z\right)  \Gamma\left(  1-z\right)  =\frac{\pi}{\sin\pi z}$, we find the
corresponding \textquotedblleft reflection\textquotedblright\ formula for
Binet's function%
\begin{equation}
\mu\left(  1-z\right)  =1-\mu\left(  z\right)  -\log\left(  2\sin\pi z\right)
-\left(  z-\frac{1}{2}\right)  \log\frac{z}{1-z}
\label{Binet_mu_reflection_formula}%
\end{equation}
which is valid for any complex number $z=x+iy$, with the exception of the
negative real axis and odd integers $x=2n+1$ with $n\in\mathbb{Z}$ at any $y$.
Since $\mu\left(  z\right)  $ is analytic for $\operatorname{Re}\left(
z\right)  >0$, the Binet reflection formula (\ref{Binet_mu_reflection_formula}%
) illustrates that Binet's function $\mu\left(  z\right)  $ has only
logarithmic singularities at negative integers $z=-k$ (with $k\in\mathbb{N}$)
including $z=0$, as shown in (\ref{mu_around_z=0}).

\medskip\emph{Duplication and multiplication formula for }$\mu\left(
z\right)  $. Combining the duplication formula $\Gamma\left(  2z\right)
=\frac{2^{2z-1}}{\sqrt{\pi}}\Gamma\left(  z\right)  \Gamma\left(  z+\frac
{1}{2}\right)  $ in \cite[6.1.18]{Abramowitz} for the Gamma function
$\Gamma\left(  z\right)  $ and the definition (\ref{def_Binet_mu(z)}) of
Binet's function $\mu\left(  z\right)  $ leads to%
\begin{equation}
\mu\left(  2z\right)  =\mu\left(  z+\frac{1}{2}\right)  +\mu\left(  z\right)
+z\log\left(  1+\frac{1}{2z}\right)  -\frac{1}{2}
\label{Binet_mu_duplication_formula}%
\end{equation}
which is generalized by Gauss's multiplication formula $\Gamma\left(
nz\right)  =\left(  2\pi\right)  ^{\frac{1}{2}\left(  1-n\right)  }%
n^{nz-\frac{1}{2}}%
%TCIMACRO{\dprod \limits_{k=0}^{n-1}}%
%BeginExpansion
{\displaystyle\prod\limits_{k=0}^{n-1}}
%EndExpansion
\Gamma\left(  z+\frac{k}{n}\right)  $ in \cite[6.1.20]{Abramowitz} as%
\begin{equation}
\mu\left(  nz\right)  =\mu\left(  z\right)  +\sum_{k=1}^{n-1}\mu\left(
z+\frac{k}{n}\right)  +\sum_{k=1}^{n-1}\left(  z+\frac{k}{n}-\frac{1}%
{2}\right)  \log\left(  z+\frac{k}{n}\right)  -\left(  n-1\right)  \left(
\frac{1}{2}+z\log z\right)  \label{Binet_mu_multiplication_formula}%
\end{equation}
After choosing $z=\frac{1}{n}$ in (\ref{Binet_mu_multiplication_formula}), we
find for $n\geq2$%
\[
\mu\left(  \frac{1}{n}\right)  =-\sum_{k=1}^{n-2}\mu\left(  \frac{k+1}%
{n}\right)  -\sum_{k=1}^{n-2}\left(  \frac{k+1}{n}-\frac{1}{2}\right)
\log\left(  \frac{k+1}{n}\right)  +\left(  n-1\right)  \left(  \frac{1}%
{2}-\frac{1}{n}\log n\right)
\]
For example, for $n=2$, we obtain $\mu\left(  \frac{1}{2}\right)  =\frac{1}%
{2}\left(  1-\log2\right)  \simeq0.1534$ but, for $n=3$, $\mu\left(  \frac
{1}{3}\right)  =1-\mu\left(  \frac{2}{3}\right)  -\frac{1}{6}\log\left(
2\right)  -\frac{1}{2}\log\left(  3\right)  $, which is also immediate from
Binet's reflection formula (\ref{Binet_mu_reflection_formula}). The
duplication formula (\ref{Binet_mu_duplication_formula}) yields a closed
expression for $z=n+\frac{1}{2}$ with integer $n\geq1$,%
\[
\mu\left(  n+\frac{1}{2}\right)  =\sum_{k=n}^{2n-1}\log k-n\left(  \log\left(
4n+2\right)  -1\right)  +\frac{1}{2}\left(  \log\left(  2\right)  +1\right)
\]

\subsection{\textbf{Gilbert's infinite product for }$\Gamma\left(  z\right)
$\textbf{ }}

Gilbert \cite[art. 3]{Gilbert1873} further investigates Gudermann's series
(\ref{Gudermann_series_mu}). After remarking that%
\begin{align*}
\sum_{k=0}^{K}\left(  z+k+\frac{1}{2}\right)  \log\left(  1+\frac{1}%
{z+k}\right)   &  =\sum_{k=0}^{K}\left(  z+k+\frac{1}{2}\right)  \left(
\log\left(  z+k+1\right)  -\log\left(  z+k\right)  \right) \\
&  =\sum_{k=1}^{K+1}\left(  z+k-\frac{1}{2}\right)  \log\left(  z+k\right)
-\sum_{k=0}^{K}\left(  z+k+\frac{1}{2}\right)  \log\left(  z+k\right)
\end{align*}
and reworking the (telescoping) sums,%
\[
\sum_{k=0}^{K}\left(  z+k+\frac{1}{2}\right)  \log\left(  1+\frac{1}%
{z+k}\right)  =\left(  z+K+\frac{1}{2}\right)  \log\left(  z+K+1\right)
-\left(  z+\frac{1}{2}\right)  \log\left(  z\right)  -\sum_{k=1}^{K}%
\log\left(  z+k\right)
\]
Gilbert \cite[art. 3]{Gilbert1873} finds%
\begin{align*}
\mu\left(  z\right)  -\mu\left(  z+K+1\right)   &  =\sum_{k=0}^{K}\left\{
\left(  z+k+\frac{1}{2}\right)  \log\left(  1+\frac{1}{z+k}\right)  -1\right\}
\\
&  =-\left(  z+\frac{1}{2}\right)  \log\left(  z\right)  +\left(  z+K+\frac
{1}{2}\right)  \log\left(  z+K+1\right)  -\left(  K+1\right)  -\sum_{k=1}%
^{K}\log\left(  z+k\right)
\end{align*}
and proceeding to the limit $K\rightarrow\infty$ yields%
\[
\mu\left(  z\right)  =-\left(  z+\frac{1}{2}\right)  \log z+\lim
_{K\rightarrow\infty}\left\{  \left(  z+K+\frac{1}{2}\right)  \log\left(
z+K+1\right)  -\left(  K+1\right)  -\sum_{k=1}^{K}\log\left(  z+k\right)
\right\}
\]
Since
\begin{align*}
\left(  z+K+\frac{1}{2}\right)  \log\left(  z+K+1\right)   &  =\left(
z+K+\frac{1}{2}\right)  \log\left(  K\left(  1+\frac{z+1}{K}\right)  \right)
\\
&  =\left(  z+K+\frac{1}{2}\right)  \log K+\left(  z+\frac{1}{2}\right)
\log\left(  1+\frac{z+1}{K}\right)  +K\log\left(  1+\frac{z+1}{K}\right)
\end{align*}
it holds that%
\[
\mu\left(  z\right)  =-\left(  z+\frac{1}{2}\right)  \log z+z+\lim
_{K\rightarrow\infty}\left\{  \left(  z+K+\frac{1}{2}\right)  \log
K-K-\sum_{k=1}^{K}\log\left(  z+k\right)  \right\}
\]
The definition $\log\Gamma\left(  z\right)  =\mu\left(  z\right)  +\left(
z-\frac{1}{2}\right)  \log z-z+\frac{1}{2}\log\left(  2\pi\right)  $ in
(\ref{def_Binet_mu(z)}) then shows that
\[
\log\Gamma\left(  z\right)  =\log\left(  \sqrt{2\pi}\right)  +\lim
_{K\rightarrow\infty}\left\{  \left(  z+K+\frac{1}{2}\right)  \log
K-K-\log\prod_{k=0}^{K}\left(  z+k\right)  \right\}
\]
After exponentiation, Gilbert \cite[art. 3]{Gilbert1873} arrives at%
\begin{equation}
\Gamma\left(  z\right)  =\sqrt{2\pi}\lim_{n\rightarrow\infty}\frac
{n^{z+n+\frac{1}{2}}e^{-n}}{\prod_{k=0}^{n}\left(  z+k\right)  }
\label{Gilbert_product_Gamma_function}%
\end{equation}
which is another product form than the Gauss product
\begin{equation}
\Gamma\left(  z\right)  =\frac{1}{z}\prod_{n=1}^{\infty}\left(  1+\frac{1}%
{n}\right)  ^{z}\left(  1+\frac{z}{n}\right)  ^{-1}
\label{Gauss_product_Gamma_function}%
\end{equation}
and the Weierstrass product in (\ref{Weierstrass_product_Gamma_function})
below. Gilbert \cite[art. 4]{Gilbert1873} adds that Stirling's formula for
large $n$%
\[
n!=\sqrt{2\pi}n^{n+\frac{1}{2}}e^{-n}\left(  1+O\left(  \frac{1}{n}\right)
\right)
\]
transforms his product (\ref{Gilbert_product_Gamma_function}) with
$n^{z}=\prod_{k=1}^{n-1}\left(  \frac{k+1}{k}\right)  ^{z}$ to the Gauss
product (\ref{Gauss_product_Gamma_function}).

From the exponentiation $\Gamma\left(  z\right)  =\sqrt{2\pi}z^{z-\frac{1}{2}%
}e^{-z+\mu\left(  z\right)  }$ of definition (\ref{def_Binet_mu(z)}) and using
the bound (\ref{bound_mu_Re(z)>0}), Gilbert \cite{Gilbert1886} derives, for
$z<n$, the bounds
\[
\frac{\sqrt{2\pi}}{\Gamma\left(  z\right)  }n^{z+n-\frac{1}{2}}e^{-n}%
e^{\left(  z-1\right)  \frac{z}{2n}-\left(  z-\frac{1}{2}\right)  \frac{z^{2}%
}{2n^{2}}}<z\left(  z+1\right)  \ldots\left(  z+n-1\right)  <\frac{\sqrt{2\pi
}}{\Gamma\left(  z\right)  }n^{z+n-\frac{1}{2}}e^{-n}e^{\left(  z-\frac{1}%
{2}\right)  \frac{z}{n}+\frac{1}{12n\left(  n+1\right)  }}%
\]

\subsection{\medskip Complex integral for Binet's function $\mu\left(
z\right)  $}

\label{sec_complex_integral_Binet_function_mu}In Appendix
\ref{sec_complex_integral_mu}, we deduce the complex integral (in two ways)
\begin{equation}
\mu\left(  z\right)  =-\frac{1}{2\pi i}\int_{c-i\infty}^{c+i\infty}\frac{\pi
}{\sin\pi s}\frac{\zeta\left(  s\right)  }{s}z^{s}ds\hspace{1cm}\text{with
}-1<c<0 \label{Binet_complex_integral_Zeta}%
\end{equation}
valid for any complex number $z=\left\vert z\right\vert e^{i\arg z}$ with
$\left\vert \arg z\right\vert <\pi$. We substitute the functional equation of
the Riemann Zeta-function $\zeta(s)=2(2\pi)^{s-1}\sin\frac{\pi s}{2}%
\Gamma(1-s)\zeta(1-s)$ in the integral (\ref{Binet_complex_integral_Zeta}),
invoke the reflection formula $\Gamma\left(  s\right)  \Gamma\left(
1-s\right)  =\frac{\pi}{\sin\pi s}$ and obtain the variant of
(\ref{Binet_complex_integral_Zeta})%
\begin{equation}
\mu\left(  z\right)  =-\frac{1}{4i}\int_{c-i\infty}^{c+i\infty}\frac{(2\pi
z)^{s}}{\cos\frac{\pi s}{2}}\frac{\zeta\left(  1-s\right)  }{\sin\pi
s\Gamma\left(  1+s\right)  }ds\hspace{1cm}\text{with }-1<c<0
\label{Binet_complex_integral_Zeta_funceq}%
\end{equation}
The variant (\ref{Binet_complex_integral_Zeta_funceq}) allows the introduction
of the Dirichlet series $\zeta\left(  s\right)  =\sum_{k=1}^{\infty}\frac
{1}{k^{s}}$ for $\operatorname{Re}\left(  s\right)  >1$ and leads, after
evaluating the resulting contour integrals, to Malmsten-Kummer's series (see
e.g. \cite{Blagouchine_2016}, \cite[p. 23]{Erdelyi_v1}) for real $0<x<1$%
\begin{equation}
\log\Gamma\left(  x\right)  =\left(  \gamma+\log\left(  2\pi\right)  \right)
\left(  \frac{1}{2}-x\right)  -\frac{1}{2}\log\frac{\sin\pi x}{\pi}+\sum
_{k=2}^{\infty}\frac{\log(k)\sin\left(  2\pi xk\right)  }{\pi k}
\label{Malmsten_Kummer_series_LogGama(x)}%
\end{equation}

We evaluate the integral (\ref{Binet_complex_integral_Zeta}) along the line
$s=c+it$, where $-1<c<0$,%
\[
\mu\left(  z\right)  =-\frac{1}{2}\int_{-\infty}^{\infty}\frac{z^{c+it}}%
{\sin\pi\left(  c+it\right)  }\frac{\zeta\left(  c+it\right)  }{c+it}dt
\]
If we choose $c=-\frac{1}{2}$, then the integral simplifies to%
\begin{equation}
\mu\left(  z\right)  =z^{-\frac{1}{2}}\int_{-\infty}^{\infty}\frac{e^{it\log
z}}{\cosh\pi t}\frac{\zeta\left(  -\frac{1}{2}+it\right)  }{-1+2it}dt
\label{Binet_complex_integral_Zeta_along_line}%
\end{equation}
Since $\zeta\left(  x+iT\right)  =O\left(  T^{\frac{1}{2}-x}\right)  $ for
$x\leq0$ and large $T$ (see e.g. \cite[Chapter V]{Titchmarshzeta}), it holds
that $\frac{\zeta\left(  -\frac{1}{2}+it\right)  }{-1+2it}=O\left(  1\right)
$ for large $t$. Hence, the integral
(\ref{Binet_complex_integral_Zeta_along_line}) can be bounded as%
\[
\left\vert \mu\left(  z\right)  \right\vert \leq\left\vert z^{-\frac{1}{2}%
}\right\vert \int_{-\infty}^{\infty}\frac{\left\vert e^{it\log z}\right\vert
}{\cosh\pi t}\left\vert \frac{\zeta\left(  -\frac{1}{2}+it\right)  }%
{-1+2it}\right\vert dt\leq C\left\vert z^{-\frac{1}{2}}\right\vert
\int_{-\infty}^{\infty}\frac{e^{-t\arg\left(  z\right)  }}{\cosh\pi t}dt
\]
where $C$ is positive real number, demonstrating existence for all complex
$z=\left\vert z\right\vert e^{i\arg z}$ provided $\left\vert \arg z\right\vert
<\pi$. In other words, the integral (\ref{Binet_complex_integral_Zeta})
defines Binet's function $\mu\left(  z\right)  $ everywhere in the complex
plane, except at the negative real axis, where $\mu\left(  z\right)  $
possesses a branch cut. The well-known Fourier integral%
\[
\int_{-\infty}^{\infty}\frac{e^{izu}}{\cosh\alpha u}\,du=\frac{\pi}{\alpha
}\frac{1}{\cosh\frac{\pi z}{2\alpha}}\hspace{1cm}\text{valid for }\left\vert
\operatorname{Im}z\right\vert <\alpha
\]
shows that $\int_{-\infty}^{\infty}\frac{e^{it\log z}}{\cosh\pi t}dt=\frac
{1}{\cosh\frac{\log z}{2}}=\frac{2z^{\frac{1}{2}}}{z+1}$ and, roughly, that
the the integral (\ref{Binet_complex_integral_Zeta_along_line}) can be
estimated as $\mu\left(  z\right)  =O\left(  \frac{1}{z}\right)  $, which
complements the bound (\ref{bound_mu_Re(z)>0}) to complex $z$, except at the
negative real axis.

\medskip\emph{Branch cut along the negative real axis}. The complex integral
(\ref{Binet_complex_integral_Zeta}) indicates with $z=re^{i\theta}$ that%
\[
\mu\left(  re^{i\theta}\right)  -\mu\left(  re^{-i\theta}\right)
=-\int_{c-i\infty}^{c+i\infty}\frac{\sin\theta s}{\sin\pi s}\frac{\zeta\left(
s\right)  }{s}r^{s}ds
\]
is purely imaginary, because $\mu\left(  re^{i\theta}\right)  -\mu\left(
re^{-i\theta}\right)  =\mu\left(  re^{i\theta}\right)  -\mu^{\ast}\left(
re^{i\theta}\right)  =2i\operatorname{Im}\mu\left(  re^{i\theta}\right)  $,
and that%
\[
\lim_{\theta\rightarrow\pi}\left(  \mu\left(  re^{i\theta}\right)  -\mu\left(
re^{-i\theta}\right)  \right)  =-\int_{c-i\infty}^{c+i\infty}\frac
{\zeta\left(  s\right)  }{s}r^{s}ds
\]
By moving the line of integration from $c<0$ to $c^{\prime}>1$, two poles at
$s=0$ and $s=1$ are enclosed and Cauchy's residue theorem leads to%
\begin{align*}
\mu\left(  re^{i\pi}\right)  -\mu\left(  re^{-i\pi}\right)   &  =2\pi i\left(
\lim_{s\rightarrow0}\zeta\left(  s\right)  r^{s}+\lim_{s\rightarrow1}%
\frac{\zeta\left(  s\right)  \left(  s-1\right)  }{s}r^{s}\right)
-\int_{c^{\prime}-i\infty}^{c^{\prime}+i\infty}\frac{\zeta\left(  s\right)
}{s}r^{s}ds\\
&  \hspace{0.5cm}+\lim_{T\rightarrow\infty}\int_{c}^{c^{\prime}}\frac
{\zeta\left(  x+iT\right)  }{x+iT}r^{x+iT}dx-\lim_{T\rightarrow\infty}\int
_{c}^{c^{\prime}}\frac{\zeta\left(  x-iT\right)  }{x-iT}r^{x-iT}dx
\end{align*}
Due to $\zeta\left(  x+iT\right)  =O\left(  T^{\frac{1}{2}-x}\right)  $ for
$x\leq0$ and large $T$ and since $c<0$ can be chosen small enough, the limits
$T\rightarrow\infty$ vanish. After using Perron's formula \cite[p.
301]{Titchmarshfunctions}, $\frac{1}{2\pi i}\int_{c^{\prime}-i\infty
}^{c^{\prime}+i\infty}\frac{\zeta\left(  s\right)  }{s}r^{s}ds=\sum_{n=1}%
^{r}1=\left[  r\right]  $, which is the integral part of $r$, we find that the
difference at both sides of the branch cut is periodic in $r\geq0$,%
\[
\mu\left(  re^{i\pi}\right)  -\mu\left(  re^{-i\pi}\right)  =2\pi i\left(
-\frac{1}{2}+r-\left[  r\right]  \right)
\]
and vanishes at $r=\frac{1}{2}+n$ with integer $n\in\mathbb{N}$.

From the complex integral (\ref{Binet_complex_integral_Zeta}), we will deduce
Stirling's asymptotic series in Section \ref{sec_Stirling_asymptotic_series}
and its convergent companion in Section
\ref{sec_Stirling_convergent_companion}.

\subsection{Stirling's asymptotic series}

\label{sec_Stirling_asymptotic_series}We cannot close the contour in
(\ref{Binet_complex_integral_Zeta}) over the entire $\operatorname{Re}\left(
s\right)  <0$ -plane, because the functional equation of the Riemann
Zeta-function $\zeta(s)=2(2\pi)^{s-1}\sin\frac{\pi s}{2}\Gamma(1-s)\zeta(1-s)$
indicates that $\zeta(-s)=O\left(  \Gamma(s)\right)  $. However, neglecting
this restriction and using $\zeta(-k)=\frac{(-1)^{k}}{k+1}B_{k+1}$ and the odd
Bernoulli numbers $B_{2k+1}=0$, for $k>0$, leads to Stirling's asymptotic
approximation\footnote{By introducing the generating function
(\ref{gf_Bernoullinumbers}) of the Bernoulli numbers in Binet's integral
(\ref{Binet_mu_integral_Bernoulli_generating_function})
\[
\mu\left(  z\right)  =\int_{0}^{\infty}\frac{e^{-zt}}{t}\left(  \frac{1}%
{e^{t}-1}-\frac{1}{t}+\frac{1}{2}\right)  dt=\int_{0}^{\infty}e^{-zt}\left(
\sum_{m=1}^{\infty}B_{2m}\,\frac{t^{2m-2}}{(2m)!}\right)  dt
\]
only valid for $\operatorname{Re}\left(  z\right)  >0$ and reversing sum and
integral, while ignoring the convergence restriction $\left\vert t\right\vert
<2\pi$ in the sum, we obtain again Stirling's asymptotic series
(\ref{Binet_asymptotic_Stirling_expansion}).} \cite[6.1.41]{Abramowitz} in the
Poincar\'{e} sense (see e.g. \cite{Paris_Kaminski})
\begin{equation}
\mu\left(  z\right)  \simeq\sum_{k=1}^{K}\frac{\left(  -1\right)  ^{k}%
\zeta\left(  -k\right)  }{k}z^{-k}=\sum_{k=1}^{K}\frac{B_{k+1}}{k\left(
k+1\right)  z^{k}}=\sum_{m=1}^{K}\frac{B_{2m}}{\left(  2m-1\right)  \left(
2m\right)  z^{2m-1}} \label{Binet_asymptotic_Stirling_expansion}%
\end{equation}
Although (\ref{Binet_asymptotic_Stirling_expansion}) diverges if
$K\rightarrow\infty$, Fig. \ref{Fig_Stirlingoptimalaccuracy} in Section
\ref{sec_Stirling_Binet} below shows that Stirling's asymptotic approximation
(\ref{Binet_asymptotic_Stirling_expansion}) is surprisingly accurate up to
some finite $K\leq K^{\ast}\left(  z\right)  $, where $K^{\ast}\left(
z\right)  $ depends upon $z$ and is roughly equal to the minimum $k$-value of
$\frac{\left\vert B_{k+1}\right\vert }{k\left(  k+1\right)  \left\vert
z\right\vert ^{k}}$.

On the other hand for $\left\vert z\right\vert <1$, the contour in
(\ref{Binet_complex_integral_Zeta}) can be closed over the $\operatorname{Re}%
\left(  s\right)  >0$ -plane, where two double poles at $s=0$ and $s=1$ are
encountered whose residues are computed in Appendix
\ref{sec_complex_integral_mu}, resulting in%
\begin{equation}
\mu\left(  z\right)  =\sum_{k=2}^{\infty}\frac{\left(  -1\right)  ^{k}%
\zeta\left(  k\right)  }{k}z^{k}-z\left(  \log z-1+\gamma\right)  -\frac{1}%
{2}\log z-\frac{1}{2}\log(2\pi) \label{Taylor_series_around_z=0}%
\end{equation}
where the Taylor series of $\log\Gamma(z+1)$ around $z=0$,
\[
\log\Gamma(z+1)=-\gamma z+\sum_{k=2}^{\infty}\frac{(-1)^{k}\,\zeta(k)}%
{k}\,z^{k}%
\]
follows directly from Weierstrass' product
\begin{equation}
\frac{1}{\Gamma(z+1)}=e^{\gamma\,z}\prod_{n=1}^{\infty}\left(  1+\frac{z}%
{n}\right)  \,e^{-z/n} \label{Weierstrass_product_Gamma_function}%
\end{equation}
In contrast to Stirling's series $\sum_{k=1}^{K}\frac{\left(  -1\right)
^{k}\zeta\left(  -k\right)  }{k}z^{-k}$ for some finite $K\leq K^{\ast}\left(
z\right)  $ in (\ref{Binet_asymptotic_Stirling_expansion}), violation of the
restriction $\left\vert z\right\vert <1$ in the Taylor series in
(\ref{Taylor_series_around_z=0}) leads to useless results.

\subsection{Convergent companion of Stirling's asymptotic series
(\ref{Binet_asymptotic_Stirling_expansion})}

\label{sec_Stirling_convergent_companion}The Taylor series of the
entire\footnote{An \emph{entire} function has no singularities in the finite
complex plane and possesses a Taylor series around any finite point with
infinitely large radius of convergence. An entire function is sometimes also
called an \emph{integral} function (as e.g. in \cite{Titchmarshfunctions}).}
function $\left(  s-1\right)  \zeta\left(  s\right)  =\sum_{m=0}^{\infty}%
g_{m}\left(  z\right)  \left(  s-z\right)  ^{m}$ around $s_{0}=z$ converges
for all finite complex $z$. After substituting the Taylor series $\left(
s-1\right)  \zeta\left(  s\right)  =\sum_{m=0}^{\infty}g_{m}\left(  1\right)
\left(  s-1\right)  ^{m}$ around $s_{0}=1$ into the complex integral
(\ref{Binet_complex_integral_Zeta}), it is shown in Appendix
\ref{sec_convergent_mu_series_Stirling} that%
\begin{equation}
\mu\left(  z\right)  =\sum_{m=1}^{\infty}g_{m}\left(  1\right)  \sum
_{v=1}^{m-1}(v-1)!(-1)^{m-1-v}\mathcal{S}_{m}^{(v+1)}\,\left(  \frac{1}%
{z}\right)  ^{v} \label{mu_series_TaylorZeta_1_op_z}%
\end{equation}
where $\mathcal{S}_{m}^{(k)}$ is the Stirling number of second Kind. The
Taylor coefficients \cite[(23.2.5)]{Abramowitz} for $k\geq0$
\begin{equation}
g_{m+1}(1)=\frac{(-1)^{m}}{m!}\lim_{K\rightarrow\infty}\left(  \sum_{n=1}%
^{K}\frac{\ln^{m}n}{n}-\frac{\ln^{m+1}K}{m+1}\right)  \label{Taylor_zeta1}%
\end{equation}
are attributed to Stieltjes, with $g_{0}\left(  1\right)  =1$ and
$g_{1}\left(  1\right)  =\gamma=\lim_{K\rightarrow\infty}\left(  \sum
_{n=1}^{K}\frac{1}{n}-\ln K\right)  $. However, computationally, Stieltjes
expression (\ref{Taylor_zeta1}) is less suited and we present fast converging
series for $g_{m}\left(  1\right)  $ in Appendix
\ref{sec_Taylor_coefficient_RZeta_around1}. Since $\left(  s-1\right)
\zeta\left(  s\right)  $ is entire, the Taylor coefficients $g_{m}\left(
1\right)  =O\left(  \frac{1}{m!}\right)  $ -- just as those of any entire
function of order 1 like $e^{z}$ -- decay rapidly in $m$ and only a few terms
in (\ref{mu_series_TaylorZeta_1_op_z}) provide accurate results for Binet's
function $\mu\left(  z\right)  $.

Although the reversal of the $m$- and $v$- sum in
(\ref{mu_series_TaylorZeta_1_op_z}) is not allowed, it is interesting to
illustrate what happens if we reverse the sums%
\begin{equation}
\mu\left(  z\right)  \overset{!}{=}\sum_{v=1}^{\infty}(v-1)!\left\{
\sum_{m=v}^{\infty}g_{m+1}\left(  1\right)  (-1)^{m-v}\mathcal{S}%
_{m+1}^{(v+1)}\right\}  \,\frac{1}{z^{v}} \label{mu_reversed_wrong_series}%
\end{equation}
We substitute the closed form (\ref{stirling2closed}) of the Stirling numbers,
$\mathcal{S}_{k}^{(m)}=\frac{1}{m!}\sum_{j=0}^{m}(-1)^{m-j}{\binom{m}{j}}%
j^{k}$, using $\mathcal{S}_{k}^{(m)}=0$ if $k<m$,
\begin{align*}
\sum_{m=v}^{\infty}g_{m+1}\left(  1\right)  (-1)^{m-v}\mathcal{S}%
_{m+1}^{(v+1)}  &  =\sum_{m=-1}^{\infty}g_{m+1}\left(  1\right)
(-1)^{m-v}\mathcal{S}_{m+1}^{(v+1)}\\
&  =\frac{1}{\left(  v+1\right)  !}\sum_{m=-1}^{\infty}g_{m+1}\left(
1\right)  (-1)^{m-v}\sum_{j=0}^{v+1}(-1)^{v+1-j}{\binom{v+1}{j}}j^{m+1}\\
&  =\frac{1}{\left(  v+1\right)  !}\sum_{j=0}^{v+1}(-1)^{j}{\binom{v+1}{j}%
}\sum_{m=-1}^{\infty}g_{m+1}\left(  1\right)  (-1)^{m+1}j^{m+1}%
\end{align*}
With $\left(  s-1\right)  \zeta\left(  s\right)  =\sum_{m=0}^{\infty}%
g_{m}\left(  1\right)  \left(  s-1\right)  ^{m}$, we have%
\[
\sum_{m=v}^{\infty}g_{m+1}\left(  1\right)  (-1)^{m-v}\mathcal{S}%
_{m+1}^{(v+1)}=\frac{1}{\left(  v+1\right)  !}\sum_{j=0}^{v+1}(-1)^{j-1}%
{\binom{v+1}{j}}j\zeta\left(  1-j\right)
\]
After substitution in the \textquotedblleft erroneous\textquotedblright%
\ series (\ref{mu_reversed_wrong_series}) for $\mu\left(  z\right)  $,
\[
\mu\left(  z\right)  \overset{!}{=}\sum_{v=1}^{\infty}\frac{(v-1)!}{\left(
v+1\right)  !}\sum_{j=0}^{v+1}(-1)^{j-1}{\binom{v+1}{j}}j\zeta\left(
1-j\right)  \,\frac{1}{z^{v}}%
\]
and using $\zeta(-k)=\frac{(-1)^{k}}{k+1}B_{k+1}$ and $\sum_{j=0}^{n}%
\binom{n+1}{j}B_{j}=0$, we obtain%
\[
\mu\left(  z\right)  \overset{!}{=}\sum_{v=1}^{\infty}\frac{B_{v+1}}{\left(
v+1\right)  v}\,\frac{1}{z^{v}}%
\]
Hence, reversal of the $m$- and $v$- sum in (\ref{mu_series_TaylorZeta_1_op_z}%
) again leads to Stirling's \emph{diverging} asymptotic series
(\ref{Binet_asymptotic_Stirling_expansion}). The series
(\ref{mu_series_TaylorZeta_1_op_z}) converges for all complex $z$ with
$\left\vert \arg z\right\vert <\pi$ and can be regarded as the convergent
companion of Stirling's asymptotic series
(\ref{Binet_asymptotic_Stirling_expansion}).

\section{Binet's investigations}

\subsection{Binet's first factorial series for $\mu\left(  z\right)  $}

\label{sec_Binet_first_factorial_expansion}We review Binet's first expansion
for $\mu\left(  z\right)  $ in \cite[Section 3, pp. 223-229]{Binet1839}.
Writing $\log\frac{z}{z+1}=\log\left(  1-\frac{1}{z+1}\right)  $, Binet
expands the right-hand side of the forward difference formula
(\ref{forwards_difference_mu})%
\[
\mu\left(  z+1\right)  -\mu\left(  z\right)  =1+\left(  z+\frac{1}{2}\right)
\log\left(  1-\frac{1}{z+1}\right)
\]
by introducing the Taylor series around $z_{0}=0$ of $\log\left(  1-z\right)
=-\sum_{n=1}^{\infty}\frac{z^{n}}{n}$, convergent for $\left\vert z\right\vert
<1$, and obtains for $\left\vert z+1\right\vert >1$,%
\begin{equation}
\mu\left(  z\right)  -\mu\left(  z+1\right)  =\frac{1}{2}\sum_{n=1}^{\infty
}\frac{n}{\left(  n+2\right)  \left(  n+1\right)  }\frac{1}{\left(
z+1\right)  ^{n+1}} \label{forwards_difference_mu_Taylor_series_1/(z+1)}%
\end{equation}
Binet replaces $z\rightarrow z+k$ in
(\ref{forwards_difference_mu_Taylor_series_1/(z+1)}), sums over all integer
$k\geq0$,%
\[
\sum_{k=0}^{N-1}\mu\left(  z+k\right)  -\mu\left(  z+k+1\right)  =\frac{1}%
{2}\sum_{n=1}^{\infty}\frac{n}{\left(  n+2\right)  \left(  n+1\right)  }%
\sum_{k=0}^{N-1}\frac{1}{\left(  z+k+1\right)  ^{n+1}}%
\]
and rewrites the telescoping series at the left-hand side%
\[
\mu\left(  z\right)  -\mu\left(  z+N\right)  =\frac{1}{2}\sum_{n=1}^{\infty
}\frac{n}{\left(  n+2\right)  \left(  n+1\right)  }\sum_{k=0}^{N-1}\frac
{1}{\left(  z+k+1\right)  ^{n+1}}%
\]
After observing that $\lim_{N\rightarrow\infty}\mu\left(  z+N\right)  =0$
(which follows e.g. from the bound (\ref{bound_mu_Re(z)>0})), Binet \cite[eq.
(58), p. 229]{Binet1839} arrives at his first convergent expansion%
\begin{equation}
\mu\left(  z\right)  =\frac{1}{2}\sum_{n=1}^{\infty}\frac{n}{\left(
n+2\right)  \left(  n+1\right)  }\sum_{k=0}^{\infty}\frac{1}{\left(
z+k+1\right)  ^{n+1}} \label{mu_convergent_expansion_polygamma_type}%
\end{equation}
Analogously, substitution of the Taylor series $\log\left(  1+z\right)
=\sum_{n=1}^{\infty}\frac{(-1)^{n-1}z^{n}}{n}$ for $\left\vert z\right\vert
<1$ in Gudermann's series (\ref{Gudermann_series_mu}) yields, after
reworking\footnote{Adding (\ref{mu_convergent_expansion_polygamma_type}) and
(\ref{mu_convergent_expansion_Gudermann}) still yields a slowly convergent
series%
\[
\mu\left(  z\right)  =\frac{1}{4}\left(  2z+1\right)  \log\left(  1+\frac
{1}{z}\right)  -\frac{1}{2}+\frac{1}{4}\sum_{n=1}^{\infty}\frac{\left(
2n-1\right)  }{\left(  2n+1\right)  n}\sum_{k=1}^{\infty}\frac{1}{\left(
z+k\right)  ^{2n}}%
\]
},
\begin{equation}
\mu\left(  z\right)  =\frac{1}{2}\sum_{n=1}^{\infty}\frac{n}{\left(
n+1\right)  \left(  n+2\right)  }\sum_{k=0}^{\infty}\frac{\left(  -1\right)
^{n+1}}{\left(  z+k\right)  ^{n+1}} \label{mu_convergent_expansion_Gudermann}%
\end{equation}
The polygamma functions $\psi^{(n)}(z)=\frac{d^{n}}{dz^{n}}\log\Gamma\left(
z\right)  $, for any integer $n\geq1$, possess the convergent series
\cite[6.4.10]{Abramowitz}
\begin{equation}
\psi^{(n)}(z)=(-1)^{n+1}n!\sum_{k=0}^{\infty}\frac{1}{(z+k)^{n+1}}
\label{polygamma_series}%
\end{equation}
In terms of the polygamma functions (\ref{polygamma_series}), we rewrite the
two variants (\ref{mu_convergent_expansion_polygamma_type}) and
(\ref{mu_convergent_expansion_Gudermann}) as
\begin{align*}
\mu\left(  z\right)   &  =\frac{1}{2}\sum_{n=1}^{\infty}\frac{n}{\left(
n+2\right)  !}\psi^{(n)}(z+1)\left(  -1\right)  ^{n+1}\\
\mu\left(  z\right)   &  =\frac{1}{2}\sum_{n=1}^{\infty}\frac{n}{\left(
n+2\right)  !}\psi^{(n)}(z)
\end{align*}

Binet \cite[{art [20], p. 232-234}]{Binet1839} then
concentrates\footnote{\label{footnote_Stirling_factorial_series}Binet invokes
the factorial expansion $\frac{1}{z-b}=\frac{\Gamma\left(  z+\alpha\right)
}{\Gamma\left(  b+\alpha\right)  }\sum_{m=0}^{\infty}\frac{\Gamma\left(
b+\alpha+m\right)  }{\Gamma\left(  z+\alpha+m+1\right)  }$, derived in
(\ref{1_op_(z-b)_factorial_series}) below, with $p=z+\alpha$ and $a=b+\alpha
$,
\[
\frac{1}{p-a}=\frac{1}{p}+\frac{a}{p\left(  p+1\right)  }+\frac{a(a+1)}%
{p\left(  p+1\right)  \left(  p+2\right)  }+\cdots=\frac{\Gamma\left(
p\right)  }{\Gamma\left(  a\right)  }\sum_{m=0}^{\infty}\frac{\Gamma\left(
a+m\right)  }{\Gamma\left(  p+m+1\right)  }%
\]
and Newton's difference expansion $f\left(  p+a\right)  =\sum_{k=0}^{\infty
}\Delta^{k}f\left(  p\right)  \binom{a}{k}$. Our factorial expansion in
(\ref{1/(z+a)_finite_factorial_series}) generalizes
(\ref{1_op_(z-b)_factorial_series}).} on the evaluation of the $k$-sum in
(\ref{mu_convergent_expansion_polygamma_type}), thus on the higher-order
derivatives $\psi^{(n)}(z)$ of the digamma function $\psi\left(  z\right)  $
and presents \cite[p. 234]{Binet1839} his first factorial expansion
(\ref{Binet_original}). Binet \cite[{art [21], p. 242}]{Binet1839} proceeds by
constructing integrals for $\mu\left(  z\right)  $. Introducing Euler's Gamma
integral $\frac{\Gamma\left(  s\right)  }{x^{s}}=\int_{0}^{\infty}%
t^{s-1}e^{-xt}dt$, valid for $\operatorname{Re}\left(  s\right)  >0$ and
$\operatorname{Re}\left(  x\right)  >0$, into
(\ref{mu_convergent_expansion_polygamma_type}) yields, for $\operatorname{Re}%
\left(  z\right)  >0$,%
\[
\mu\left(  z\right)  =\frac{1}{2}\sum_{n=1}^{\infty}\frac{n}{\left(
n+2\right)  \left(  n+1\right)  }\frac{1}{n!}\sum_{k=0}^{\infty}\int
_{0}^{\infty}t^{n}e^{-\left(  z+k+1\right)  t}dt
\]
After reversal of integral and $k$-summation,%
\[
\mu\left(  z\right)  =\frac{1}{2}\sum_{n=1}^{\infty}\frac{n}{\left(
n+2\right)  \left(  n+1\right)  }\frac{1}{n!}\int_{0}^{\infty}e^{-zt}%
\frac{t^{n}}{e^{t}-1}dt
\]
and using $\frac{n}{\left(  n+2\right)  \left(  n+1\right)  }=\frac{2}%
{n+2}-\frac{1}{n+1}$ in the $n$-sum, which results in $\frac{1}{2}\sum
_{n=1}^{\infty}\frac{n}{\left(  n+2\right)  \left(  n+1\right)  }\frac{t^{n}%
}{n!}=\frac{e^{t}t-2e^{t}+2+t}{t^{2}}$, then leads to Binet's integral
(\ref{Binet_mu_integral_Bernoulli_generating_function}). Via an integral due
to Poisson, $\frac{e^{t}+1}{e^{t}-1}-\frac{1}{2t}=\frac{1}{4}\int_{0}^{\infty
}\frac{\sin tx\,dx}{e^{2\pi x}-1}$, Binet also derives
(\ref{Binet_mu_integral_arctan}). Binet writes at length and reconsiders
previous derivations, but his great Memoire definitely contains the
foundations about his function $\mu\left(  z\right)  $.

\subsection{Binet's second factorial expansion}

\label{sec_Binet_expansion}

\begin{theorem}
\label{theorem_Binet_convergent_expansion}A second convergent factorial series
of Binet's function $\mu\left(  z\right)  $ is
\begin{equation}
\mu\left(  z\right)  =\frac{1}{z}\sum_{m=1}^{\infty}\frac{c_{m}}{\prod
_{k=1}^{m-1}(z+k)}\hspace{1cm}\text{for }\operatorname{Re}\left(  z\right)  >0
\label{Binet_function_faculteits_expansie}%
\end{equation}
where the rational coefficients are%
\begin{equation}
c_{m}=\frac{\left(  -1\right)  ^{m-1}}{2m}\sum_{k=1}^{m}\frac{kS_{m}^{(k)}%
}{\left(  k+2\right)  \left(  k+1\right)  } \label{def_Binet_coefficients}%
\end{equation}
and $S_{m}^{(k)}$ is the Stirling Number of the First Kind.
\end{theorem}

We essentially follow the steps in Binet's original proof in \cite[p.
339]{Binet1839}. In Section \ref{sec_1_factorial_series}, we formalize Binet's
proof as a recipe in five steps.

\textbf{Proof (Binet):} Binet \cite[p. 339]{Binet1839} substitutes
$e^{-t}=1-u$ or $t=-\log\left(  1-u\right)  $ in the integral
(\ref{Binet_mu_integral_Bernoulli_generating_function}),%
\begin{equation}
2\mu\left(  z\right)  =-\int_{0}^{1}\frac{\left(  1-u\right)  ^{z-1}}%
{u}\left(  \frac{2-u}{\log\left(  1-u\right)  }+\frac{2u}{\log^{2}\left(
1-u\right)  }\right)  du \label{Binet_function_integral}%
\end{equation}
and proceeds to expand
\[
\frac{2-u}{\log\left(  1-u\right)  }+\frac{2u}{\log^{2}\left(  1-u\right)
}=2\left\{  \frac{1}{\log\left(  1-u\right)  }+\frac{u}{\log^{2}\left(
1-u\right)  }\right\}  -\frac{u}{\log\left(  1-u\right)  }%
\]
in a Taylor series around $u=0$. Instead of following Binet, who has used
integrals rather than Stirling numbers $S_{m}^{(k)}$, we invoke the Taylor
expansion (\ref{Taylor_series_x/Log^n(1+x)}) for $n=2$, derived in the
Appendix \ref{sec_Taylor_series_x/log(1+x)} and convergent for $\left\vert
u\right\vert <1$,
\[
\frac{1}{\log\left(  1-u\right)  }+\frac{u}{\log^{2}\left(  1-u\right)
}=-\frac{1}{2}-\sum_{m=1}^{\infty}\left(  \sum_{k=1}^{m}\frac{k!S_{m}^{(k)}%
}{\left(  k+2\right)  !}\right)  \frac{\left(  -u\right)  ^{m}}{m!}%
\]
and the Taylor series (\ref{Taylor_series_x/Log(1+x)})%
\[
\frac{-u}{\log\left(  1-u\right)  }=1+\sum_{m=1}^{\infty}\left(  \sum
_{k=1}^{m}\frac{S_{m}^{(k)}}{k+1}\right)  \frac{\left(  -u\right)  ^{m}}{m!}%
\]
to obtain the Taylor series, valid for $\left\vert u\right\vert <1$,%
\begin{equation}
\frac{2-u}{\log\left(  1-u\right)  }+\frac{2u}{\log^{2}\left(  1-u\right)
}=\sum_{m=1}^{\infty}\left(  \sum_{k=1}^{m}\frac{kS_{m}^{(k)}}{\left(
k+2\right)  \left(  k+1\right)  }\right)  \frac{\left(  -u\right)  ^{m}}{m!}
\label{Taylor_series_bijna_cm_coefficienten}%
\end{equation}
Introducing (\ref{Taylor_series_bijna_cm_coefficienten}) in Binet's function
(\ref{Binet_function_integral}) and reversing summation and integral,
justified because a Taylor series can be term-wise integrated within its
radius of convergence,%
\[
2\mu\left(  z\right)  =\sum_{m=1}^{\infty}\left(  \sum_{k=1}^{m}\frac
{kS_{m}^{(k)}}{\left(  k+2\right)  \left(  k+1\right)  }\right)  \frac{\left(
-1\right)  ^{m-1}}{m!}\int_{0}^{1}\left(  1-u\right)  ^{z-1}u^{m-1}du
\]
using the Beta integral $\int_{0}^{1}u^{p-1}\left(  1-u\right)  ^{q-1}%
du=\frac{\Gamma\left(  p\right)  \Gamma\left(  q\right)  }{\Gamma\left(
p+q\right)  }$, valid for $\operatorname{Re}\left(  p\right)  >0$ and
$\operatorname{Re}\left(  q\right)  >0$, yields a converging series, for
$\operatorname{Re}\left(  z\right)  >0$,%
\[
\mu\left(  z\right)  =\frac{1}{2}\sum_{m=1}^{\infty}\left(  \sum_{k=1}%
^{m}\frac{kS_{m}^{(k)}}{\left(  k+2\right)  \left(  k+1\right)  }\right)
\frac{\left(  -1\right)  ^{m-1}}{m}\frac{\Gamma\left(  z\right)  }%
{\Gamma\left(  z+m\right)  }%
\]
With $\frac{\Gamma\left(  z+m\right)  }{\Gamma\left(  z\right)  }=\prod
_{k=0}^{m-1}(z+k)$, we arrive at Binet's second factorial series
(\ref{Binet_function_faculteits_expansie}).\hfill$\square\medskip$

The first few coefficients $c_{m}$ in (\ref{def_Binet_coefficients}) are
$c_{1}=\frac{1}{12}$, $c_{2}=0$, $c_{3}=-\frac{1}{360}$, $c_{4}=-\frac{1}%
{120}$, $c_{5}=-\frac{5}{168}$, $c_{6}=-\frac{11}{84}$, $c_{7}=-\frac
{3499}{5040}$, which are smaller in absolute value than 1, but $c_{8}%
=-\frac{1039}{240}$, $c_{9}=-\frac{369689}{11880}$ exceed 1 in absolute value.
It holds that $\left\vert c_{m}\right\vert >1$ for $m>8$ as shown below after
Theorem \ref{theorem_Binet_coefficient_integral}.

The generating function of the Stirling numbers $S_{m}^{(k)}$ of the First
Kind \cite[Sec. 24.1.3 and 24.1.4]{Abramowitz},
\begin{equation}
m!\binom{x}{m}=\frac{\Gamma(x+1)}{\Gamma(x+1-m)}=\prod_{k=0}^{m-1}%
(x-k)=\sum_{k=0}^{m}S_{m}^{(k)}\;x^{k} \label{genfunc_stirling}%
\end{equation}
indicates that $\sum_{k=1}^{m}S_{m}^{(k)}=\delta_{1m}$. Thus, we can add
$a\sum_{k=1}^{m}S_{m}^{(k)}$ to $\sum_{k=1}^{m}\frac{kS_{m}^{(k)}}{\left(
k+2\right)  \left(  k+1\right)  }$ in (\ref{def_Binet_coefficients}) and find
that the coefficient equals%
\[
c_{m}=\frac{\left(  -1\right)  ^{m-1}}{2m}\sum_{k=1}^{m}\frac{ak^{2}+\left(
3a-1\right)  k+2a}{\left(  k+1\right)  \left(  k+2\right)  }S_{m}^{(k)}%
\hspace{1cm}\text{for }m>1\text{ and any }a\in\mathbb{C}%
\]
For example, for $a=\frac{1}{6}$ and $m>1$, we find
\[
c_{m}=\frac{\left(  -1\right)  ^{m-1}}{12m}\sum_{k=1}^{m}\frac{(k-1)\left(
k-2\right)  }{\left(  k+1\right)  \left(  k+2\right)  }S_{m}^{(k)}%
\]
illustrating that $c_{2}=0$. The second generating function of the Stirling
numbers $S_{m}^{\left(  k\right)  }$, convergent for $\left\vert u\right\vert
<1$, is (see e.g \cite[24.1.3.A]{Abramowitz})%
\begin{equation}
\log^{k}(1+u)=k!\sum_{m=k}^{\infty}S_{m}^{(k)}\;\frac{u^{m}}{m!}
\label{genfunc_stirlinginf}%
\end{equation}

\begin{theorem}
\label{theorem_Binet_coefficient_integral}The rational coefficients $c_{m}$ in
the second factorial series (\ref{Binet_function_faculteits_expansie}) of
Binet's function $\mu\left(  z\right)  $ can be represented by an integral
\begin{equation}
c_{m}=\frac{1}{m}\int_{0}^{1}\left(  u-\frac{1}{2}\right)  u\prod_{k=1}%
^{m-1}(k-u)du \label{Binet_coefficient_integral}%
\end{equation}
Moreover, all coefficients $c_{m}$, except for the first two, are negative,
i.e. $c_{m}<0$ for all $m>2$.
\end{theorem}

\textbf{Proof}: Using $\frac{j}{\left(  j+2\right)  \left(  j+1\right)
}=\frac{2}{j+2}-\frac{1}{j+1}$, the coefficient $c_{m}$ in
(\ref{def_Binet_coefficients}) is written as
\[
c_{m}=\frac{\left(  -1\right)  ^{m}}{2m}\sum_{j=1}^{m}\frac{S_{m}^{(j)}}%
{j+1}-\frac{\left(  -1\right)  ^{m}}{2m}\sum_{j=1}^{m}\frac{2S_{m}^{(j)}}{j+2}%
\]
Multiplying both sides of the generating function (\ref{genfunc_stirling}) of
the Stirling numbers $S_{m}^{(k)}$ by $x^{q-1}$ and integrating yields
\begin{equation}
\int_{a}^{b}u^{q-1}\prod_{k=0}^{m-1}(u-k)du=\sum_{l=0}^{m}S_{m}^{(l)}%
\;\frac{b^{l+q}-a^{l+q}}{l+q} \label{integrated_Stirling_gf}%
\end{equation}
After substituting the case for $q=1$ and $q=2$, we obtain%
\[
c_{m}=\frac{\left(  -1\right)  ^{m}}{2m}\int_{0}^{1}\left(  1-2u\right)
\prod_{k=0}^{m-1}(u-k)du
\]
from which (\ref{Binet_coefficient_integral}) follows. For $m=1$ in
(\ref{Binet_coefficient_integral}), we find $c_{1}=\int_{0}^{1}\left(
u-\frac{1}{2}\right)  udu=\frac{1}{12}$.

In the second part, we will demonstrate that $c_{m}<0$ for $m>2$. Since
$u\prod_{k=1}^{m-1}(k-u)\geq0$ for $u\in\lbrack0,1]$, we split the integration
interval in (\ref{Binet_coefficient_integral}),
\[
mc_{m}=\int_{0}^{\frac{1}{2}}\left(  u-\frac{1}{2}\right)  u\prod_{k=1}%
^{m-1}(k-u)du+\int_{\frac{1}{2}}^{1}\left(  u-\frac{1}{2}\right)  u\prod
_{k=1}^{m-1}(k-u)du
\]
After making the substitution $u=1-w$ in the last integral, we arrive at%
\[
mc_{m}=\int_{0}^{\frac{1}{2}}\left(  \frac{1}{2}-u\right)  \left(  1-u\right)
u\left\{  \prod_{k=2}^{m-1}(k-\left(  1-u\right)  )-\prod_{k=2}^{m-1}%
(k-u)\right\}  du
\]
For $m=2$, the right-hand side is zero, because the two products are equal to
1. Since $1-u>u$ for $0\leq u<\frac{1}{2}$, the product $\prod_{k=2}%
^{m-1}(k-\left(  1-u\right)  )<\prod_{k=2}^{m-1}(k-u)$ for $u\in\lbrack
0,\frac{1}{2})$ and for $m>2$. Hence, we conclude that $c_{m}<0$ for
$m>2$.\hfill$\square\medskip$

The important fact in Theorem \ref{theorem_Binet_coefficient_integral}, that
all coefficients $c_{m}<0$ for $m>2$, implies that any truncation at $m=K>2$
terms in (\ref{Binet_function_faculteits_expansie}) upper bounds Binet's
function $\mu\left(  z\right)  $. In appendix
\ref{sec_integral_Binet_coefficient}, we derive an other integral
(\ref{Integral_Binet_coefficient_c[m]}) for the coefficients $c_{m}$.

\subsection{Growth of the coefficient $c_{m}$ with $m$}

\label{sec_growth_cm_with_m}Since $\prod_{k=1}^{m}(k-u)=\left(  m-u\right)
\prod_{k=1}^{m-1}(k-u)$ for $m\geq1$, the integral
(\ref{Binet_coefficient_integral}) becomes
\begin{align*}
\left(  m+1\right)  c_{m+1}  &  =\int_{0}^{1}\left(  u-\frac{1}{2}\right)
u\left(  m-u\right)  \prod_{k=1}^{m-1}(k-u)du\\
&  =m\int_{0}^{1}\left(  u-\frac{1}{2}\right)  u\prod_{k=1}^{m-1}%
(k-u)du-\int_{0}^{1}\left(  u-\frac{1}{2}\right)  u^{2}\prod_{k=1}%
^{m-1}(k-u)du
\end{align*}
The last integral is smaller in absolute value than $m\left\vert
c_{m}\right\vert $, because $u^{2}\leq u$ for $u\in\left[  0,1\right]  $.
However, unlike the proof of Theorem \ref{theorem_Binet_coefficient_integral},
the last integral is positive. Indeed,%
\[
\int_{0}^{1}\left(  u-\frac{1}{2}\right)  u^{2}\prod_{k=1}^{m-1}%
(k-u)du=\int_{0}^{\frac{1}{2}}\left(  \frac{1}{2}-u\right)  u\left(
1-u\right)  \left\{  \left(  1-u\right)  \prod_{k=2}^{m-1}(k-\left(
1-u\right)  )-u\prod_{k=2}^{m-1}(k-u)\right\}  dw
\]
and $\left(  1-u\right)  \prod_{k=2}^{m-1}(k-\left(  1-u\right)
)>u\prod_{k=2}^{m-1}\left(  k-u\right)  $ for $0\leq u<\frac{1}{2}$. Thus, we
find the inequality%
\[
\left(  m+1\right)  c_{m+1}\geq m^{2}c_{m}-m\left\vert c_{m}\right\vert
=m\left(  m+1\right)  c_{m}%
\]
Iterating this recursion inequality, $c_{m}\geq\left(  m-1\right)  c_{m-1}$,
yields%
\[
c_{m}\geq\left(  m-1\right)  c_{m-1}\geq\cdots\geq\left(  m-1\right)
\cdots\left(  m-p\right)  c_{m-p}%
\]
With $c_{3}=-\frac{1}{360}$ and $p=m-3$, we obtain $c_{m}\geq-\frac{\left(
m-1\right)  !}{720}$. The recursion inequality demonstrates that $\left\vert
c_{m}\right\vert $ increases strictly with $m$ for $m\geq2$.

The logarithmic behavior (\ref{mu_around_z=0}) of $\mu\left(  z\right)  $
around $z=0$ shows that $\lim_{z\rightarrow0}z\mu\left(  z\right)  =0$.
Binet's second factorial series (\ref{Binet_function_faculteits_expansie}),
written as%
\[
z\mu\left(  z\right)  =\sum_{m=1}^{\infty}\frac{c_{m}}{\prod_{k=1}^{m-1}%
(z+k)}=\frac{1}{12}-\sum_{m=2}^{\infty}\frac{\left\vert c_{m+1}\right\vert
}{\prod_{k=1}^{m}(z+k)}%
\]
illustrates, with $\lim_{z\rightarrow0}z\mu\left(  z\right)  =0$ that%
\[
0=\sum_{m=1}^{\infty}\frac{c_{m}}{\left(  m-1\right)  !}=\frac{1}{12}%
-\sum_{m=3}^{\infty}\frac{\left\vert c_{m}\right\vert }{\left(  m-1\right)  !}%
\]
where $\sum_{m=1}^{K}\frac{c_{m}}{\left(  m-1\right)  !}$ converges very
slowly with increasing $K$. For positive real $z$, it holds that $\sum
_{m=3}^{\infty}\frac{\left\vert c_{m}\right\vert }{\prod_{k=1}^{m-1}(z+k)}%
\leq\sum_{m=3}^{\infty}\frac{\left\vert c_{m}\right\vert }{\left(  m-1\right)
!}=\frac{1}{12}$, which agrees with the bound (\ref{bound_mu_Re(z)>0}). Since
all $c_{m}<0$ for $2<m$ by Theorem \ref{theorem_Binet_coefficient_integral},
the convergence indicates that $\frac{c_{m}}{\left(  m-1\right)  !}=O\left(
\frac{1}{m^{1+\varepsilon}}\right)  $ for $\varepsilon>0$. Alternatively, with
$\prod_{k=1}^{m-1}(k-u)=\left(  m-1\right)  !\prod_{k=1}^{m-1}(1-\frac{u}{k}%
)$, the integral (\ref{Binet_coefficient_integral}) is%
\[
\frac{mc_{m}}{\left(  m-1\right)  !}=\int_{0}^{1}\left(  u-\frac{1}{2}\right)
u\prod_{k=1}^{m-1}\left(  1-\frac{u}{k}\right)  du
\]
Since all factors in the last integrand are in absolute value smaller than or
equal to 1, the integral decreases in absolute value with $m$ and we conclude
that $\frac{\left\vert c_{m}\right\vert }{\left(  m-1\right)  !}<\frac{1}{m}$,
which is a prerequisite for convergence of $\sum_{m=1}^{\infty}\frac{c_{m}%
}{\left(  m-1\right)  !}$. The asymptotic behavior of $c_{m}$ for large $m$ is
computed in the Appendix \ref{sec_asymptotic_bm(alpha)_large_m}.

\section{Gilbert's investigations}

\label{sec_Gilbert}

\subsection{Gilbert's expansion (\ref{mu_Gilbert_Gudermann}) of $\mu\left(
z\right)  $}

By substitution in Binet's integral $\mu\left(  z\right)  =\int_{0}^{\infty
}\frac{e^{-zt}}{t}\left(  \frac{1}{e^{t}-1}-\frac{1}{t}+\frac{1}{2}\right)
dt$ in (\ref{Binet_mu_integral_Bernoulli_generating_function_2}) of the
partial fraction expansion\footnote{Cauchy's integral $\frac{1}{e^{t}-1}%
=\frac{1}{2\pi i}\int_{C(t)}\frac{dw}{(w-t)\;(e^{w}-1)}$, where the contour
$C\left(  t\right)  $ encloses only the point $w=t$, leads, after deforming
the contour to enclose the entire plane except for a small region around
$w=t$, to (\ref{partialfract_inv(exp_t-1)}).}
\begin{equation}
\frac{1\,}{e^{t}-1}=\frac{1}{t}-\frac{1}{2}+2t\sum_{k=1}^{\infty}\frac
{1}{4k^{2}\pi^{2}+t^{2}} \label{partialfract_inv(exp_t-1)}%
\end{equation}
Gilbert \cite[art. 6]{Gilbert1873} obtains%
\begin{equation}
\mu\left(  z\right)  =2\sum_{k=1}^{\infty}\int_{0}^{\infty}\frac{e^{-zt}%
}{4k^{2}\pi^{2}+t^{2}}dt \label{mu_Cauchy}%
\end{equation}
which was due to Cauchy. Using the Laplace transform $\int_{0}^{\infty}%
e^{-zu}\sin\left(  au\right)  du=\frac{a}{a^{2}+z^{2}}$ for $\operatorname{Re}%
\left(  z\right)  >0$, Gilbert \cite[art. 6]{Gilbert1873} reformulates the
integral%
\begin{align*}
\int_{0}^{\infty}\frac{e^{-zt}}{a^{2}+t^{2}}dt  &  =\frac{1}{a}\int
_{0}^{\infty}dt\;e^{-zt}\int_{0}^{\infty}e^{-tu}\sin\left(  au\right)
du=\frac{1}{a}\int_{0}^{\infty}\left(  \int_{0}^{\infty}dt\;e^{-\left(
z+u\right)  t}\right)  \sin\left(  au\right)  du\\
&  =\frac{1}{a}\int_{0}^{\infty}\frac{\sin\left(  au\right)  }{z+u}du
\end{align*}
and obtains%
\begin{equation}
\mu\left(  z\right)  =\frac{1}{\pi}\sum_{k=1}^{\infty}\frac{1}{k}\int
_{0}^{\infty}\frac{\sin\left(  2\pi ku\right)  }{z+u}du \label{mu_Gilbert}%
\end{equation}
which can be written as%
\[
\mu\left(  z\right)  =\frac{1}{\pi}\int_{0}^{\infty}\frac{du}{z+u}\sum
_{k=1}^{\infty}\frac{\sin\left(  2\pi ku\right)  }{k}%
\]
Gilbert \cite[art. 8]{Gilbert1873} then introduces the Fourier series
$\sum_{k=1}^{\infty}\frac{\sin\left(  ku\right)  }{k}=\frac{\pi-u}{2}$ for
$u\in\left[  0,2\pi\right]  $ and finds, after some manipulations, that%
\begin{equation}
\mu\left(  z\right)  =\sum_{k=0}^{\infty}\int_{0}^{1}\frac{\left(  \frac{1}%
{2}-x\right)  dx}{z+k+x} \label{mu_Gilbert_Gudermann}%
\end{equation}
Next, Gilbert \cite[art. 9]{Gilbert1873} shows that
(\ref{mu_Gilbert_Gudermann}) can be obtained from Gudermann's series
(\ref{Gudermann_series_mu}), while (\ref{mu_Gilbert_Gudermann}) leads to
Binet's series (\ref{mu_convergent_expansion_polygamma_type}) by observing
that $\int_{0}^{1}\frac{\left(  \frac{1}{2}-x\right)  dx}{z+k+x}=\int_{0}%
^{1}\frac{\left(  u-\frac{1}{2}\right)  du}{z+k+1-u}$, expanding the
denominator into a geometric series and using the integral $\int_{0}%
^{1}\left(  z-\frac{1}{2}\right)  z^{n}dz=\frac{n}{2\left(  n+1\right)
\left(  n+2\right)  }$.

Gilbert \cite[art. 14-16]{Gilbert1873} derives from (\ref{mu_Gilbert}) the
Stirling series (\ref{Binet_asymptotic_Stirling_expansion}) and the
Malmsten-Kummer series (\ref{Malmsten_Kummer_series_LogGama(x)}). Analogous to
the theory of the Riemann Zeta-function \cite{Titchmarshzeta}, Gilbert
\cite[art. 17]{Gilbert1873} integrates the argument $f\left(  t\right)
=\frac{e^{-zt}}{t}\left(  \frac{1}{e^{t}-1}-\frac{1}{t}+\frac{1}{2}\right)  $
in Binet's integral (\ref{Binet_mu_integral_Bernoulli_generating_function_2})
along a contour that starts at the origin the complex $t$-plane along the real
$t$-axis, until $t=R$, travels along a circle $t=Re^{i\theta}$ to the
imaginary $t$-axis, from which the contour returns to the origin by passing
the poles of $\frac{1}{e^{t}-1}$ at $t=2mi\pi$ along a small semicircles in
the positive $\operatorname{Re}\left(  t\right)  $-plane. After letting
$R\rightarrow\infty$ and invoking Cauchy's integral theorem, Gilbert obtains,
after some manipulations,%
\[
\mu\left(  z\right)  =\frac{1}{2}\sum_{k=1}^{\infty}\frac{\cos\left(  2\pi
kz\right)  }{k}+\frac{1}{2}\int_{0}^{\infty}\left(  \frac{1}{x}-\cot x\right)
\frac{\sin\left(  2zx\right)  }{x}dx
\]
as well as $\sum_{k=1}^{\infty}\frac{\sin\left(  2\pi kz\right)  }{k}=\int
_{0}^{\infty}\left(  \frac{1}{x}-\cot x\right)  \frac{\cos\left(  2zx\right)
}{x}dx$, from which he again deduces the Malmsten-Kummer series
(\ref{Malmsten_Kummer_series_LogGama(x)}) in \cite[art. 30]{Gilbert1873}.

\subsection{Gilbert's generalized factorial series for Binet's function}

\label{sec_Gilbert_generalized_factorial_series}Gilbert started from the
factorial series (see footnote \ref{footnote_Stirling_factorial_series}) due
to Stirling. We slightly generalize Gilbert's derivations in
\cite{Gilbert1886} by starting from our general factorial series of $\frac
{1}{z+a}$ in (\ref{1/(z+a)_finite_factorial_series}).

In the identity $\frac{1}{z+a}=\frac{1}{z+b_{k}}+\frac{b_{k}-a}{\left(
z+b_{k}\right)  }\frac{1}{\left(  z+a\right)  }$ for arbitrary numbers $b_{k}%
$, we recursively replace $\frac{1}{z+a}$ in each iteration $k$
\begin{align*}
\frac{1}{z+a}  &  =\frac{1}{z+b_{1}}+\frac{b_{1}-a}{\left(  z+b_{1}\right)
\left(  z+a\right)  }\\
&  =\frac{1}{z+b_{1}}+\frac{b_{1}-a}{\left(  z+b_{1}\right)  }\frac{1}%
{z+b_{2}}+\frac{b_{1}-a}{\left(  z+b_{1}\right)  }\frac{b_{2}-a}{\left(
z+b_{2}\right)  }\frac{1}{\left(  z+a\right)  }%
\end{align*}
and obtain, after $p$ iterations, the finite factorial series%
\begin{equation}
\frac{1}{z+a}=\sum_{l=1}^{p}\frac{1}{z+b_{l}}\prod_{j=1}^{l-1}\frac{b_{j}%
-a}{b_{j}+z}+\frac{1}{z+a}\prod_{j=1}^{p}\frac{b_{j}-a}{b_{j}+z}
\label{1/(z+a)_finite_factorial_series}%
\end{equation}
If $a$, $z$ and the set $\left\{  b_{j}\right\}  _{1\leq j\leq p}$ are
positive real numbers, then (\ref{1/(z+a)_finite_factorial_series}) converges
for $p\rightarrow\infty$, because $\prod_{j=1}^{p}\frac{b_{j}-a}{b_{j}%
+z}\rightarrow0$ as $\left\vert \frac{b_{j}-a}{b_{j}+z}\right\vert <1$. The
factorial series (\ref{1/(z+a)_finite_factorial_series}) for $\left(
z+a\right)  ^{-1}$ demonstrates that for functions, possessing a series
$\sum_{k=0}^{\infty}a_{k}z^{-k}$ such as Laplace transforms in Section
\ref{sec_factorial_series_Laplace_transforms}, infinitely many factorial
series are possible. Leaving convergence considerations aside for an arbitrary
set $\left\{  b_{j}\right\}  _{1\leq j}$ of complex numbers when
$p\rightarrow\infty$ in (\ref{1/(z+a)_finite_factorial_series}), Cauchy's
integral $f\left(  z\right)  =\frac{1}{2\pi i}\int_{C\left(  z\right)  }%
\frac{f\left(  w\right)  }{w-z}dw$ becomes with
(\ref{1/(z+a)_finite_factorial_series})
\[
f\left(  z\right)  =\frac{1}{2\pi i}\sum_{l=1}^{\infty}\prod_{j=1}%
^{l-1}\left(  b_{j}+z\right)  \int_{C\left(  z\right)  }\frac{f\left(
w\right)  }{\prod_{j=1}^{l}b_{j}+w}dw
\]
If the contour can be closed over the entire $w$-plane and all $b_{j}$ are
different, then $\frac{1}{2\pi i}\int_{C\left(  z\right)  }\frac{f\left(
w\right)  }{\prod_{j=1}^{l}b_{j}+w}dw=\sum_{k=1}^{l}\frac{f\left(
-b_{k}\right)  }{\prod_{j=1;j\neq k}^{l}b_{j}-b_{k}}$ and we formally arrive
at a generalization of Taylor's series%
\begin{equation}
f\left(  z\right)  =\sum_{l=1}^{\infty}\sum_{k=1}^{l}\frac{f\left(
-b_{k}\right)  }{\prod_{j=1;j\neq k}^{l}\left(  b_{j}-b_{k}\right)  }%
\prod_{j=1}^{l-1}\left(  b_{j}+z\right)  \label{Taylor_series_generalized}%
\end{equation}
Explicitly,%
\begin{align*}
f\left(  z\right)   &  =f\left(  -b_{1}\right)  +\frac{f\left(  -b_{1}\right)
-f\left(  -b_{2}\right)  }{b_{2}-b_{1}}\left(  b_{1}+z\right) \\
&  \hspace{0.5cm}+\left(  \frac{f\left(  -b_{1}\right)  }{\left(  b_{2}%
-b_{1}\right)  \left(  b_{3}-b_{1}\right)  }-\frac{f\left(  -b_{2}\right)
}{\left(  b_{2}-b_{1}\right)  \left(  b_{3}-b_{2}\right)  }+\frac{f\left(
-b_{3}\right)  }{\left(  b_{3}-b_{1}\right)  \left(  b_{3}-b_{2}\right)
}\right)  \left(  b_{1}+z\right)  \left(  b_{2}+z\right) \\
&  \hspace{0.5cm}+\left(  \sum_{k=1}^{4}\frac{f\left(  -b_{k}\right)  }%
{\prod_{j=1;j\neq k}^{4}\left(  b_{j}-b_{k}\right)  }\right)  \left(
b_{1}+z\right)  \left(  b_{2}+z\right)  \left(  b_{3}+z\right)  +\cdots
\end{align*}
We omit here the further exploration of (\ref{Taylor_series_generalized}) and
continue with Gilbert's method.

Substitution of (\ref{1/(z+a)_finite_factorial_series}) with $z\rightarrow
z+k$ and $a=x$ in (\ref{mu_Gilbert_Gudermann}) yields, for positive real $z$
and $b_{j}$, a general factorial series%
\[
\mu\left(  z\right)  =\sum_{k=0}^{\infty}\sum_{l=2}^{p}\frac{\int_{0}%
^{1}\left(  \frac{1}{2}-x\right)  \prod_{j=1}^{l-1}\left(  b_{j}-x\right)
dx}{\prod_{j=1}^{l}\left(  b_{j}+z+k\right)  }+\sum_{k=0}^{\infty}\frac
{1}{\prod_{j=1}^{p}\left(  b_{j}+z+k\right)  }\int_{0}^{1}dx\frac{\left(
\frac{1}{2}-x\right)  \prod_{j=1}^{p}\left(  b_{j}-x\right)  }{z+k+x}dx
\]
because $\int_{0}^{1}\left(  \frac{1}{2}-x\right)  dx=0$. If $p\rightarrow
\infty$, the last sum, which is bounded\footnote{As shown below, for
$b_{j}=j-1$, the series reduces to Binet's second series
(\ref{Binet_function_faculteits_expansie}) and the remainder in the integer
$p$,%
\[
R_{p}=\sum_{k=0}^{\infty}\frac{1}{\prod_{j=0}^{p-1}\left(  j+z+k\right)  }%
\int_{0}^{1}dx\frac{\left(  \frac{1}{2}-x\right)  \prod_{j=0}^{p-1}\left(
j-x\right)  }{z+k+x}dx
\]
is nicely bounded by Gilbert \cite{Gilbert1886} as $\left\vert R_{p}%
\right\vert <\frac{1}{64z^{2}}\frac{1}{1+\gamma z+z\log\left(  p-1\right)  }$,
where $\gamma$ is Euler's constant. Gilbert also presents another bound
$\left\vert R_{p}\right\vert <\frac{\Gamma\left(  z\right)  }{64z}%
\frac{e^{\frac{1}{2p}\left(  z+\frac{1}{6}\right)  }}{p^{z}}$.} in
\cite{Gilbert1886}, vanishes and we obtain, for any set $\left\{
b_{j}\right\}  _{1\leq j}$ of positive real numbers, the general factorial
series
\begin{equation}
\mu\left(  z\right)  =\sum_{l=2}^{\infty}\int_{0}^{1}\left(  \frac{1}%
{2}-x\right)  \prod_{j=1}^{l-1}\left(  b_{j}-x\right)  dx\sum_{k=0}^{\infty
}\frac{1}{\prod_{j=1}^{l}\left(  b_{j}+z+k\right)  }
\label{mu_more_general_factorial_series_than_Gilbert}%
\end{equation}

Since $\frac{1}{\prod_{j=1}^{l}\left(  b_{j}+z+k\right)  }=\frac{1}{\left(
b_{1}+z+k\right)  \left(  b_{l}+z+k\right)  \prod_{j=2}^{l-1}\left(
b_{j}+z+k\right)  }$ and with the partial fraction $\frac{1}{\left(
b_{1}+z+k\right)  \left(  b_{l}+z+k\right)  }=\frac{1}{\left(  b_{l}%
-b_{1}\right)  }\left(  \frac{1}{\left(  b_{1}+z+k\right)  }-\frac{1}{\left(
b_{l}+z+k\right)  }\right)  $, we rewrite the infinite $k$-sum as%
\begin{align*}
\sum_{k=0}^{\infty}\frac{1}{\prod_{j=1}^{l}\left(  b_{j}+z+k\right)  }  &
=\sum_{k=0}^{\infty}\frac{1}{\prod_{j=2}^{l-1}\left(  b_{j}+z+k\right)  }%
\frac{1}{\left(  b_{l}-b_{1}\right)  }\left(  \frac{1}{\left(  b_{1}%
+z+k\right)  }-\frac{1}{\left(  b_{l}+z+k\right)  }\right) \\
&  =\sum_{k=0}^{\infty}\frac{\frac{1}{\left(  b_{l}-b_{1}\right)  }}%
{\prod_{j=1}^{l-1}\left(  b_{j}+z+k\right)  }-\sum_{k=0}^{\infty}\frac
{\frac{1}{\left(  b_{l}-b_{1}\right)  }}{\prod_{j=2}^{l}\left(  b_{j}%
+z+k\right)  }\\
&  =\sum_{k=0}^{\infty}\frac{\frac{1}{\left(  b_{l}-b_{1}\right)  }}%
{\prod_{j=1}^{l-1}\left(  b_{j}+z+k\right)  }-\sum_{k=0}^{\infty}\frac
{\frac{1}{\left(  b_{l}-b_{1}\right)  }}{\prod_{j=1}^{l-1}\left(
b_{j+1}+z+k\right)  }%
\end{align*}
If $b_{j+1}=b_{j}+m$, where $m$ is an integer, then the second sum is
\[
\sum_{k=0}^{\infty}\frac{\frac{1}{\left(  b_{l}-b_{1}\right)  }}{\prod
_{j=1}^{l-1}\left(  b_{j+1}+z+k\right)  }=\sum_{k=0}^{\infty}\frac{\frac
{1}{\left(  b_{l}-b_{1}\right)  }}{\prod_{j=1}^{l-1}\left(  b_{j}%
+z+k+m\right)  }=\sum_{k=m}^{\infty}\frac{\frac{1}{\left(  b_{l}-b_{1}\right)
}}{\prod_{j=1}^{l-1}\left(  b_{j}+z+k\right)  }%
\]
and a finite series is found%
\[
\sum_{k=0}^{\infty}\frac{1}{\prod_{j=1}^{l}\left(  b_{j}+z+k\right)  }%
=\frac{1}{\left(  b_{l}-b_{1}\right)  }\sum_{k=0}^{m-1}\frac{1}{\prod
_{j=1}^{l-1}\left(  b_{j}+z+k\right)  }%
\]
Since the solution of the difference equation $b_{j+1}=b_{j}+m$ is
$b_{j}=jm+b_{0}$, we conclude that only if $b_{j}$ is linear in $j$ with
highest coefficient an integer, the infinite $k$-sum can be written as finite
sum. Gilbert \cite[art. 13]{Gilbert1873} has started from
(\ref{1/(z+a)_finite_factorial_series}) with the choice $b_{j}=b+\left(
j-1\right)  p$, in which case%
\[
\sum_{k=0}^{\infty}\frac{1}{\prod_{j=1}^{l}\left(  b_{j}+z+k\right)  }%
=\frac{1}{\left(  l-1\right)  p}\sum_{k=0}^{p-1}\frac{1}{\prod_{j=0}%
^{l-2}\left(  b+jp+z+k\right)  }%
\]
Substituted into the general series
(\ref{mu_more_general_factorial_series_than_Gilbert}), Gilbert arrives at a
general, double sum factorial series%
\begin{equation}
\mu\left(  z\right)  =\frac{1}{p}\sum_{l=1}^{\infty}\sum_{k=0}^{p-1}%
\frac{\frac{1}{l}\int_{0}^{1}\left(  \frac{1}{2}-x\right)  \prod_{j=0}%
^{l-1}\left(  b+jp-x\right)  dx}{\prod_{j=0}^{l-1}\left(  b+jp+z+k\right)  }
\label{Gilbert_factorial_series_Binet_functions}%
\end{equation}
which complements $\varphi(z)=\beta\sum_{m=0}^{\infty}\frac{m!\phi_{m}\left(
\alpha,\beta\right)  }{\prod_{k=0}^{m}(\beta z+\alpha+k)}$ in
(\ref{factorial_series_Laplace_transform_in_alpha_en_beta}) below, for $p>1$.
For $b=0$ and $p=1$, Gilbert's factorial series
(\ref{Gilbert_factorial_series_Binet_functions}) reduces to
\[
\mu\left(  z\right)  =\sum_{l=1}^{\infty}\frac{\frac{1}{l}\int_{0}^{1}\left(
\frac{1}{2}-x\right)  \prod_{j=0}^{l-1}\left(  j-x\right)  dx}{\prod
_{j=0}^{l-1}\left(  j+z\right)  }%
\]
which is Binet's second factorial series
(\ref{Binet_function_faculteits_expansie}). For $b=1$ and $p=1$, Gilbert's
factorial series (\ref{Gilbert_factorial_series_Binet_functions}) reduces to
\[
\mu\left(  z\right)  =\sum_{l=1}^{\infty}\frac{\frac{1}{l}\int_{0}^{1}\left(
\frac{1}{2}-x\right)  \prod_{j=0}^{l-1}\left(  1+j-x\right)  dx}{\prod
_{j=0}^{l-1}\left(  1+j+z\right)  }%
\]
which is Binet's first factorial series (\ref{Binet_original}).

\section{Gilbert's factorial series for the Binet function $\mu\left(
z\right)  $}

\label{sec_infinitely_many_factorial series}We investigate Gilbert's factorial
series (\ref{Gilbert_factorial_series_Binet_functions}) for $p=1$ further.
Several forms and properties of the coefficients in Gilbert's factorial series
are deduced.\hfill

\begin{theorem}
\label{theorem_mu_infinite_convergent_factorial_expansions}Binet's function
$\mu\left(  z\right)  $ possesses infinitely many factorial expansions in the
complex parameter $\alpha$, for $\operatorname{Re}\left(  z\right)  >0$ and
$\operatorname{Re}\left(  \alpha\right)  >-\operatorname{Re}\left(  z\right)
$,%
\[
\mu\left(  z\right)  =\sum_{m=1}^{\infty}\frac{b_{m}\left(  \alpha\right)
}{\prod_{k=0}^{m-1}(z+\alpha+k)}%
\]
where the Binet polynomials in $\alpha$ are%
\begin{equation}
b_{m}(\alpha)=\frac{1}{m}\sum_{k=1}^{m}\left\{  \frac{\alpha^{k+2}-\left(
\alpha-1\right)  ^{k+2}}{k+2}+\frac{\left(  \frac{1}{2}-\alpha\right)  \left(
\alpha^{k+1}-\left(  \alpha-1\right)  ^{k+1}\right)  }{k+1}\right\}  \left(
-1\right)  ^{k-m}S_{m}^{(k)} \label{generalized_Binet_numbers_alpha}%
\end{equation}
and $S_{m}^{(k)}$ is the Stirling Number of the First Kind. Another expression
in terms of the coefficients $c_{k}=b_{k}\left(  0\right)  $ in
(\ref{def_Binet_coefficients}) is%
\begin{equation}
b_{m}\left(  \alpha\right)  =c_{m}+\alpha\sum_{k=1}^{m-1}\binom{m-1}{k-1}%
\prod_{j=1}^{m-k-1}(\alpha+j)c_{k} \label{b_m(alpha)_factorial_sum_cm}%
\end{equation}
The corresponding integral representation is%
\begin{equation}
b_{m}(\alpha)=\frac{1}{m}\int_{\alpha-1}^{\alpha}\left(  x+\left(  \frac{1}%
{2}-\alpha\right)  \right)  \prod_{k=0}^{m-1}(k+x)dx
\label{generalized_Binet_numbers_alpha_integral}%
\end{equation}
In particular, $b_{1}(\alpha)=\frac{1}{12}$ and $b_{2}\left(  \alpha\right)
=\frac{\alpha}{12}$.
\end{theorem}

\textbf{Proof: }We write Binet's integral in (\ref{Binet_function_integral}),
valid for $\operatorname{Re}\left(  z\right)  >0$, as%
\[
\mu\left(  z\right)  =\int_{0}^{1}\frac{\left(  1-u\right)  ^{z+\alpha-1}}%
{u}\left(  1-u\right)  ^{-\alpha}\left(  \frac{\frac{u}{2}-1}{\log\left(
1-u\right)  }-\frac{u}{\log^{2}\left(  1-u\right)  }\right)  du
\]
After substituting the Taylor series
(\ref{exponential_generating_function_b_m(alpha)}) in Appendix
\ref{sec_Taylor_series_c[m]}
\begin{equation}
g_{\alpha}\left(  u\right)  =\left(  1-u\right)  ^{-\alpha}\left(  \frac
{\frac{u}{2}-1}{\log\left(  1-u\right)  }-\frac{u}{\log^{2}\left(  1-u\right)
}\right)  =\sum_{m=1}^{\infty}\frac{b_{m}\left(  \alpha\right)  }{\left(
m-1\right)  !}u^{m} \label{Taylor_series_g_alpha(u)}%
\end{equation}
and following the same steps as in the proof of Theorem
\ref{theorem_Binet_convergent_expansion}, we arrive at
(\ref{mu_factorial_alpha_expansions}). Executing the Cauchy product of the two
Taylor series of $\left(  1-u\right)  ^{-\alpha}$ and $g_{0}\left(  u\right)
=\sum_{m=1}^{\infty}\frac{c_{m}}{\left(  m-1\right)  !}u^{m}$ and equating
corresponding powers in $u$, leads to the factorial expansion
(\ref{b_m(alpha)_factorial_sum_cm}) of the Binet polynomial $b_{m}\left(
\alpha\right)  $ in terms of the coefficients $c_{k}=b_{k}\left(  0\right)  $
in (\ref{def_Binet_coefficients}).

We proceed by deducing (\ref{generalized_Binet_numbers_alpha}). Introducing
the series (\ref{mu_factorial_alpha_expansions}) in the difference formula
(\ref{forwards_difference_mu}), provides a factorial expansion for the function%

\begin{equation}
1+\left(  z+\frac{1}{2}\right)  \log\frac{z}{z+1}=-\sum_{m=1}^{\infty}%
\frac{mb_{m}(\alpha)}{\prod_{k=0}^{m}(z+\alpha+k)}
\label{Binet_expansion_a_log_function}%
\end{equation}
which we rewrite, after denoting $y=z+\alpha$, as%
\[
1-\left(  y+\frac{1}{2}-\alpha\right)  \left\{  \log\left(  1-\frac{\alpha
-1}{y}\right)  -\log\left(  1-\frac{\alpha}{y}\right)  \right\}  =-\sum
_{m=1}^{\infty}\frac{mb_{m}(\alpha)}{\prod_{k=0}^{m}(y+k)}%
\]
We expand now both sides of (\ref{Binet_expansion_a_log_function}) into powers
of $\frac{1}{y}$. The Taylor series around $z_{0}=0$ of $\log\left(
1-x\right)  =-\sum_{k=1}^{\infty}\frac{x^{k}}{k}$, convergent for $\left\vert
z\right\vert <1$, in the left-hand side of
(\ref{Binet_expansion_a_log_function}), leads, for $\max\left(  \left\vert
\alpha-1\right\vert ,\left\vert \alpha\right\vert \right)  <\left\vert
y\right\vert $, to%
\[
1-\left(  y+\frac{1}{2}-\alpha\right)  \log\frac{\left(  1-\frac{\alpha-1}%
{y}\right)  }{\left(  1-\frac{\alpha}{y}\right)  }=-\sum_{k=2}^{\infty}\left(
\frac{\alpha^{k+1}-\left(  \alpha-1\right)  ^{k+1}}{k+1}+\frac{\left(
\frac{1}{2}-\alpha\right)  \left(  \alpha^{k}-\left(  \alpha-1\right)
^{k}\right)  }{k}\right)  \frac{1}{y^{k}}%
\]
Nielsen \cite[band I, p. 68]{NielsenChelsea} derives%
\begin{equation}
\frac{1}{\prod_{k=0}^{m}\left(  z+k\right)  }=\sum_{k=m}^{\infty}\left(
-1\right)  ^{m-k}\mathcal{S}_{k}^{(m)}\frac{1}{z^{k+1}}
\label{Nielsen_factorial_expansion_inverse_powers}%
\end{equation}
where $\mathcal{S}_{k}^{(n)}$ denotes the Stirling Number of the Second Kind
\cite[Sec. 24.1.3 and 24.1.4]{Abramowitz}, which we use in the right-hand side
of (\ref{Binet_expansion_a_log_function})%
\begin{align*}
\sum_{m=1}^{\infty}\frac{mb_{m}(\alpha)}{\prod_{k=0}^{m}(y+k)}  &  =\sum
_{m=1}^{\infty}mb_{m}(\alpha)\sum_{k=m}^{\infty}\left(  -1\right)
^{m-k}\mathcal{S}_{k}^{(m)}\frac{1}{y^{k+1}}\\
&  =\sum_{k=1}^{\infty}\sum_{m=1}^{k}mb_{m}(\alpha)\left(  -1\right)
^{m-k}\mathcal{S}_{k}^{(m)}\frac{1}{y^{k+1}}%
\end{align*}
Equating corresponding powers in $\frac{1}{y}$ of both sides in
(\ref{Binet_expansion_a_log_function}) yields, for $k\geq1$,%
\begin{equation}
\sum_{m=1}^{k}mb_{m}(\alpha)\left(  -1\right)  ^{m-k}\mathcal{S}_{k}%
^{(m)}=\frac{\alpha^{k+2}-\left(  \alpha-1\right)  ^{k+2}}{k+2}+\frac{\left(
\frac{1}{2}-\alpha\right)  \left(  \alpha^{k+1}-\left(  \alpha-1\right)
^{k+1}\right)  }{k+1} \label{recursion_general_coeff_b_m(alpha)_StirlingS2}%
\end{equation}
Finally, after multiplying both sides in
(\ref{recursion_general_coeff_b_m(alpha)_StirlingS2}) by $\left(  -1\right)
^{k}S_{j}^{(k)}$, summing over $k\in\left[  1,j\right]  $, we have%
\[
\sum_{k=1}^{j}\sum_{m=1}^{k}mb_{m}(\alpha)\left(  -1\right)  ^{m}S_{j}%
^{(k)}\mathcal{S}_{k}^{(m)}=\sum_{k=1}^{j}\left\{  \frac{\alpha^{k+2}-\left(
\alpha-1\right)  ^{k+2}}{k+2}+\frac{\left(  \frac{1}{2}-\alpha\right)  \left(
\alpha^{k+1}-\left(  \alpha-1\right)  ^{k+1}\right)  }{k+1}\right\}  \left(
-1\right)  ^{k}S_{j}^{(k)}%
\]
We reverse the $k$- and $m$-summation in the double sum at the left-hand side
\[
\sum_{k=1}^{j}\sum_{m=1}^{k}mb_{m}(\alpha)\left(  -1\right)  ^{m}S_{j}%
^{(k)}\mathcal{S}_{k}^{(m)}=\sum_{m=1}^{j}mb_{m}(\alpha)\left(  -1\right)
^{m}\left(  \sum_{k=m}^{j}S_{j}^{(k)}\mathcal{S}_{k}^{(m)}\right)
\]
invoke the second orthogonality relation for the Stirling numbers \cite[sec.
24.1.4]{Abramowitz}%
\begin{equation}
\sum_{k=m}^{j}S_{j}^{(k)}\mathcal{S}_{k}^{(m)}=\delta_{mj}
\label{orthogonality_StirlingNumbers_2}%
\end{equation}
and obtain $\sum_{k=1}^{j}\sum_{m=1}^{k}mb_{m}(\alpha)\left(  -1\right)
^{m}S_{j}^{(k)}\mathcal{S}_{k}^{(m)}=jc_{j}\left(  a\right)  \left(
-1\right)  ^{j}$, which demonstrates (\ref{generalized_Binet_numbers_alpha}).

The corresponding integral representation of the coefficient $b_{m}(\alpha)$
is translated, via (\ref{integrated_Stirling_gf}), as%
\[
mb_{m}(\alpha)=\int_{-\left(  \alpha-1\right)  }^{-\alpha}\left(  u-\left(
\frac{1}{2}-\alpha\right)  \right)  \prod_{k=0}^{m-1}(k-u)du
\]
After substitution of $x=-u$, we arrive at
(\ref{generalized_Binet_numbers_alpha_integral}).\hfill$\square\medskip$

We now discuss implications of Theorem
\ref{theorem_mu_infinite_convergent_factorial_expansions}. There are two
particularly interesting cases of the Binet polynomial $b_{m}(\alpha)$ in
(\ref{generalized_Binet_numbers_alpha}): for $\alpha=0$,%
\[
b_{m}(0)=c_{m}=\frac{1}{2m}\sum_{k=1}^{m}\left\{  \frac{k}{\left(  k+1\right)
\left(  k+2\right)  }\right\}  \left(  -1\right)  ^{m}S_{m}^{(k)}%
\]
but $\alpha=1$ leads to the original Binet coefficients
(\ref{Binet_coefficient_original}),%
\[
b_{m}(1)=\beta_{m}=\frac{1}{2m}\sum_{k=1}^{m}\left\{  \frac{k}{\left(
k+1\right)  \left(  k+2\right)  }\right\}  \left(  -1\right)  ^{k-m}%
S_{m}^{(k)}%
\]
Since $\left(  -1\right)  ^{k-j}S_{j}^{(k)}$ is a non-negative integer, it
follows that the original Binet coefficients $\beta_{m}=b_{m}(1)$ are all
positive, in contrast to $b_{m}(0)=c_{m}$ in Theorem
\ref{theorem_Binet_coefficient_integral}, whose sum
(\ref{def_Binet_coefficients}) is alternating and does not obviously to lead
to conclusions about the sign. As mentioned earlier, any truncation at $m=K>2$
terms in (\ref{Binet_function_faculteits_expansie}) upper bounds Binet's
function $\mu\left(  z\right)  $, whereas any truncation of $m=K$ terms in
Binet's original expansion (\ref{Binet_original}) lower bounds $\mu\left(
z\right)  $. Hence, for any finite integer $\dot{K}>2$, it holds that%
\[
\sum_{m=1}^{K}\frac{b_{m}\left(  1\right)  }{\prod_{k=1}^{m}(z+k)}<\mu\left(
z\right)  <\sum_{m=1}^{K}\frac{b_{m}\left(  0\right)  }{\prod_{k=0}%
^{m-1}(z+k)}%
\]
which suggests that there may exist a tighter value of $\alpha$ between 0 and
1, explored in Section \ref{sec_Stirling_Binet}.

It follows from (\ref{generalized_Binet_numbers_alpha})%
\[
b_{m}(\alpha)=\frac{\left(  -1\right)  ^{m}}{m}\sum_{k=1}^{m}\left\{
\frac{\alpha^{k+2}-\left(  \alpha-1\right)  ^{k+2}}{k+2}+\frac{\left(
\frac{1}{2}-\alpha\right)  \left(  \alpha^{k+1}-\left(  \alpha-1\right)
^{k+1}\right)  }{k+1}\right\}  \left(  -1\right)  ^{k}S_{m}^{(k)}%
\]
that%
\[
b_{m}(1-\alpha)=\frac{\left(  -1\right)  ^{m}}{m}\sum_{k=1}^{m}\left\{
\frac{\alpha^{k+2}-\left(  \alpha-1\right)  ^{k+2}}{k+2}+\frac{\left(
\frac{1}{2}-\alpha\right)  \left(  \alpha^{k+1}-\left(  \alpha-1\right)
^{k+1}\right)  }{k+1}\right\}  S_{m}^{(k)}%
\]
illustrating, with $\left(  -1\right)  ^{m-k}S_{m}^{(k)}\geq0$, the absence of
symmetry around $\alpha=\frac{1}{2}$. A second observation of
(\ref{generalized_Binet_numbers_alpha}) for $b_{m}(\alpha)$ and the fact that
Stirling numbers are integers is that, if $\alpha\in\mathbb{Q}$ is a rational
number, i.e. $\alpha=\frac{l}{k}$ for integers $l$ and $k$, then the Binet
polynomial $b_{m}(\alpha)$ is also rational.

\subsection{Properties of the Binet polynomial $b_{m}\left(  \alpha\right)  $}

\begin{property}
\label{property_b_m(a)_polynomial_in_a} The Binet polynomial $b_{m}\left(
\alpha\right)  $ in (\ref{generalized_Binet_numbers_alpha}) is a polynomial of
degree $m-1$ in $\alpha$,
\begin{equation}
b_{m}(\alpha)=\sum_{j=0}^{m-1}p_{j}\left(  m\right)  \alpha^{j}
\label{b_m(a)_polynomial_a_order_m}%
\end{equation}
where the coefficients%
\begin{equation}
p_{j}\left(  m\right)  =\frac{1}{j!}\left.  \frac{d^{j}b_{m}(\alpha)}%
{d\alpha^{j}}\right\vert _{\alpha=0}=\frac{\left(  -1\right)  ^{m-1-j}%
}{2m\;j!}\sum_{k=j+1}^{m}\frac{k!\left(  k-j\right)  }{\left(  k+2-j\right)
!}S_{m}^{(k)} \label{coefficient_j_in_polynomial_b_m(a)}%
\end{equation}
from which $p_{m-1}\left(  m\right)  =\frac{1}{12}$ and $p_{m-2}\left(
m\right)  =\frac{1}{12}\binom{m-1}{2}$. An integral form is%
\begin{equation}
p_{j}\left(  m\right)  =\frac{1}{mj!}\int_{-1}^{0}\left(  u+\frac{1}%
{2}\right)  \frac{d^{j}}{du^{j}}\prod_{k=0}^{m-1}(k+u)du
\label{p_j(m)_integral}%
\end{equation}

\end{property}

\textbf{Proof}: Substitution of%
\[
\frac{\alpha^{k+2}}{k+2}-\frac{\left(  \alpha-1\right)  ^{k+2}}{k+2}%
+\frac{\left(  \frac{1}{2}-\alpha\right)  \alpha^{k+1}}{k+1}-\frac{\left(
\frac{1}{2}-\alpha\right)  \left(  \alpha-1\right)  ^{k+1}}{k+1}=-\frac{1}%
{2}\sum_{j=0}^{k}\frac{k!\left(  k-j\right)  }{j!\left(  k+2-j\right)
!}\left(  -1\right)  ^{k-j}\alpha^{j}%
\]
into the Binet polynomial (\ref{generalized_Binet_numbers_alpha}) and using
$S_{m}^{(0)}=0$ for $m>0$ yields
\[
b_{m}(\alpha)=-\frac{1}{2m}\sum_{k=0}^{m}\left\{  \sum_{j=0}^{k}%
\frac{k!\left(  k-j\right)  }{j!\left(  k+2-j\right)  !}\left(  -1\right)
^{k-j}\alpha^{j}\right\}  \left(  -1\right)  ^{k-m}S_{m}^{(k)}%
\]
We reverse the $j$- and $k$-sum, verify that $p_{m}\left(  m\right)  =0$, and
arrive at (\ref{b_m(a)_polynomial_a_order_m}) and
(\ref{coefficient_j_in_polynomial_b_m(a)}).

The integral form (\ref{p_j(m)_integral}) is immediate from the integral
(\ref{generalized_Binet_numbers_alpha_integral}) of $b_{m}\left(
\alpha\right)  $ after substitution of $u=x-\alpha$ as%
\[
p_{j}\left(  m\right)  =\frac{1}{j!}\left.  \frac{d^{j}b_{m}(\alpha)}%
{d\alpha^{j}}\right\vert _{\alpha=0}=\frac{1}{mj!}\int_{-1}^{0}\left(
u+\frac{1}{2}\right)  \left.  \frac{d^{j}}{d\alpha^{j}}\prod_{k=0}%
^{m-1}(k+u+\alpha)\right\vert _{\alpha=0}du
\]
because $\frac{d^{j}}{d\alpha^{j}}\prod_{k=0}^{m-1}(k+u+\alpha)=\frac{d^{j}%
}{du^{j}}\prod_{k=0}^{m-1}(k+u+\alpha)$. Introducing the $j$-th derivative of
the generating function (\ref{genfunc_stirling}), $\frac{1}{j!}\frac{d^{j}%
}{du^{j}}\prod_{k=0}^{m-1}(k+u)=\sum_{k=j}^{m}\binom{k}{j}S_{m}^{(k)}\;\left(
-1\right)  ^{m-k}u^{k-j}$, into (\ref{p_j(m)_integral}) alternatively leads to
(\ref{coefficient_j_in_polynomial_b_m(a)}). $\hfill\square\medskip$

Clearly, if $\alpha=0$, then we find the coefficients $c_{m}=b_{m}\left(
0\right)  =p_{0}\left(  m\right)  $ in (\ref{def_Binet_coefficients}) of
Binet's second factorial expansion again.

\begin{corollary}
\label{corollary_derivates_mu(z)} The $n$-th derivative of Binet's function
$\mu\left(  z\right)  $ is
\begin{equation}
\frac{d^{n}\mu\left(  z\right)  }{dz^{n}}=\left(  -1\right)  ^{n}\sum
_{m=n}^{\infty}\frac{\frac{d^{n}b_{m}\left(  \alpha\right)  }{d\alpha^{n}}%
}{\prod_{k=0}^{m-1}(z+\alpha+k)} \label{mu(z)_derivative_n_in_alpha}%
\end{equation}

\end{corollary}

\textbf{Proof}: Taking the derivative of (\ref{mu_factorial_alpha_expansions})
with respect to $\alpha$ yields%
\[
0=\sum_{m=1}^{\infty}\frac{1}{\prod_{k=0}^{m-1}(z+\alpha+k)}\frac
{db_{m}\left(  \alpha\right)  }{d\alpha}+\sum_{m=1}^{\infty}b_{m}\left(
\alpha\right)  \frac{d}{d\alpha}\frac{1}{\prod_{k=0}^{m-1}(z+\alpha+k)}%
\]
and since $\frac{d}{d\alpha}\frac{1}{\prod_{k=0}^{m-1}(z+\alpha+k)}=\frac
{d}{dz}\frac{1}{\prod_{k=0}^{m-1}(z+\alpha+k)}$ , we conclude that%
\[
\frac{d\mu\left(  z\right)  }{dz}=-\sum_{m=1}^{\infty}\frac{\frac
{db_{m}\left(  \alpha\right)  }{d\alpha}}{\prod_{k=0}^{m-1}(z+\alpha+k)}%
\]
Repeating the argument for an integer $n\geq0$, using $\frac{d^{n}b_{m}\left(
\alpha\right)  }{d\alpha^{n}}=0$ for $n>m$ by Property
\ref{property_b_m(a)_polynomial_in_a}, leads to
(\ref{mu(z)_derivative_n_in_alpha}).\hfill$\square\medskip$

Substituting the polynomial (\ref{b_m(a)_polynomial_a_order_m}) for
$b_{m}\left(  \alpha\right)  $ in Property
\ref{property_b_m(a)_polynomial_in_a} into
(\ref{mu_factorial_alpha_expansions}) indicates, with $b_{0}\left(
\alpha\right)  =0$ and $p_{j}\left(  0\right)  =0$, that%
\[
\mu\left(  z\right)  =\sum_{m=0}^{\infty}\frac{\sum_{j=0}^{m}p_{j}\left(
m\right)  \alpha^{j}}{\prod_{k=0}^{m-1}(z+\alpha+k)}=\sum_{j=0}^{\infty
}\left(  \sum_{m=j}^{\infty}\frac{p_{j}\left(  m\right)  }{\prod_{k=0}%
^{m-1}(z+\alpha+k)}\right)  \alpha^{j}%
\]
Replacing $\alpha\rightarrow z_{0}-z$ and using $p_{m}\left(  m\right)  =0$
yields the Taylor series of $\mu\left(  z\right)  $ around $z_{0}$,%
\[
\mu\left(  z\right)  =\sum_{j=0}^{\infty}\left(  \left(  -1\right)  ^{j}%
\sum_{m=j+1}^{\infty}\frac{p_{j}\left(  m\right)  }{\prod_{k=0}^{m-1}%
(z_{0}+k)}\right)  \left(  z-z_{0}\right)  ^{j}%
\]
The derivatives of Binet's function for $\operatorname{Re}\left(  z\right)
>0$ follow from the Taylor coefficient $\frac{1}{j!}\left.  \frac{d^{j}%
\mu\left(  z\right)  }{dz^{j}}\right\vert _{z=z_{0}}$ or from
(\ref{mu(z)_derivative_n_in_alpha}) and from
(\ref{Binet_mu_integral_Bernoulli_generating_function}),%
\begin{equation}
\frac{d^{j}\mu\left(  z\right)  }{dz^{j}}=\left(  -1\right)  ^{j}%
j!\sum_{m=j+1}^{\infty}\frac{p_{j}\left(  m\right)  }{\prod_{k=0}^{m-1}%
(z+k)}=\frac{\left(  -1\right)  ^{j}}{2}\int_{0}^{\infty}t^{j-1}e^{-zt}\left(
\frac{1+e^{-t}}{1-e^{-t}}-\frac{2}{t}\right)  dt
\label{mu(z)_derivative_j_in_alpha_in_p_integral}%
\end{equation}

\begin{property}
\label{property_b_m(a)_polynomial_coeff_p_j(m)_in_cm} The coefficients
$p_{j}\left(  m\right)  $ of the Binet polynomial $b_{m}(\alpha)=\sum
_{j=0}^{m-1}p_{j}\left(  m\right)  \alpha^{j}$ can be expressed in terms of
the coefficients $c_{m}=b_{m}\left(  0\right)  $ in
(\ref{def_Binet_coefficients}) as%
\begin{equation}
p_{j}\left(  m\right)  =\sum_{l=j}^{m-1}\binom{m-1}{l}S_{l}^{(j)}\left(
-1\right)  ^{l-j}c_{m-l} \label{p_j(m)_in_terms_ck}%
\end{equation}

\end{property}

\textbf{Proof:} Introducing the generating function (\ref{genfunc_stirling})
in (\ref{b_m(alpha)_factorial_sum_cm}) yields%

\[
b_{m}\left(  \alpha\right)  =\sum_{k=1}^{m}\binom{m-1}{k-1}\prod_{j=0}%
^{m-k-1}(\alpha+j)c_{k}=\sum_{k=1}^{m}\binom{m-1}{k-1}c_{k}\sum_{j=0}%
^{m-k}S_{m-k}^{(j)}\left(  -1\right)  ^{m-k-j}\;\alpha^{j}%
\]
Letting $l=m-k$, reversing the sums,%
\[
b_{m}\left(  \alpha\right)  =\sum_{j=0}^{m-1}\left(  \sum_{l=j}^{m-1}%
\binom{m-1}{l}c_{m-l}S_{l}^{(j)}\left(  -1\right)  ^{l-j}\right)  \;\alpha^{j}%
\]
and comparing with the definition $b_{m}(\alpha)=\sum_{j=0}^{m-1}p_{j}\left(
m\right)  \alpha^{j}$ in (\ref{b_m(a)_polynomial_a_order_m}) leads to
(\ref{p_j(m)_in_terms_ck}). $\hfill\square\medskip$

We split the sum (\ref{p_j(m)_in_terms_ck}),%
\[
p_{j}\left(  m\right)  =\frac{1}{12}S_{m-1}^{(j)}\left(  -1\right)
^{m-1-j}-\sum_{l=j}^{m-3}\binom{m-1}{l}\left\vert c_{m-l}\right\vert
S_{l}^{(j)}\left(  -1\right)  ^{l-j}%
\]
which is the difference of two positive numbers that turns out to be positive
for $j>0$. This observation is further substantiated in Property
\ref{property_b_m(a)_polynomial_coeff_p_j(m)_positive} below. If $j=0$, then
$S_{l}^{(0)}=\delta_{0l}$ and the first positive term vanishes for $m>1$ and
$p_{0}\left(  m\right)  =c_{m}<0$. If we substitute the explicit form
(\ref{def_Binet_coefficients}) of $c_{m}$ into (\ref{p_j(m)_in_terms_ck}), we
again arrive at (\ref{coefficient_j_in_polynomial_b_m(a)}) after using the
formula\footnote{Equate corresponding powers of the Taylor series in
$\log^{k+m}(1+u)=\log^{m}(1+u)\log^{k}(1+u)$ from the second generating
function of the Stirling numbers $S_{m}^{\left(  k\right)  }$ in
(\ref{genfunc_stirlinginf}).}
\begin{equation}
S_{j}^{(k+m)}=\frac{m!k!}{\left(  k+m\right)  !}\sum_{l=0}^{j}\binom{j}%
{l}S_{l}^{(m)}S_{j-l}^{(k)} \label{Stirling_binomial_sum_prod_two_Stirling}%
\end{equation}

\begin{property}
\label{property_b_m(a)_polynomial_coeff_p_j(m)_positive} Except for one
negative coefficient, $p_{0}\left(  m\right)  =c_{m}$ for $m>2$, all other
coefficients $p_{j}\left(  m\right)  $ of the Binet polynomial $b_{m}%
(\alpha)=\sum_{j=0}^{m-1}p_{j}\left(  m\right)  \alpha^{j}$ are positive, i.e.
$p_{j}\left(  m\right)  >0$ for $j>0$ and $m>2$. Generally, for $j>1$, it
holds that%
\begin{equation}
p_{j}\left(  m\right)  =\frac{\left(  -1\right)  ^{m-j}}{j\left(  j-1\right)
m}\left(  \left(  m-j+1\right)  S_{m-1}^{(j-2)}+\left(  j-1\right)  \left(
\frac{m}{2}-1\right)  S_{m-1}^{(j-1)}\right)
\label{coefficient_p_j(m)_in_polynomial_b_m(a)_j>1}%
\end{equation}
and also%
\begin{equation}
p_{j}\left(  m\right)  =\frac{\left(  -1\right)  ^{m-j}}{j\left(  j-1\right)
m}\left(  \left(  m-j+1\right)  S_{m}^{(j-1)}+m\left(  m-\frac{j}{2}-\frac
{1}{2}\right)  S_{m-1}^{(j-1)}\right)
\label{coefficient_p_j(m)_in_polynomial_b_m(a)_j>1_versie2}%
\end{equation}
In particular, $p_{2}\left(  m\right)  =\frac{1}{4}\left(  m-2\right)
!\left(  1-\frac{2}{m}\right)  $ and $p_{3}\left(  m\right)  =\frac{1}%
{6}\left(  m-2\right)  !\left(  1-\frac{2}{m}\right)  \sum_{l=2}^{m-2}\frac
{1}{l}$.
\end{property}

\textbf{Proof:} The derivative of the integral representation
(\ref{generalized_Binet_numbers_alpha_integral}) of $b_{m}(\alpha)$ is%
\[
m\frac{db_{m}\left(  \alpha\right)  }{d\alpha}=\left(  \frac{m}{2}%
-1+\alpha\right)  \alpha\prod_{k=1}^{m-2}(k+\alpha)-\int_{\left(
\alpha-1\right)  }^{\alpha}\prod_{k=0}^{m-1}(k+x)dx
\]
from which for $\alpha=0$,%
\[
p_{1}\left(  m\right)  =\left.  \frac{db_{m}\left(  \alpha\right)  }{d\alpha
}\right\vert _{\alpha=0}=\frac{1}{m}\int_{0}^{1}u\prod_{k=1}^{m-1}(k-u)du>0
\]
An additional derivation yields%
\[
m\frac{d^{2}b_{m}\left(  \alpha\right)  }{d\alpha^{2}}=\frac{d}{d\alpha
}\left(  \left(  \frac{m}{2}-1+\alpha\right)  \prod_{k=0}^{m-2}(k+\alpha
)\right)  -\prod_{k=0}^{m-1}(k+\alpha)+\prod_{k=0}^{m-1}(k+\alpha-1)
\]
After employing the logarithmic derivative $\frac{df\left(  x\right)  }%
{dx}=f\left(  x\right)  \frac{d\log f\left(  x\right)  }{dx}$, we find
\[
m\frac{d^{2}b_{m}\left(  \alpha\right)  }{d\alpha^{2}}=\left(  \frac{m}%
{2}-1\right)  \prod_{k=1}^{m-2}(k+\alpha)+\alpha\sum_{j=1}^{m-2}\left(
\frac{m}{2}-1-j\right)  \prod_{k=1;k\neq j}^{m-2}(k+\alpha)
\]
from which, for $m\geq2$, it follows that $p_{2}\left(  m\right)  =\frac{1}%
{2}\left.  \frac{d^{2}b_{m}\left(  \alpha\right)  }{d\alpha^{2}}\right\vert
_{\alpha=0}=\left(  m-2\right)  !\left(  \frac{m-2}{4m}\right)  >0$. We may
continue this tedious process of differentiations to discover closed
expressions for other $p_{j}\left(  m\right)  $ with $j>2$. For example, from
$m\frac{d^{3}b_{m}\left(  \alpha\right)  }{d\alpha^{3}}=$ $\sum_{j=0}%
^{m-2}\left(  \frac{m}{2}-1-j\right)  \sum_{l=0;l\neq j}\prod_{k=0;k\neq
\{j,l\}}^{m-2}(k+\alpha)$, we find $p_{3}\left(  m\right)  =\frac{1}{6}\left(
m-2\right)  !\left(  1-\frac{2}{m}\right)  \sum_{l=2}^{m-2}\frac{1}{l}$, but
that process essentially boils down to computing the Stirling numbers
$S_{m}^{(k)}$ in closed form. The first derivative still contains an integral,
while higher order derivatives are sum of products.

Instead of computing the derivatives $\left.  \frac{d^{n}b_{m}\left(
\alpha\right)  }{d\alpha^{n}}\right\vert _{\alpha=0}$ for $n>1$ from the
integral representation (\ref{generalized_Binet_numbers_alpha_integral}), they
can be deduced more elegantly from the derivative $\frac{d^{n}\mu\left(
z\right)  }{dz^{n}}$ in (\ref{mu(z)_derivative_n_in_alpha}) and from Binet's
integral (\ref{Binet_mu_integral_Bernoulli_generating_function}) as%
\[
\frac{d^{n}\mu\left(  z\right)  }{dz^{n}}=\left(  -1\right)  ^{n}\frac{1}%
{2}\int_{0}^{\infty}t^{n-1}e^{-zt}\left(  \frac{1+e^{-t}}{1-e^{-t}}-\frac
{2}{t}\right)  dt
\]
We mimic Binet's method in the proof of Theorem
\ref{theorem_Binet_convergent_expansion} and invoke Binet's substitution
$e^{-t}=1-u$,%
\[
\frac{d^{n}\mu\left(  z\right)  }{dz^{n}}=-\frac{1}{2}\int_{0}^{1}%
\frac{\left(  1-u\right)  ^{z-1}}{u}\left(  \left(  2-u\right)  \left(
\log\left(  1-u\right)  \right)  ^{n-1}+2u\left(  \log\left(  1-u\right)
\right)  ^{n-2}\right)  du
\]
Since $n>1$, we now use the second generating function
(\ref{genfunc_stirlinginf}) of the Stirling numbers $S_{m}^{\left(  k\right)
}$ and obtain the Taylor series, convergent for $\left\vert u\right\vert <1$,
of the integrand%
\begin{align*}
h\left(  u\right)   &  =\left(  2-u\right)  \left(  \log\left(  1-u\right)
\right)  ^{n-1}+2u\left(  \log\left(  1-u\right)  \right)  ^{n-2}\\
&  =\sum_{m=n-1}^{\infty}\left\{  2\left(  n-1\right)  !\;\frac{S_{m}^{(n-1)}%
}{m}+\left(  n-1\right)  !S_{m-1}^{(n-1)}\;-2\left(  n-2\right)
!S_{m-1}^{(n-2)}\;\right\}  \frac{\left(  -1\right)  ^{m}u^{m}}{\left(
m-1\right)  !}%
\end{align*}
Substituting this Taylor series into the integral, reversing the integration
and summation, invoking the Beta integral and $\frac{\Gamma\left(  z+m\right)
}{\Gamma\left(  z\right)  }=\prod_{k=0}^{m-1}(z+k)$, leads with
(\ref{mu(z)_derivative_n_in_alpha}) for $n>1$ to%
\[
\left.  \frac{d^{n}b_{m}\left(  \alpha\right)  }{d\alpha^{n}}\right\vert
_{\alpha=0}=\left(  -1\right)  ^{m-n-1}\left(  n-2\right)  !\left(
\frac{\left(  n-1\right)  S_{m}^{(n-1)}}{m}-S_{m-1}^{(n-2)}+\frac{\left(
n-1\right)  }{2}S_{m-1}^{(n-1)}\right)
\]
from which $p_{j}\left(  m\right)  =\frac{1}{j!}\left.  \frac{d^{j}%
b_{m}\left(  \alpha\right)  }{d\alpha^{j}}\right\vert _{\alpha=0}$ in
(\ref{coefficient_p_j(m)_in_polynomial_b_m(a)_j>1}) and
(\ref{coefficient_p_j(m)_in_polynomial_b_m(a)_j>1_versie2}) follow, after
eliminating $S_{m}^{(n-1)}$ and $S_{m-1}^{(n-2)}$ respectively by the
recursion $S_{m+1}^{(n)}=S_{m}^{(n-1)}-mS_{m}^{(n)}$ (see e.g.
\cite[24.1.3.II.A]{Abramowitz}). We may verify, by computation of
(\ref{coefficient_p_j(m)_in_polynomial_b_m(a)_j>1}) and
(\ref{coefficient_p_j(m)_in_polynomial_b_m(a)_j>1_versie2}), that
$p_{j}\left(  m\right)  >0$ for $j>1$ and $m>2$\textbf{.}\hfill$\square
\medskip$

The integral (\ref{generalized_Binet_numbers_alpha_integral}) for the Binet
polynomial $b_{m}(\alpha)$ becomes after substitution $u=x+\frac{1}{2}-\alpha$
and reduction of the integration interval $\left[  -\frac{1}{2},\frac{1}%
{2}\right]  $ to $\left[  0,\frac{1}{2}\right]  $,%
\begin{equation}
b_{m}\left(  \alpha+\frac{1}{2}\right)  =\frac{1}{m}\int_{0}^{\frac{1}{2}%
}u\left\{  \prod_{k=0}^{m-1}\left(  k+\alpha+u\right)  -\prod_{k=0}%
^{m-1}\left(  k+\alpha-u\right)  \right\}  du
\label{generalized_Binet_numbers_alpha_integral_alpha+1/2}%
\end{equation}
If $\alpha>0$, then all coefficients are positive, i.e. $b_{m}\left(
\alpha+\frac{1}{2}\right)  >0$ for $m\geq1$, because $\prod_{k=0}^{m-1}\left(
k+\alpha+u\right)  >\prod_{k=0}^{m-1}\left(  k+\alpha-u\right)  $. For larger
negative $\alpha<0$, the Binet polynomial $b_{m}\left(  \alpha\right)  $
starts oscillating and determining the sign is more complicated. The
asymptotic behavior of $b_{m}\left(  \alpha+\frac{1}{2}\right)  $ for large
$m$ is analyzed in Appendix \ref{sec_asymptotic_bm(alpha)_large_m}.

\begin{property}
\label{property_b_m(a)_polynomial_real_zeros} The Binet polynomial%
\[
b_{m}(\alpha)=\sum_{j=0}^{m-1}p_{j}\left(  m\right)  \alpha^{j}=\frac{1}%
{12}\prod_{j=1}^{m-1}\left(  \alpha-\xi_{j}\left(  m\right)  \right)
\]
has $m-1$ real, distinct zeros $\xi_{1}\left(  m\right)  >\xi_{2}\left(
m\right)  >\cdots>\xi_{m-1}\left(  m\right)  $.
\end{property}

\textbf{Proof}: We apply the generalized mean-value theorem\footnote{If
$\varphi\left(  x\right)  $ is non-negative for $x\in$ $\left[  a,b\right]  $,
then $\int_{a}^{b}f\left(  x\right)  \varphi\left(  x\right)  dx=f\left(
\theta\right)  \int_{a}^{b}\varphi\left(  x\right)  dx$ for some $\theta$
obeying $a<\theta<b$.} \cite[p. 321]{Hardy_pure_math} to the integral
(\ref{generalized_Binet_numbers_alpha_integral_alpha+1/2}) for $\alpha\geq0$,
\begin{equation}
b_{m}\left(  \frac{1}{2}+\alpha\right)  =\frac{1}{8m}\left.  \left(
\prod_{k=0}^{m-1}(k+\alpha+\theta_{m})-\prod_{k=0}^{m-1}(k+\alpha-\theta
_{m})\right)  \right\vert _{0<\theta_{m}<\frac{1}{2}}
\label{b_m(alpha+1/2)_generalized_mean_value_theorem}%
\end{equation}
which is the central difference $\delta_{h}f(x)=f\left(  x+\frac{h}{2}\right)
-f\left(  x-\frac{h}{2}\right)  $ with step $h=2\theta_{m}<1$ of the
polynomial $f\left(  x\right)  =\prod_{k=0}^{m-1}(k+x)$ of degree $m$ in $x$
with real zeros at the integers $x=-k$, for $0\leq k\leq m-1$. In a region
containing the distinct zeros, the polynomial $f\left(  x\right)  $ oscillates
below and above the real axis, as well as its shifted companion $f\left(
x+h\right)  $ with step $h$ smaller than the distance between the zeros. This
means that $f\left(  x\right)  $ and $f\left(  x+h\right)  $ will intersect
$m-1$ times at distinct points, implying that $b_{m}\left(  \alpha\right)  $
has $m-1$ real zeros in the interval $\left[  \frac{1}{2}-m,\frac{1}%
{2}\right]  $.\hfill$\square\medskip$

The sum of the zeros equals $\sum_{j=1}^{m-1}\xi_{j}\left(  m\right)  =$
$-\binom{m-1}{2}$ and their product $c_{m}=\frac{1}{12}\prod_{j=1}%
^{m-1}\left(  -\xi_{j}\left(  m\right)  \right)  $. We found that the Binet
polynomial $b_{2m}\left(  \alpha\right)  $ has a \textquotedblleft center
zero\textquotedblright\ equal to $\xi_{m}\left(  2m\right)  =-\left(
m-1\right)  $ and that $\xi_{m-1}\left(  m\right)  >-m-1$. Since, for $m>2$,
all coefficients $p_{j}\left(  m\right)  >0$ for $j>0$ and $p_{0}\left(
m\right)  <0$ and all zeros are real, we conclude \cite[art. 218, p.
289]{PVM_graphspectra} that the largest zero is positive, $\xi_{1}\left(
m\right)  >0$, while all others are negative, $\xi_{j}\left(  m\right)  <0$
for $j>1$.

We end this section with Property
\ref{property_generalized_Binet_Bernoulli_polynomials} that relates the Binet
polynomials $b_{m}\left(  \alpha\right)  $ to Bernoulli polynomials
$B_{n}\left(  \alpha\right)  $,

\begin{property}
\label{property_generalized_Binet_Bernoulli_polynomials}The Binet polynomial
$b_{m}\left(  \alpha\right)  $ in (\ref{generalized_Binet_numbers_alpha}) of
the generalized factorial series (\ref{mu_factorial_alpha_expansions}) of
Binet's function $\mu\left(  z\right)  $ for $\operatorname{Re}\left(
z\right)  >0$ can be expressed in terms of Bernoulli polynomials as%
\begin{equation}
b_{m}\left(  \alpha\right)  =\left(  -1\right)  ^{m-1}\sum_{k=0}^{m-1}%
\frac{S_{m-1}^{(k)}\left(  -1\right)  ^{k}}{\left(  k+1\right)  \left(
k+2\right)  }\left(  B_{k+2}(\alpha)-\alpha^{k+2}+\frac{1}{2}\left(
k+2\right)  \alpha^{k+1}\right)
\label{generalized_Binet_numbers_alpha_Bernoulli_polynomials}%
\end{equation}

\end{property}

\textbf{Proof}: After substituting Nielsen's expansion
(\ref{Nielsen_factorial_expansion_inverse_powers}) into Binet's convergent
series (\ref{Binet_function_faculteits_expansie})
\[
\mu\left(  z\right)  =\sum_{m=1}^{\infty}\frac{b_{m}\left(  \alpha\right)
}{\prod_{k=0}^{m-1}(z+\alpha+k)}=\sum_{m=0}^{\infty}b_{m+1}\left(
\alpha\right)  \sum_{k=m}^{\infty}\left(  -1\right)  ^{m-k}\mathcal{S}%
_{k}^{(m)}\frac{1}{\left(  z+\alpha\right)  ^{k+1}}%
\]
Using $\frac{1}{\left(  z+\alpha\right)  ^{k+1}}=\sum_{n=k}^{\infty}\binom
{n}{k}\left(  -\alpha\right)  ^{n-k}\frac{1}{z^{n+1}}$ in \cite[24.1.1.B]%
{Abramowitz} leads to%
\begin{align*}
\mu\left(  z\right)   &  =\sum_{m=0}^{\infty}b_{m+1}\left(  \alpha\right)
\sum_{k=0}^{\infty}\left(  -1\right)  ^{m-k}\mathcal{S}_{k}^{(m)}\sum
_{n=k}^{\infty}\binom{n}{k}\left(  -\alpha\right)  ^{n-k}\frac{1}{z^{n+1}}\\
&  =\sum_{m=0}^{\infty}b_{m+1}\left(  \alpha\right)  \sum_{n=0}^{\infty
}\left(  \sum_{k=0}^{n}\binom{n}{k}\left(  -1\right)  ^{m+n}\mathcal{S}%
_{k}^{(m)}\alpha^{n-k}\right)  \frac{1}{z^{n+1}}%
\end{align*}
and, since $\mathcal{S}_{k}^{(m)}=0$ if $m>k$, to
\[
\mu\left(  z\right)  =\sum_{n=0}^{\infty}\left(  \sum_{k=0}^{n}\binom{n}%
{k}\left\{  \sum_{m=0}^{k}b_{m+1}\left(  \alpha\right)  \left(  -1\right)
^{m+n}\mathcal{S}_{k}^{(m)}\alpha^{n-k}\right\}  \right)  \frac{1}{z^{n+1}}%
\]

Equating corresponding powers of $\frac{1}{z}$ in the above and in Stirling's
asymptotic series $\mu\left(  z\right)  =\sum_{n=1}^{\infty}\frac{B_{n+1}%
}{n\left(  n+1\right)  z^{n}}$ in (\ref{Binet_asymptotic_Stirling_expansion})
indicates that%
\[
\frac{B_{n+2}\left(  -1\right)  ^{n}\left(  \frac{1}{\alpha}\right)  ^{n}%
}{\left(  n+1\right)  \left(  n+2\right)  }=\sum_{k=0}^{n}\binom{n}{k}\left(
\frac{1}{\alpha}\right)  ^{k}\sum_{m=0}^{k}b_{m+1}\left(  \alpha\right)
\left(  -1\right)  ^{m}\mathcal{S}_{k}^{(m)}%
\]
Binomials inversion \cite[chap. 2]{Riordan}, $a_{n}=\sum_{k=0}^{n}{\binom
{n}{k}}b_{k}\Leftrightarrow b_{n}=\left(  -1\right)  ^{n}\sum_{k=0}%
^{n}(-1)^{k}{\binom{n}{k}}a_{k}$, yields%
\[
\left(  \frac{1}{\alpha}\right)  ^{k}\sum_{m=0}^{k}b_{m+1}\left(
\alpha\right)  \left(  -1\right)  ^{m}\mathcal{S}_{k}^{(m)}=\left(  -1\right)
^{k}\sum_{l=0}^{k}\binom{k}{l}\frac{B_{l+2}\left(  \frac{1}{\alpha}\right)
^{l}}{\left(  l+1\right)  \left(  l+2\right)  }%
\]
and%
\begin{align*}
\sum_{m=0}^{k}b_{m+1}\left(  \alpha\right)  \left(  -1\right)  ^{m}%
\mathcal{S}_{k}^{(m)}  &  =\left(  -1\right)  ^{k}\sum_{l=0}^{k}\binom{k}%
{l}\frac{B_{l+2}\alpha^{k-l}}{\left(  l+1\right)  \left(  l+2\right)  }\\
&  =\frac{\left(  -1\right)  ^{k}}{\left(  k+1\right)  \left(  k+2\right)
}\sum_{l=2}^{k+2}\binom{k+2}{l}B_{l}\alpha^{k+2-l}%
\end{align*}
With the definition \cite[23.1.7]{Abramowitz}\ of the Bernoulli polynomials,
$B_{n}(\alpha)=\sum_{l=0}^{n}\binom{n}{l}B_{l}\alpha^{n-l}$, we
obtain\footnote{The Bernoulli numbers can be written in terms of the Stirling
numbers $\mathcal{S}_{k}^{(m)}$ of the second Kind \cite[p. 220]{Comtet},
\[
B_{k}=\sum_{m=1}^{k}\frac{(-1)^{m}m!\,\mathcal{S}_{k}^{(m)}}{m+1}\hspace
{1cm}\text{for }k>0
\]
}
\begin{equation}
\sum_{m=0}^{k}b_{m+1}\left(  \alpha\right)  \left(  -1\right)  ^{m}%
\mathcal{S}_{k}^{(m)}=\frac{\left(  -1\right)  ^{k}}{\left(  k+1\right)
\left(  k+2\right)  }\left(  B_{k+2}(\alpha)-\alpha^{k+2}+\frac{1}{2}\left(
k+2\right)  \alpha^{k+1}\right)
\label{Bernoulli_polynomial_as_generalized_Binet_coeffient}%
\end{equation}
Formula (\ref{Bernoulli_polynomial_as_generalized_Binet_coeffient}) expresses
the Bernoulli polynomial $B_{n}\left(  \alpha\right)  $ in terms of the Binet
polynomial $b_{m}\left(  \alpha\right)  $.

We invert relation (\ref{Bernoulli_polynomial_as_generalized_Binet_coeffient})
to find the Binet polynomial $b_{m}\left(  \alpha\right)  $. After multiplying
both sides by $S_{j}^{(k)}$, summing over $k\in\left[  0,j\right]  $, we have%
\[
\sum_{k=0}^{j}\sum_{m=0}^{k}b_{m+1}\left(  \alpha\right)  \left(  -1\right)
^{m}S_{j}^{(k)}\mathcal{S}_{k}^{(m)}=\sum_{k=0}^{j}\frac{S_{j}^{(k)}\left(
-1\right)  ^{k}}{\left(  k+1\right)  \left(  k+2\right)  }\left(
B_{k+2}(\alpha)-\alpha^{k+2}+\frac{1}{2}\left(  k+2\right)  \alpha
^{k+1}\right)
\]
Reversing the $k$- and $m$-sum and applying the second orthogonality formula
(\ref{orthogonality_StirlingNumbers_2}) yields%
\[
\sum_{k=0}^{j}\sum_{m=0}^{k}b_{m+1}\left(  \alpha\right)  \left(  -1\right)
^{m}S_{j}^{(k)}\mathcal{S}_{k}^{(m)}=\sum_{m=0}^{j}b_{m+1}\left(
\alpha\right)  \left(  -1\right)  ^{m}\left(  \sum_{k=m}^{j}S_{j}%
^{(k)}\mathcal{S}_{k}^{(m)}\right)  =b_{j+1}\left(  \alpha\right)  \left(
-1\right)  ^{j}%
\]
from which Property \ref{property_generalized_Binet_Bernoulli_polynomials}
follows.\hfill$\square\medskip$

\section{Digamma and polygamma functions}

\label{sec_derivatives_mu_digamma_polygamma}We present the corresponding
factorial series for the digamma function $\psi\left(  z\right)  =\frac
{\Gamma^{\prime}\left(  z\right)  }{\Gamma\left(  z\right)  }$ and for the
polygamma function, defined as $\psi^{(n)}(z)=$ $\frac{d^{n+1}\ln\Gamma
(z)}{dz^{n+1}}$ with $\psi^{(0)}(z)=\psi(z)$. We confine ourselves to the case
$\alpha=0$.

Differentiation of (\ref{def_Binet_mu(z)}) with respect to $z$ expresses the
digamma function $\psi\left(  z\right)  =\frac{\Gamma^{\prime}\left(
z\right)  }{\Gamma\left(  z\right)  }$ in terms of the Binet function
$\mu\left(  z\right)  $ as
\[
\psi\left(  z\right)  =\log z-\frac{1}{2z}+\mu^{\prime}\left(  z\right)
\]
Introducing the factorial expansion
(\ref{mu(z)_derivative_j_in_alpha_in_p_integral}) for $j=1$ gives%
\begin{equation}
\psi\left(  z\right)  =\log z-\frac{1}{2z}-\sum_{m=2}^{\infty}\frac
{p_{1}\left(  m\right)  }{\prod_{k=0}^{m-1}(z+k)}
\label{digamma_factorial_expansion}%
\end{equation}
Explicitly with a few coefficients (\ref{coefficient_j_in_polynomial_b_m(a)})
of $p_{1}\left(  m\right)  $,
\[
\psi\left(  z\right)  =\log z-\frac{1}{2z}-\frac{1}{12z\left(  z+1\right)
}-\frac{1}{12z\left(  z+1\right)  \left(  z+2\right)  }-\frac{19}{120z\left(
z+1\right)  \left(  z+2\right)  \left(  z+3\right)  }-\sum_{m=5}^{\infty}%
\frac{p_{1}\left(  m\right)  }{\prod_{k=0}^{m-1}(z+k)}%
\]
is the convergent companion for $\operatorname{Re}\left(  z\right)  >0$ of the
asymptotic series \cite[6.3.18]{Abramowitz}
\begin{equation}
\psi\left(  z\right)  \sim\log z-\frac{1}{2z}+\sum_{m=1}^{\infty}\frac{B_{2m}%
}{2mz^{2m}}=\log z-\frac{1}{2z}-\frac{1}{12z^{2}}+\frac{1}{120z^{4}}-\frac
{1}{252z^{6}}+\cdots\label{digamma_asymptotic}%
\end{equation}
Since $p_{1}\left(  m\right)  >0$, truncation of $\sum_{m=2}^{K}\frac
{p_{1}\left(  m\right)  }{\prod_{k=0}^{m-1}(z+k)}$ in
(\ref{digamma_factorial_expansion}) after any \thinspace$K\geq2$ provides an
upper bound for the digamma function $\psi\left(  z\right)  $. For integer
$z=n$ in (\ref{digamma_factorial_expansion}) for which $\psi\left(  n\right)
=\sum_{k=1}^{n-1}\frac{1}{k}-\gamma$ (see \cite[6.3.2]{Abramowitz}), the
harmonic series is
\begin{equation}
\sum_{k=1}^{n}\frac{1}{k}=\log n+\gamma+\frac{1}{2n}-\lim_{K\rightarrow\infty
}\sum_{m=2}^{K}\frac{\left(  n-1\right)  !p_{1}\left(  m\right)  }{\left(
n-1+m\right)  !} \label{Harmonic_series}%
\end{equation}
and the $m$-sum converges rapidly, also for relatively small $n$, but
increasingly fast for larger $n$. For example, with $K=5$ terms evaluated in
(\ref{Harmonic_series}), the error is less that $10^{-6}$ for $n=10$ and less
than $10^{-11}$ for $n=100$.

Since $b_{m}\left(  \alpha\right)  $ and the derivatives $\frac{db_{m}\left(
\alpha\right)  }{d\alpha}$ contain integrals as illustrated in the proof of
Property \ref{property_b_m(a)_polynomial_coeff_p_j(m)_positive}, the
functional regime of $\psi^{(n)}(z)$ for $n>2$ is different than for $n\leq2$.
Consequently, asymptotic series for $n>2$ disappear and convergent factorial
series loose their attractiveness, because convergent power series exist.
Indeed, in terms of the Binet function and starting from $\psi^{(0)}%
(z)=\psi\left(  z\right)  =\log z-\frac{1}{2z}+\mu^{\prime}\left(  z\right)
$, it holds for $n\geq1$ that%
\[
\psi^{(n)}(z)=\left(  -1\right)  ^{n-1}\left(  n-1\right)  !z^{-n}+\frac{1}%
{2}\left(  -1\right)  ^{n-1}n!z^{-n-1}+\frac{d^{n+1}}{dz^{n+1}}\mu\left(
z\right)
\]
Introducing the factorial expansion
(\ref{mu(z)_derivative_j_in_alpha_in_p_integral}), valid for $n\geq1$, yields%
\[
\psi^{(n)}(z)=\left(  n-1\right)  !\left(  -1\right)  ^{n-1}\left(  \frac
{1}{z^{n}}+\frac{1}{2}\frac{n}{z^{n+1}}+\left(  n+1\right)  n\sum
_{m=n+2}^{\infty}\frac{p_{n+1}\left(  m\right)  }{\prod_{k=0}^{m-1}%
(z+k)}\right)
\]
The factorial series for $\psi^{(n)}(z)$ converges slower for larger $n$, is
more complicated and only valid for $\operatorname{Re}\left(  z\right)  >0$ in
contrast to $\psi^{(n)}(z)=n!(-1)^{n-1}\sum_{k=0}^{\infty}\frac{1}%
{(z+k)^{n+1}}$ in (\ref{polygamma_series}) that converges for all complex $z$,
except for the poles at $z=-k$ for integer $k\geq0$.

\section{Stirling's asymptotic and Binet's generalized factorial series}

\label{sec_Stirling_Binet}

In this section, we will compare the accuracy of the generalized Binet
factorial expansion (\ref{mu_factorial_alpha_expansions}) in terms of the
error
\[
e_{\alpha}\left(  z,K\right)  =\left\vert \mu\left(  z\right)  -\sum_{m=1}%
^{K}\frac{b_{m}\left(  \alpha\right)  }{\prod_{k=0}^{m-1}(z+\alpha
+k)}\right\vert
\]
which is a function of $z$, the \textquotedblleft free\textquotedblright%
\ parameter $\alpha$ and the number $K$ of terms evaluated. We assume here
that $K$ is finite. Similarly, we denote the error of Stirling's asymptotic
series (\ref{Binet_asymptotic_Stirling_expansion}) by
\[
e_{\text{Stirling}}\left(  z,K\right)  =\left\vert \mu\left(  z\right)
-\sum_{m=1}^{K}\frac{B_{2m}}{\left(  2m-1\right)  \left(  2m\right)  z^{2m-1}%
}\right\vert
\]
We are interested in the \textquotedblleft best\textquotedblright\ value
$\alpha^{\ast}$ of the parameter $\alpha$ for which the error $e_{\alpha
}\left(  z,K\right)  $ is minimal. The asymptotic, diverging nature of the
Stirling approximation allows us to compute the number $K^{\ast}\left(
z\right)  $ of terms that minimizes the error, i.e. for any real $z$,
$e_{\text{Stirling}}\left(  z,K\right)  \geq$ $e_{\text{Stirling}}\left(
z,K^{\ast}\left(  z\right)  \right)  $. Thus, we will take the best possible
performance with minimal error $e_{\text{Stirling}}\left(  z,K^{\ast}\left(
z\right)  \right)  $ as a benchmark to compare $e_{\alpha}\left(  z,K^{\ast
}\left(  z\right)  \right)  $ as a function of $z$ and $\alpha$.

Stirling's asymptotic expansion (\ref{Binet_asymptotic_Stirling_expansion})%
\begin{equation}
\mu\left(  z\right)  \simeq\frac{1}{12z}-\frac{1}{360z^{3}}+\frac{1}%
{1260z^{5}}+\sum_{m=4}^{K}\frac{B_{2m}}{\left(  2m-1\right)  \left(
2m\right)  z^{2m-1}} \label{Stirling_explicit}%
\end{equation}
can be compared to the generalized Binet factorial series
(\ref{mu_factorial_alpha_expansions}) with a same number $K$ of terms,%
\begin{equation}
\mu\left(  z\right)  \simeq\frac{1}{12(z+\alpha)}+\frac{\alpha}{12(z+\alpha
)(z+\alpha+1)}+\sum_{m=3}^{K}\frac{b_{m}\left(  \alpha\right)  }{\prod
_{k=0}^{m-1}(z+\alpha+k)} \label{Binet_algemeen_explicit}%
\end{equation}
In particular, for $\alpha=0$ and $b_{m}\left(  0\right)  =c_{m}$, where
$c_{2}=0$,%
\begin{equation}
\mu\left(  z\right)  \simeq\frac{1}{12z}-\frac{1}{360z\left(  z+1\right)
\left(  z+2\right)  }-\frac{1}{120z\left(  z+1\right)  \left(  z+2\right)
\left(  z+3\right)  }+\sum_{m=5}^{K}\frac{c_{m}}{\prod_{k=0}^{m-1}(z+k)}
\label{Binet_explicit}%
\end{equation}
we observe that the first two coefficients in Stirling's asymptotic
(\ref{Stirling_explicit}) and Binet's convergent (\ref{Binet_explicit})
expansion are the same. Moreover, Stirling's asymptotic
(\ref{Stirling_explicit}) has alternating terms -- the Bernoulli numbers
$B_{2m}=\left(  -1\right)  ^{m-1}\left\vert B_{2m}\right\vert $ alternate --,
in contrast to (\ref{Binet_algemeen_explicit}), where $b_{m}\left(
\alpha\right)  $ changes sign at most once with increasing $m$. While
Stirling's expansion (\ref{Stirling_explicit}) is an asymptotic and
approximate series with a best possible, non-zero error $e_{\text{Stirling}%
}\left(  z,K^{\ast}\left(  z\right)  \right)  $, Binet's factorial, convergent
series (\ref{Binet_algemeen_explicit}) can always beat the accuracy
$e_{\text{Stirling}}\left(  z,K^{\ast}\left(  z\right)  \right)  $ for any
$\operatorname{Re}\left(  z\right)  >0$ if the number $K$ of terms is
sufficiently large. Therefore, we investigate whether Binet's series
(\ref{Binet_algemeen_explicit}) with the same number $K^{\ast}\left(
z\right)  $ of terms can achieve a similar accuracy as Stirling's expansion
(\ref{Stirling_explicit}) with optimal number $K^{\ast}\left(  z\right)  $ of
terms.
%TCIMACRO{\FRAME{fhFU}{16.1276cm}{10.0803cm}{0pt}{\Qcb{Optimal accuracy of
%Stirling's asymptotic expansion: the optimal number of terms $K^{\ast}\left(
%z\right)  $ (left axis) to achieve the lowest error $e_{\text{Stirling}%
%}\left(  z,K^{\ast}\left(  z\right)  \right)  $ (right axis). The error
%$e_{0}\left(  z,K^{\ast}\left(  z\right)  \right)  $ and $e_{1}\left(
%z,K^{\ast}\left(  z\right)  \right)  $ of Binet's two factorial expansions for
%$\alpha=0$ and $\alpha=1$ are added (right axis). The legend shows linear fits
%of the data.}}{\Qlb{Fig_Stirlingoptimalaccuracy}}{stirlingoptimalaccuracy.eps}%
%{\special{ language "Scientific Word";  type "GRAPHIC";
%maintain-aspect-ratio TRUE;  display "USEDEF";  valid_file "F";
%width 16.1276cm;  height 10.0803cm;  depth 0pt;  original-width 6.2881in;
%original-height 3.9202in;  cropleft "0";  croptop "1";  cropright "1";
%cropbottom "0";
%filename 'StirlingOptimalAccuracy.eps';file-properties "XNPEU";}}}%
%BeginExpansion
\begin{figure}
[h]
\begin{center}
\includegraphics[
height=10.0803cm,
width=16.1276cm
]%
{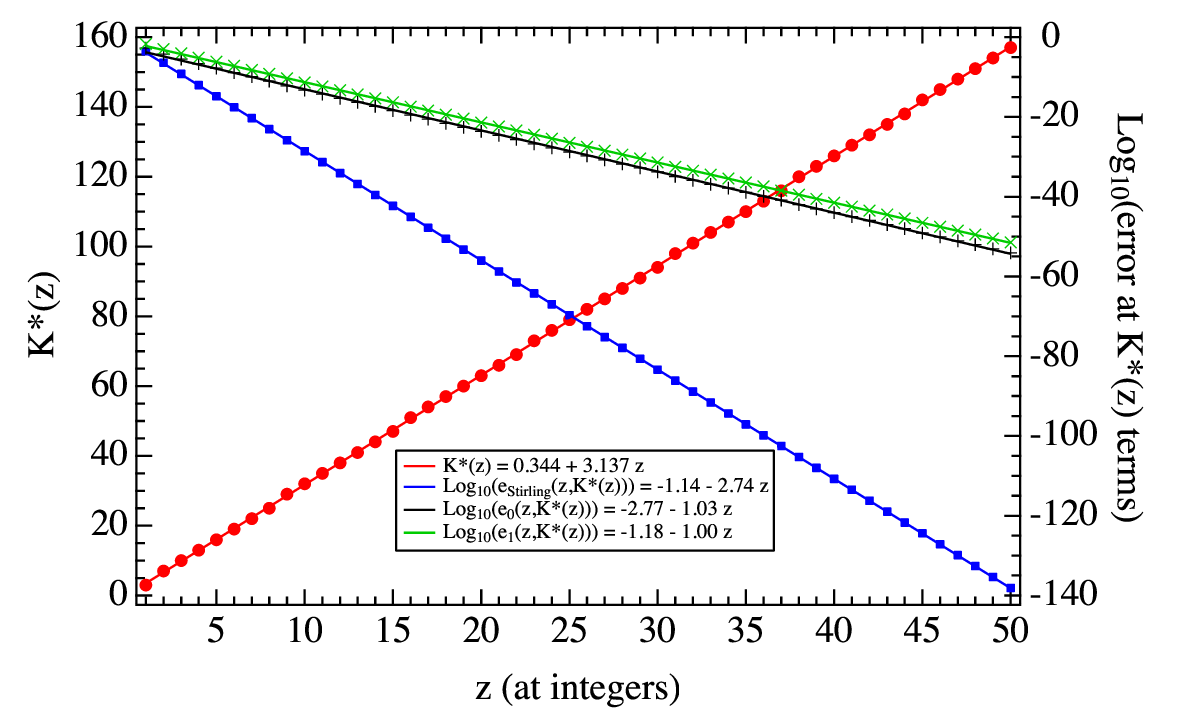}%
\caption{Optimal accuracy of Stirling's asymptotic expansion: the optimal
number of terms $K^{\ast}\left(  z\right)  $ (left axis) to achieve the lowest
error $e_{\text{Stirling}}\left(  z,K^{\ast}\left(  z\right)  \right)  $
(right axis). The error $e_{0}\left(  z,K^{\ast}\left(  z\right)  \right)  $
and $e_{1}\left(  z,K^{\ast}\left(  z\right)  \right)  $ of Binet's two
factorial expansions for $\alpha=0$ and $\alpha=1$ are added (right axis). The
legend shows linear fits of the data.}%
\label{Fig_Stirlingoptimalaccuracy}%
\end{center}
\end{figure}
%EndExpansion
Fig. \ref{Fig_Stirlingoptimalaccuracy} shows the logarithm (in basis 10) of
the error (right axis), evaluated at the optimal number of terms $K^{\ast
}\left(  z\right)  $ versus $z$ (left axis). The comparison of Stirling's
asymptotic series with Binet's two factorial expansions (\ref{Binet_original})
for $\alpha=1$ and (\ref{Binet_function_faculteits_expansie}) for $\alpha=0$
clearly illustrates the amazing superiority of Stirling's asymptotic series.

\begin{figure}
[ptb]
\begin{center}
\includegraphics[
height=3.6781in,
width=6.1194in
]%
{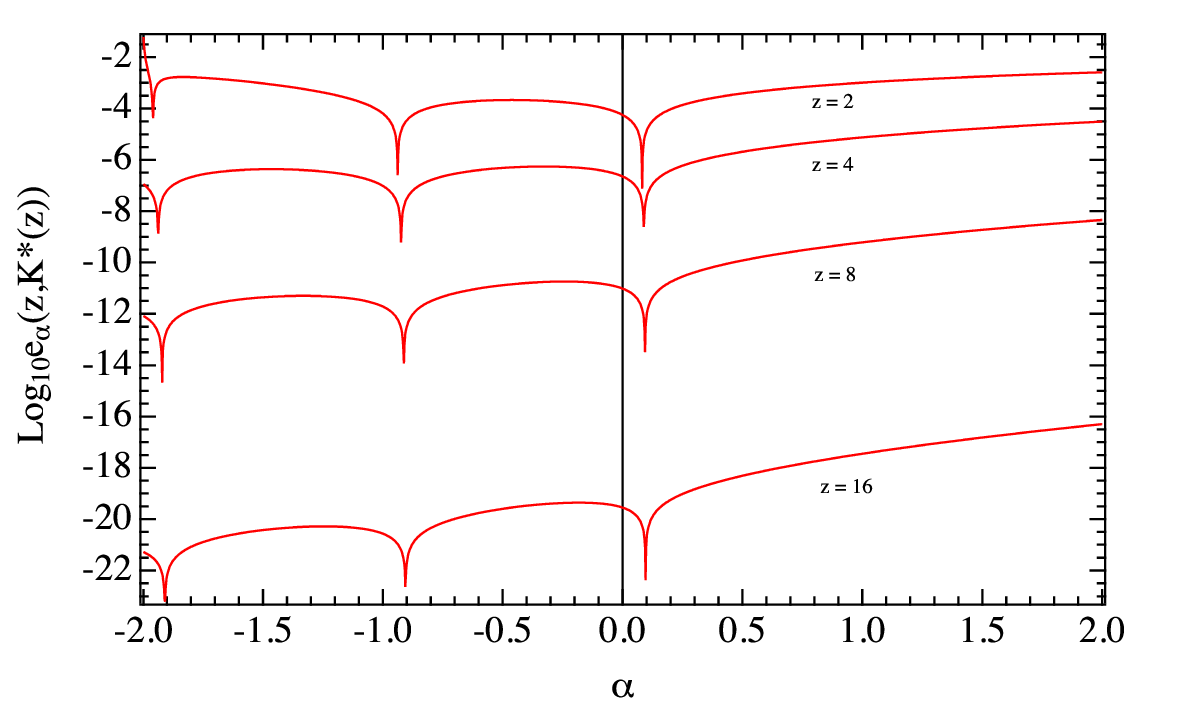}%
\caption{The logarithm of the error $e_{\alpha}\left(  z,K^{\ast}\left(
z\right)  \right)  $ versus $\alpha$ for various $z=\left\{  2,4,8,16\right\}
$.}%
\label{Fig_Binetalpha}%
\end{center}
\end{figure}

Fig. \ref{Fig_Binetalpha} illustrates that the logarithm of the error
$e_{\alpha}\left(  z,K^{\ast}\left(  z\right)  \right)  $ versus $\alpha
\in\left[  -2,2\right]  $ varies considerably. In particular around the zeros
$\xi_{1}\left(  m\right)  $, $\xi_{2}\left(  m\right)  $ and $\xi_{3}\left(
m\right)  $ of the polynomial $b_{m}\left(  \alpha\right)  $ in the interval
$\left[  -2,2\right]  $, the error of the generalized Binet factorial
expansion (\ref{Binet_algemeen_explicit}) decreases sharply. The minimal
errors, peaked around the zeros, only shift little in $\alpha$ for various
$K^{\ast}\left(  z\right)  $ and $z$.
%TCIMACRO{\FRAME{ftbpFU}{6.1194in}{3.6781in}{0pt}{\Qcb{The logarithm of the
%error $e_{\alpha}\left(  z,K^{\ast}\left(  z\right)  \right)  $ versus
%$\alpha$ for various $z=\left\{  2,4,8,16\right\}  $.}}{\Qlb{Fig_Binetalpha}%
%}{binetalpha.eps}{\special{ language "Scientific Word";  type "GRAPHIC";
%maintain-aspect-ratio TRUE;  display "USEDEF";  valid_file "F";
%width 6.1194in;  height 3.6781in;  depth 0pt;  original-width 7.862in;
%original-height 4.7115in;  cropleft "0";  croptop "1";  cropright "1";
%cropbottom "0";  filename 'Binetalpha.eps';file-properties "XNPEU";}}}%
%BeginExpansion
%EndExpansion

Numerically for $m>10$, we found that the largest and only positive zero
$\xi_{1}\left(  m\right)  \in\left[  0.08,0.1\right]  $ and that $\xi
_{1}\left(  m\right)  $ attains its largest value 0.0963016 at $m=72$ and
$\xi_{1}\left(  m\right)  $ slowly decreases for $m>72$. For any finite $K$,
there exists a value of $\alpha$ around $\xi_{1}\left(  K\right)  $ that
minimizes the error%
\begin{align*}
e_{\alpha}\left(  z,K\right)   &  =\left\vert \mu\left(  z\right)  -\sum
_{j=0}^{K}\left(  \sum_{m=j+1}^{K}\frac{p_{j}\left(  m\right)  }{\prod
_{k=0}^{m-1}(z+\alpha+k)}\right)  \alpha^{j}\right\vert \\
&  =\left\vert \mu\left(  z\right)  -\frac{1}{12\left(  z+\alpha\right)  }%
\sum_{m=1}^{K}\prod_{j=1}^{m-1}\frac{\alpha-\xi_{j}\left(  m\right)
}{z+\alpha+j}\right\vert
\end{align*}
Given $z$, the last form is related to a Pad\'{e} approximant of order
$\left[  m-1/m+1\right]  $ in $\alpha$.

If $K=K^{\ast}\left(  z\right)  $, then we can find a value $\alpha^{\ast}$
that has a comparable error $\left\vert e_{\alpha}\left(  K\right)
\right\vert $ than Stirling's asymptotic approximation. For $z=10$ and
$\alpha\in\left\{  \frac{94909394316015843}{10^{18}},\frac{94909394316015845}%
{10^{18}}\right\}  $, the error $\log_{10}\left\vert e_{\alpha}\left(
10,K^{\ast}\left(  10\right)  \right)  \right\vert <-30$, while Stirling's
lowest error is $\log_{10}\left\vert e_{\text{Stirling}}\left(  10,K^{\ast
}\left(  10\right)  \right)  \right\vert =-28.5834$. The myth that Stirling's
approximation is \emph{always} better than any factorial, convergent expansion
for $\mu\left(  z\right)  $ with a \emph{same number of terms evaluated} is
not true, as illustrated by this counter example. If $\alpha^{\ast}$ is
approximated by a rational number, then all coefficients of $b_{m}\left(
\alpha\right)  $ are rational numbers, just as the Bernoulli numbers in the
Stirling approximation.

The major advantage of Binet's series (\ref{Binet_algemeen_explicit}) over
Stirling's asymptotic (\ref{Binet_asymptotic_Stirling_expansion}) lies in its
convergence, for all $\operatorname{Re}\left(  z\right)  >0$ and
$\operatorname{Re}\left(  \alpha\right)  >-\operatorname{Re}\left(  z\right)
$, towards $\mu\left(  z\right)  $, which allows incorporation into integrals
and series and may lead to other sharp bounds and approximations. Moreover,
the \textquotedblleft free\textquotedblright\ parameter $\alpha$ can be tuned
to achieve a similar accuracy as the best accuracy of Stirling's asymptotic
series (\ref{Stirling_explicit}). Finally, the coefficients $\left\vert
b_{m}\left(  \alpha\right)  \right\vert $ are monotonously increasing in $m>2$
and, for $\alpha>\xi_{1}\left(  m\right)  $, $b_{m}\left(  \alpha\right)  $ is
positive, but only changes the sign once for $\xi_{2}\left(  m\right)
<\alpha<\xi_{1}\left(  m\right)  $. For finite $K$, it might be interesting to
know the smallest possible error $e_{\alpha}\left(  z,K\right)  $ after
optimization of $\alpha$.

\section{Factorial series for Laplace transforms}

\label{sec_factorial_series_Laplace_transforms}The Laplace transform of a real
function $f\left(  t\right)  $ is defined (see e.g. \cite{Titchmarshfourier},
\cite[Chapter VII]{Evgrafov_1965}, \cite{Widder})\ for complex $z$ as
\begin{equation}
\varphi(z)=\mathcal{L}\left[  f(t)\right]  =\int_{0}^{\infty}e^{-zt}f(t)dt
\label{def_pgf_continuousrv}%
\end{equation}
with the inverse transform,
\begin{equation}
f(t)=\mathcal{L}^{-1}\left[  \varphi(z)\right]  =\frac{1}{2\pi i}%
\int_{c-i\infty}^{c+i\infty}\varphi(z)e^{zt}dz \label{inverse_Laplace}%
\end{equation}
where $c$ is the smallest real value of $\operatorname{Re}(z)$ for which the
integral in (\ref{def_pgf_continuousrv}) converges. Many functions can be
defined by an integral as (\ref{def_pgf_continuousrv}) as well as the
probability generating function \cite[Sec. 2.3.3]{PVM_PAComplexNetsCUP} of a
continuous random variable. For example, Binet's integral
(\ref{Binet_mu_integral_Bernoulli_generating_function}) is a Laplace transform
(\ref{def_pgf_continuousrv}), where $\varphi(z)=\mu\left(  z\right)  $ and the
integrand $f\left(  t\right)  =\frac{1}{t}\left(  \frac{1}{e^{t}-1}-\frac
{1}{t}+\frac{1}{2}\right)  $ with $f\left(  0\right)  =\frac{1}{12}$.

Factorial series are hardly studied. By starting from Gudermann's series
(\ref{Gudermann_series_mu}), Jensen \cite[art. 14]{Jensen_Gronwall1916} has
demonstrated, without using integrals, that Binet's factorial expansions are
absolutely and uniformly convergent in some region of $z$. Temme has written a
literature overview \cite{Temme1967}, of which parts were incorporated by
Lauwerier in his book \cite[p. 33-45]{Lauwerier1974}, that devotes one chapter
to factorial series. Parts of the content of \cite{Temme1967} and \cite[p.
33-45]{Lauwerier1974} are here absorbed. Apart from reviewing the literature
more extensively than here, Delabaere and Rasoamanana \cite{Delabaere2007}
have presented similar results, but their method and exposition is rather
different. In this Section \ref{sec_factorial_series_Laplace_transforms}, we
generalize the idea of Binet's proof of Theorem
\ref{theorem_Binet_convergent_expansion} as far as possible.

\subsection{The analogon of Stirling's asymptotic series}

\label{sec_analogon_Stirling}If we assume that the Taylor series $f\left(
t\right)  =\sum_{k=0}^{\infty}f_{k}t^{k}$, with Taylor coefficients
$f_{k}=\frac{1}{k!}\left.  \frac{d^{k}f(t)}{dt^{k}}\right\vert _{t=0}$, has an
infinite radius of convergence, then the Laplace transform
(\ref{def_pgf_continuousrv}) can be expanded as%
\[
\varphi(z)=\int_{0}^{\infty}e^{-zt}f(t)dt=\sum_{k=0}^{\infty}f_{k}\int
_{0}^{\infty}e^{-zt}t^{k}dt
\]
into a Laurent series \cite[Sec. 2.7]{Titchmarshfunctions}%
\begin{equation}
\varphi(z)=\sum_{k=0}^{\infty}\frac{k!f_{k}}{z^{k+1}}
\label{Laplace_transform_inverse_powers_z}%
\end{equation}
An infinite radius of convergence implies that the Taylor coefficient $f_{k}$
converges faster to zero than any exponential function, i.e. $f_{k}=o\left(
e^{-ak}\right)  $ for any finite $a$ and $k\rightarrow\infty$, and that
$f\left(  t\right)  $ is an entire function. Most functions, however, are not
entire functions. If the requirement of an infinite radius of convergence is
ignored, then (\ref{Laplace_transform_inverse_powers_z}) may represent an
asymptotic series $\varphi(z)=\sum_{k=0}^{K}\frac{k!f_{k}}{z^{k}}$, which
diverges when $K\rightarrow\infty$. Indeed, repeated partial integration of
the Laplace transform (\ref{def_pgf_continuousrv}) yields%
\begin{equation}
\varphi(z)=\sum_{k=0}^{K}\frac{k!f_{k}}{z^{k}}+\frac{1}{z^{K+1}}\int
_{0}^{\infty}e^{-zt}f^{\left(  K+1\right)  }(t)dt
\label{Laplace_transform_finite_Laurent_series}%
\end{equation}
illustrating\footnote{By the generalized mean-value theorem \cite[p.
321]{Hardy_pure_math}, there exists a positive real $\theta$ for which%
\[
\int_{0}^{\infty}e^{-zt}f^{\left(  K+1\right)  }(t)dt=f^{\left(  K+1\right)
}(\theta)\int_{0}^{\infty}e^{-zt}dt=\frac{1}{z}f^{\left(  K+1\right)  }%
(\theta)
\]
}, for any $\operatorname{Re}\left(  z\right)  >c$ but $z\neq0$ and real
$t\geq0$, that $\lim_{K\rightarrow\infty}\frac{f^{\left(  K+1\right)  }%
(t)}{z^{K+1}}$ must vanish to obtain a convergent Laurent series
(\ref{Laplace_transform_inverse_powers_z}).

\subsection{Binet's method as a recipe in five steps}

\label{sec_1_factorial_series}We formalize Binet's proof of Theorem
\ref{theorem_Binet_convergent_expansion}.

\begin{enumerate}
\item \emph{Binet's substitution }$e^{-t}=1-u$\emph{ or }$t=-\log\left(
1-u\right)  $ in the integral (\ref{def_pgf_continuousrv}) yields%
\begin{equation}
\varphi(z)=\int_{0}^{1}\left(  1-u\right)  ^{z-1}f(-\log\left(  1-u\right)
)du \label{Laplace_transform_after_Binet_substitution}%
\end{equation}
which converges for $\operatorname{Re}\left(  z\right)  >c$. Binet's
substitution is rather unusual, certainly in the study of Laplace integrals.
In the early days, Euler represented the Gamma function in the form of
(\ref{Laplace_transform_after_Binet_substitution}), after letting $w=1-u$, as
$\Gamma\left(  z\right)  =\int_{0}^{1}\log^{z-1}\left(  \frac{1}{w}\right)
dw$ (see e.g. \cite[{art. [1]}]{Binet1839}). Lauwerier \cite[p. 33-35]%
{Lauwerier1974} and Temme \cite{Temme1967} explain the success of the
substitution $u=1-e^{-t}$ by comparing $u=re^{i\theta}$ confined to the unit
disk $\left\vert u-1\right\vert <1$ at $u_{0}=1$ and the map $t=-\log\left(
1-re^{i\theta}\right)  $. In particular, the circle $\left\vert u-1\right\vert
=1$ at $u_{0}=1$ with radius 1 maps into the curve
\begin{align*}
t  &  =-\log\left(  1-\cos\theta-i\sin\theta\right)  =-\log\left(
\sqrt{2\left(  1-\cos\theta\right)  }e^{-i\arccos\frac{1-\cos\theta}%
{\sqrt{2\left(  1-\cos\theta\right)  }}}\right) \\
&  =-\log\left(  2\sin\frac{\theta}{2}\right)  +i\arccos\left(  \sin
\frac{\theta}{2}\right)  =-\log\left(  2\sin\frac{\theta}{2}\right)  +i\left(
\frac{\pi-\theta}{2}\right)
\end{align*}
Hence, $\operatorname{Re}t=-\log\left(  2\sin\frac{\theta}{2}\right)  $ and
$\operatorname{Im}t=\frac{\pi-\theta}{2}$ with $\theta\in\left[
0,2\pi\right]  $. Elimination of $\theta=\pi-2\operatorname{Im}t$ yields%
\begin{equation}
\operatorname{Re}t=-\log\left(  2\cos\left(  \operatorname{Im}t\right)
\right)  \hspace{1cm}\text{for }-\frac{\pi}{2}\leq\operatorname{Im}t\leq
\frac{\pi}{2} \label{map_curve}%
\end{equation}
The curve (\ref{map_curve}) is symmetric around the $\operatorname{Re}t$-axis
due to $\cos\left(  \operatorname{Im}t\right)  =\cos\left(  -\operatorname{Im}%
t\right)  $ and $\cos\left(  \operatorname{Im}t\right)  \geq0$ confines the
curve (\ref{map_curve}) to the strip $-\frac{\pi}{2}\leq\operatorname{Im}%
t\leq\frac{\pi}{2}$, where the minimum occurs at $-\log2$ for
$\operatorname{Im}t=0$ and $\operatorname{Re}t$ grows boundlessly if
$\operatorname{Im}t\rightarrow\pm\frac{\pi}{2}$. Thus, the map $t=-\log\left(
1-u\right)  $ of the unit disk $\left\vert u-1\right\vert <1$ appears as the
interior $t$-region bounded by the curve (\ref{map_curve}). That $t$-region is
considerably broader than the unit disk, which accounts for better results in
the sense that the resulting factorial series in
(\ref{factorial_expansion_alpha=0}) below converges for more functions than
its corresponding Laurent series (\ref{Laplace_transform_inverse_powers_z}).

\item The second step involves the Taylor expansion of $f(-\log\left(
1-u\right)  )$ around $u_{0}=0$, where $\log\left(  1-u_{0}\right)  =0$. After
Binet's substitution $e^{-t}=1-u$, the Taylor series $f\left(  t\right)
=\sum_{k=0}^{\infty}f_{k}t^{k}$ becomes $f\left(  -\log\left(  1-u\right)
\right)  =\sum_{k=0}^{\infty}f_{k}\left(  -\log\left(  1-u\right)  \right)
^{k}$. Introducing the second generating function (\ref{genfunc_stirlinginf})
of the Stirling numbers and reversing the $m$- and $k$-sum leads to%
\begin{equation}
f(-\log\left(  1-u\right)  )=\sum_{m=0}^{\infty}\left(  \sum_{k=0}^{m}\left.
\frac{d^{k}f(t)}{dt^{k}}\right\vert _{t=0}\left(  -1\right)  ^{m-k}S_{m}%
^{(k)}\,\right)  \frac{u^{m}}{m!} \label{Taylor_series_f(-log(1-u))_rond_u0=0}%
\end{equation}
Clearly, the Stirling numbers, which are integers, play a key role in Binet's
transformation $u=1-e^{-t}$.

\item Crucially for the third step, we assume that the Taylor series
(\ref{Taylor_series_f(-log(1-u))_rond_u0=0}) converges for $\left\vert
u\right\vert <1$. Hence, the radius of convergence of the Taylor series
(\ref{Taylor_series_f(-log(1-u))_rond_u0=0}) should be at least equal to one.
After substitution of the Taylor series
(\ref{Taylor_series_f(-log(1-u))_rond_u0=0}) in the integral
(\ref{Laplace_transform_after_Binet_substitution}) and reversing the summation
and integration operator, justified because a Taylor series can be term-wise
integrated within its radius of convergence, yields%
\[
\varphi(z)=\sum_{m=0}^{\infty}\left(  \sum_{k=0}^{m}\left.  \frac{d^{k}%
f(t)}{dt^{k}}\right\vert _{t=0}\left(  -1\right)  ^{m-k}S_{m}^{(k)}\,\right)
\frac{1}{m!}\int_{0}^{1}u^{m}\left(  1-u\right)  ^{z-1}du
\]

\item The fourth step uses the Beta integral $\int_{0}^{1}u^{p-1}\left(
1-u\right)  ^{q-1}du=\frac{\Gamma\left(  p\right)  \Gamma\left(  q\right)
}{\Gamma\left(  p+q\right)  }$, valid for $\operatorname{Re}\left(  p\right)
>0$ and $\operatorname{Re}\left(  q\right)  >0$, and leads to%
\[
\varphi(z)=\sum_{m=0}^{\infty}\left(  \sum_{k=0}^{m}\left.  \frac{d^{k}%
f(t)}{dt^{k}}\right\vert _{t=0}\left(  -1\right)  ^{m-k}S_{m}^{(k)}\,\right)
\frac{\Gamma\left(  z\right)  }{\Gamma\left(  m+1+z\right)  }%
\]

\item The fifth and final step replaces $\frac{\Gamma\left(  z+m\right)
}{\Gamma\left(  z\right)  }=\prod_{k=0}^{m-1}(z+k)$ and we arrive at the
factorial expansion, valid for $\operatorname{Re}\left(  z\right)  >c$,%
\begin{equation}
\varphi(z)=\sum_{m=0}^{\infty}\frac{\sum_{k=0}^{m}\left.  \frac{d^{k}%
f(t)}{dt^{k}}\right\vert _{t=0}\left(  -1\right)  ^{m-k}S_{m}^{(k)}}%
{\prod_{k=0}^{m}(z+k)} \label{factorial_expansion_alpha=0}%
\end{equation}

\end{enumerate}

The factorial expansion (\ref{factorial_expansion_alpha=0}) generalizes
Binet's second factorial expansion (\ref{Binet_function_faculteits_expansie})
in Theorem \ref{theorem_Binet_convergent_expansion}.

\subsection{Infinitely many factorial series for $\varphi(z)$}

\label{sec_infinitely_many_factorial_series_Laplance_transforms}We generalize
the recipe with five steps in Section \ref{sec_1_factorial_series}, as we did
for the particular case $\varphi(z)=\mu\left(  z\right)  $ in Theorem
\ref{theorem_mu_infinite_convergent_factorial_expansions}, together with an
additional $\beta$-scaling inspired by Temme \cite[p. 11]{Temme1967}. Our main
result is: \hfill

\begin{theorem}
\label{theorem_Laplace_transform_infinite_convergent_factorial_expansions}Only
if the Taylor series%
\begin{equation}
f(-\beta\log\left(  1-u\right)  )=\sum_{m=0}^{\infty}\left(  \sum_{k=0}%
^{m}k!f_{k}\beta^{k}\left(  -1\right)  ^{m-k}S_{m}^{(k)}\,\right)  \frac
{u^{m}}{m!} \label{Taylor_series_f(-beta*log(1-u))_rond_u0=0}%
\end{equation}
has a radius of convergence at least equal to 1, then the Laplace transform
$\varphi(z)$ of the function $f\left(  t\right)  $ possesses infinitely many
factorial series in the complex parameter $\alpha$ and real $\beta>0$, for
$\operatorname{Re}\left(  z\right)  >c$ and $\frac{\operatorname{Re}\left(
\alpha\right)  }{\beta}>-\operatorname{Re}\left(  z\right)  $,%
\begin{equation}
\varphi(z)=\beta\sum_{m=0}^{\infty}\frac{m!\phi_{m}\left(  \alpha
,\beta\right)  }{\prod_{k=0}^{m}(\beta z+\alpha+k)}
\label{factorial_series_Laplace_transform_in_alpha_en_beta}%
\end{equation}
where%
\begin{equation}
\phi_{m}\left(  \alpha,\beta\right)  =\frac{1}{m!}\sum_{k=0}^{m}\left(
\sum_{j=0}^{m-k}\left(  k+j\right)  !S_{m}^{(k+j)}\left(  -1\right)
^{m-\left(  k+j\right)  }f_{j}\beta^{j}\right)  \frac{\alpha^{k}}{k!}
\label{factorial_polynomial_in_alpha_en_beta}%
\end{equation}
Further, $\phi_{m}\left(  \alpha,\beta\right)  $ is a polynomial of degree $m$
in $\alpha$ with highest order term $\frac{f(0)}{m!}\alpha^{m}$~and $\phi
_{0}\left(  \alpha,\beta\right)  =f(0)$. A complex integral for $\phi
_{m}\left(  \alpha,\beta\right)  $, where the contour $C\left(  0\right)  $
encloses the point $\omega_{0}=0$, is%
\begin{equation}
\phi_{m}\left(  \alpha,\beta\right)  =\frac{\left(  -1\right)  ^{m}}{2\pi
i}\int_{C(0)}\frac{f(-\beta\log\left(  1+\omega\right)  )d\omega}{\left(
1+\omega\right)  ^{\alpha}\omega^{m+1}}
\label{factorial_polynomial_in_alpha_complex_integral}%
\end{equation}
while another compact form is%
\begin{equation}
\phi_{m}\left(  \alpha,\beta\right)  =\frac{1}{m!}\sum_{j=0}^{m}\left.
\frac{d^{j}}{dt^{j}}\left(  f(t)e^{\frac{\alpha}{\beta}t}\right)  \right\vert
_{t=0}S_{m}^{(j)}\left(  -1\right)  ^{m-j}\beta^{j}
\label{factorial_polynomial_in_alpha_en_beta_Leibniz}%
\end{equation}

\end{theorem}

Clearly, the factorial series in
(\ref{factorial_series_Laplace_transform_in_alpha_en_beta}) with
(\ref{factorial_polynomial_in_alpha_en_beta_Leibniz}) reduces to
(\ref{factorial_expansion_alpha=0}) for $\alpha=0$ and $\beta=1$.

\textbf{Proof}: We repeat the five steps in Section
\ref{sec_1_factorial_series}, but we rewrite the Laplace transform
(\ref{Laplace_transform_after_Binet_substitution}) of $\varphi(z)$ after a
generalization of Binet's substitution $e^{-t}=\left(  1-u\right)  ^{\beta}$
or $t=-\beta\log\left(  1-u\right)  $ as%
\[
\varphi(z)=\beta\int_{0}^{1}\left(  1-u\right)  ^{\beta z+\alpha-1}\left(
1-u\right)  ^{-\alpha}f(-\beta\log\left(  1-u\right)  )du
\]
which still converges for $\operatorname{Re}\left(  z\right)  >c$, in spite of
the introduction of the \textquotedblleft free\textquotedblright\ parameter
$\alpha$ and the real $\beta>0$.

The second step now involves the Taylor expansion of
\begin{equation}
g\left(  u;\alpha,\beta\right)  =\left(  1-u\right)  ^{-\alpha}f(-\beta
\log\left(  1-u\right)  )=\sum_{m=0}^{\infty}\phi_{m}\left(  \alpha
,\beta\right)  u^{m} \label{Taylor_series_g_alpa(u)}%
\end{equation}
around $u_{0}=0$. The Taylor series $f(-\beta\log\left(  1-u\right)  )$ in
(\ref{Taylor_series_f(-beta*log(1-u))_rond_u0=0}) follows from
(\ref{Taylor_series_f(-log(1-u))_rond_u0=0}). The Taylor series $\left(
1-u\right)  ^{-\alpha}=\sum_{k=0}^{\infty}\binom{-\alpha}{k}\left(  -1\right)
^{k}u^{k}$ is valid for any complex $\alpha$ provided $\left\vert u\right\vert
<1$. Thus, the radius of convergence of (\ref{Taylor_series_g_alpa(u)}) is
limited to 1 by $\left(  1-u\right)  ^{-\alpha}$, which is just sufficient for
the reversal of summation and integration in step three, provided the radius
of convergence of $g\left(  u;0,\beta\right)  =f(-\beta\log\left(  1-u\right)
)$ in (\ref{Taylor_series_f(-log(1-u))_rond_u0=0}) is at least equal to 1.
From Cauchy's integral theorem \cite{Titchmarshfunctions} we directly find the
integral representation (\ref{factorial_polynomial_in_alpha_complex_integral}%
). In addition, the Taylor coefficient $\phi_{m}\left(  \alpha,\beta\right)  $
in (\ref{Taylor_series_g_alpa(u)}) follows from the Cauchy product%
\begin{equation}
\phi_{m}\left(  \alpha,\beta\right)  =\sum_{j=0}^{m}\binom{-\alpha}{m-j}%
\frac{\left(  -1\right)  ^{m-j}}{j!}\sum_{l=0}^{j}l!f_{l}\beta^{l}\left(
-1\right)  ^{j-l}S_{j}^{(l)} \label{Taylor_coefficient_phi_m(alpha,beta)}%
\end{equation}

The remaining steps in Binet's method of Section \ref{sec_1_factorial_series},
consisting of the substitution of the Taylor series
(\ref{Taylor_series_g_alpa(u)}) of $g\left(  u;\alpha,\beta\right)  $ in the
integral, the reversal of integral and summation and the explicit evaluation
of the remaining Beta integral,
\[
\varphi(z)=\beta\sum_{m=0}^{\infty}\phi_{m}\left(  \alpha,\beta\right)
\int_{0}^{1}\left(  1-u\right)  ^{\beta z+\alpha-1}u^{m}du=\beta\sum
_{m=0}^{\infty}\phi_{m}\left(  \alpha,\beta\right)  \frac{m!\Gamma\left(
\beta z+\alpha\right)  }{\Gamma\left(  \beta z+\alpha+m+1\right)  }%
\]
lead to the factorial expansion
(\ref{factorial_series_Laplace_transform_in_alpha_en_beta}), valid for
$\operatorname{Re}\left(  z\right)  >c$ and $\frac{\operatorname{Re}\left(
\alpha\right)  }{\beta}>-\operatorname{Re}\left(  z\right)  $.

The remainder of the proof consists of simplifying the Taylor coefficient
$\phi_{m}\left(  \alpha,\beta\right)  $ in
(\ref{Taylor_coefficient_phi_m(alpha,beta)}). It is convenient to reverse the
summations,%
\[
\phi_{m}\left(  \alpha,\beta\right)  =\sum_{l=0}^{m}l!f_{l}\beta^{l}\sum
_{j=l}^{m}\frac{\Gamma\left(  1-\alpha\right)  }{\Gamma\left(  1-\alpha
-m+j\right)  }\frac{\left(  -1\right)  ^{m-l}}{\left(  m-j\right)  !j!}%
S_{j}^{(l)}%
\]
We introduce the generating function (\ref{genfunc_stirling}) of the Stirling
numbers $S_{m}^{(k)}$,%
\[
\phi_{m}\left(  \alpha,\beta\right)  =\frac{1}{m!}\sum_{l=0}^{m}l!f_{l}%
\beta^{l}\sum_{j=l}^{m}\binom{m}{j}\left(  -1\right)  ^{m-l}S_{j}^{(l)}%
\sum_{k=0}^{m-j}S_{m-j}^{(k)}\left(  -1\right)  ^{k}\;\alpha^{k}%
\]
Let $q=m-j$ in the double sum%
\[
T=\sum_{j=l}^{m}\binom{m}{j}\left(  -1\right)  ^{m-l}S_{j}^{(l)}\sum
_{k=0}^{m-j}S_{m-j}^{(k)}\;\left(  -1\right)  ^{k}\alpha^{k}=\left(
-1\right)  ^{m-l}\sum_{q=0}^{m-l}\binom{m}{q}S_{m-q}^{(l)}\sum_{k=0}^{q}%
S_{q}^{(k)}\;\left(  -1\right)  ^{k}\alpha^{k}%
\]
and after reversing the sums, we obtain $T=\left(  -1\right)  ^{m-l}\sum
_{k=0}^{m-l}\left(  \sum_{q=k}^{m-l}\binom{m}{q}S_{q}^{(k)}S_{m-q}%
^{(l)}\right)  \left(  -1\right)  ^{k}\alpha^{k}$. Using
(\ref{Stirling_binomial_sum_prod_two_Stirling}) and $S_{q}^{(m)}=0$ if $m>q$
gives $\sum_{q=k}^{m-l}\binom{m}{q}S_{q}^{(k)}S_{m-q}^{(l)}=\sum_{q=0}%
^{m}\binom{m}{q}S_{q}^{(k)}S_{m-q}^{(l)}=\frac{\left(  k+l\right)  !}%
{l!k!}S_{m}^{(k+l)}$, resulting in%
\[
T=\frac{\left(  -1\right)  ^{m-l}}{l!}\sum_{k=0}^{m-l}\frac{\left(
k+l\right)  !}{k!}S_{m}^{(k+l)}\left(  -1\right)  ^{k}\alpha^{k}%
\]
Hence, the Taylor coefficient (\ref{Taylor_coefficient_phi_m(alpha,beta)})
becomes%
\begin{equation}
\phi_{m}\left(  \alpha,\beta\right)  =\frac{1}{m!}\sum_{l=0}^{m}f_{l}\beta
^{l}\sum_{k=0}^{m-l}\frac{\left(  k+l\right)  !}{k!}S_{m}^{(k+l)}\left(
-1\right)  ^{m-\left(  l+k\right)  }\alpha^{k}
\label{factorial_polynomial_in_alpha_reversed_for_Taylor_coefficients}%
\end{equation}
which we can express as a polynomial in $\alpha$ by letting $n=m-l$%
\[
\phi_{m}\left(  \alpha,\beta\right)  =\frac{1}{m!}\sum_{n=0}^{m}f_{m-n}%
\beta^{m-n}\sum_{k=0}^{n}\frac{\left(  k+m-n\right)  !}{k!}S_{m}%
^{(k+m-n)}\left(  -1\right)  ^{n+k}\alpha^{k}%
\]
After reversing the sums and letting $j=m-n$, we arrive at
(\ref{factorial_polynomial_in_alpha_en_beta}). Reversing the sums in
(\ref{factorial_polynomial_in_alpha_en_beta})%
\[
\phi_{m}\left(  \alpha,\beta\right)  =\frac{1}{m!}\sum_{j=0}^{m}\left(
\sum_{k=0}^{j}\binom{j}{k}\left.  \frac{d^{j-k}f(t)}{dt^{j-k}}\right\vert
_{t=0}\left(  \frac{\alpha}{\beta}\right)  ^{k}\right)  S_{m}^{(j)}\left(
-1\right)  ^{m-j}\beta^{j}%
\]
substituting $\left(  \frac{\alpha}{\beta}\right)  ^{k}=\left.  \frac{d^{k}%
}{dt^{k}}e^{\frac{\alpha}{\beta}t}\right\vert _{t=0}$ into the $k$-sum in
brackets%
\[
\sum_{k=0}^{j}\binom{j}{k}\left.  \frac{d^{j-k}f(t)}{dt^{j-k}}\right\vert
_{t=0}\left(  \frac{\alpha}{\beta}\right)  ^{k}=\sum_{k=0}^{j}\binom{j}%
{k}\left.  \frac{d^{j-k}f(t)}{dt^{j-k}}\right\vert _{t=0}\left.  \frac{d^{k}%
}{dt^{k}}e^{\frac{\alpha}{\beta}t}\right\vert _{t=0}=\left.  \frac{d^{j}%
}{dt^{j}}\left(  f(t)e^{\frac{\alpha}{\beta}t}\right)  \right\vert _{t=0}%
\]
where Leibniz's rule has been used, demonstrates
(\ref{factorial_polynomial_in_alpha_en_beta_Leibniz}) and proves Theorem
\ref{theorem_Laplace_transform_infinite_convergent_factorial_expansions}%
.\hfill$\square\medskip$

A verification of Theorem
\ref{theorem_Laplace_transform_infinite_convergent_factorial_expansions} by
computing the inverse Laplace transform is given in Appendix
\ref{sec_inverse_Laplace_transform_factorial_series}. Lauwerier \cite[p.
35]{Lauwerier1974} interestingly mentions that, due to the asymptotic relation
$\frac{\Gamma\left(  m\right)  }{\Gamma\left(  z+m\right)  }\sim m^{-z}$ for
large $m$, the factorial series
(\ref{factorial_series_Laplace_transform_in_alpha_en_beta}) for $\beta=1$,
rewritten as%
\[
\varphi(z)=\Gamma\left(  z+\alpha\right)  \sum_{m=1}^{\infty}\frac
{\Gamma\left(  m\right)  }{\Gamma\left(  z+\alpha+m\right)  }\phi_{m-1}\left(
\alpha,1\right)
\]
and the corresponding Dirichlet series $\varphi_{D}(z)=\sum_{m=1}^{\infty
}\frac{m^{\alpha}\phi_{m-1}\left(  \alpha,1\right)  }{m^{-z}}$ have the same
converge range for $\operatorname{Re}\left(  z\right)  >c$. Consequently, the
rich theory of Dirichlet series (see e.g.
\cite{Titchmarshfunctions,Titchmarshzeta}) directly applies to convergence
aspect of the factorial series
(\ref{factorial_series_Laplace_transform_in_alpha_en_beta}). Reviewing
Landau's work on the factorial series, Temme \cite{Temme1967} adds Newton's
series $\sum_{m=0}^{\infty}\left(  -1\right)  ^{m}\binom{z-1}{m}\phi
_{m-1}\left(  \alpha,1\right)  $ to the factorial and Dirichlet series as the
third type of series with the same convergence range.

Lauwerier \cite[p. 42-43]{Lauwerier1974} gives examples of factorial series
(\ref{factorial_series_Laplace_transform_in_alpha_en_beta}) with $\alpha=0$
and $\beta=1$. Temme \cite{Temme1967} provides even more interesting examples
such as $\varphi(z)=\int_{0}^{\infty}\frac{e^{-zt}dt}{\left(  1+t\right)
^{\nu}}$. Here, we add:

\medskip\textbf{Example 1} If $f\left(  t\right)  =e^{bt}$, then
$\mathcal{L}\left[  e^{bt}\right]  =\frac{1}{z-b}$ and the corresponding
factorial polynomial $\phi_{m}\left(  \alpha,1\right)  $ in
(\ref{factorial_polynomial_in_alpha_en_beta_Leibniz}) is, with $\left.
\frac{d^{k}}{dt^{k}}\left(  f(t)e^{\alpha t}\right)  \right\vert
_{t=0}=\left.  \frac{d^{k}}{dt^{k}}\left(  e^{\left(  \alpha+b\right)
t}\right)  \right\vert _{t=0}=\left(  \alpha+b\right)  ^{k}$ and the
generating function (\ref{genfunc_stirling}),%
\[
m!\left.  \phi_{m}\left(  \alpha,1\right)  \right\vert _{e^{bt}}=\sum
_{k=0}^{m}\left(  \alpha+b\right)  ^{k}\left(  -1\right)  ^{m-k}S_{m}%
^{(k)}=\prod_{k=0}^{m-1}(\alpha+b+k)
\]
The factorial series
(\ref{factorial_series_Laplace_transform_in_alpha_en_beta}) becomes, for
$\operatorname{Re}\left(  z\right)  >b$,%
\begin{equation}
\frac{1}{z-b}=\sum_{m=0}^{\infty}\frac{\prod_{k=0}^{m-1}(b+\alpha+k)}%
{\prod_{k=0}^{m}(z+\alpha+k)}=\frac{\Gamma\left(  z+\alpha\right)  }%
{\Gamma\left(  b+\alpha\right)  }\sum_{m=0}^{\infty}\frac{\Gamma\left(
b+\alpha+m\right)  }{\Gamma\left(  z+\alpha+m+1\right)  }
\label{1_op_(z-b)_factorial_series}%
\end{equation}
which is known for $\alpha=0$ (see e.g. Nielsen \cite[band I, p.
77]{NielsenChelsea}). Indeed, Gauss's classical result \cite[15.1.20]%
{Abramowitz} for the hypergeometric series at $z=1$ is, for $c\neq-k$ ($k$
integer) and Re$\left(  c-a-b\right)  >0$,%
\begin{equation}
F\left(  a,b;c;1\right)  =\frac{\Gamma(c)\Gamma(c-a-b)}{\Gamma(c-a)\Gamma
(c-b)}=\frac{\Gamma(c)}{\Gamma(a)\Gamma(b)}\sum_{m=0}^{\infty}\frac
{\Gamma(a+m)\Gamma(b+m)}{\Gamma(c+m)m!} \label{hypergeometric_z=1_Gauss}%
\end{equation}
For $a=1$, $b\rightarrow b+\alpha$ and $c=z+\alpha+1$, Gauss's formula
(\ref{hypergeometric_z=1_Gauss}) reduces to (\ref{1_op_(z-b)_factorial_series}).

\medskip\textbf{Example 2} Let $f\left(  t\right)  =E_{a,b}\left(  t\right)
=\sum_{k=0}^{\infty}\frac{t^{k}}{\Gamma\left(  b+ak\right)  }$, which is the
Mittag-Leffler function \cite{Gorenflo_2020}. The Laplace transform (see e.g.
\cite[art. 20]{PVM_Mittag_Leffler}) is%
\begin{equation}
\int_{0}^{\infty}e^{-zt}t^{\gamma-1}E_{a,b}\left(  xt^{\beta}\right)
dt=\sum_{k=0}^{\infty}\frac{\Gamma\left(  \gamma+\beta k\right)  }%
{\Gamma\left(  b+ak\right)  }\frac{x^{k}}{z^{\gamma+k\beta}}
\label{Laplace_transform_general}%
\end{equation}
valid for $\left\vert \beta\right\vert \leq\left\vert a\right\vert $ and
$\gamma>0$. For $\gamma=\beta=x=1$, the Laplace transform
(\ref{Laplace_transform_general}) simplifies to a Laurent series
(\ref{Laplace_transform_inverse_powers_z}) in $z$. The corresponding factorial
polynomial (\ref{factorial_polynomial_in_alpha_en_beta}) can be written as
\[
m!\left.  \phi_{m}\left(  \alpha,\beta\right)  \right\vert _{E_{a,b}\left(
t\right)  }=\sum_{j=0}^{m}\left(  \sum_{k=0}^{j}\frac{j!}{\Gamma\left(
b+a\left(  j-k\right)  \right)  k!}\left(  \frac{\alpha}{\beta}\right)
^{k}\right)  S_{m}^{(j)}\left(  -1\right)  ^{m-j}\beta^{j}%
\]
where the $k$-sum reduces to $\left(  1+\frac{\alpha}{\beta}\right)  ^{j}$ for
$a=b=1$, in which case $E_{1,1}\left(  t\right)  =e^{t}$ simplifying to
example 1. Unfortunately, for arbitrary $a$ and $b$, we could not simply
$\left.  \phi_{m}\left(  \alpha,\beta\right)  \right\vert _{E_{a,b}\left(
t\right)  }$ for the Mittag-Leffler function $E_{a,b}\left(  t\right)  $.

On the other hand, after replacing $z\rightarrow kz$ in factorial series
(\ref{factorial_series_Laplace_transform_in_alpha_en_beta}), multiplying both
sides by $x^{k}$ and adding over all $k\geq0$, we formally obtain a
Mittag-Leffler transformation%
\[
\sum_{k=0}^{\infty}\frac{x^{k}\varphi(zk)}{\Gamma\left(  \beta zk+\alpha
\right)  }=\int_{0}^{\infty}f(t)E_{\beta z,\alpha}\left(  xe^{-zt}\right)
dt=\beta\sum_{m=0}^{\infty}\phi_{m}\left(  \alpha,\beta\right)  m!E_{\beta
z,\alpha+m+1}\left(  x\right)
\]

\subsection{Open question}

\label{sec_considerations}The factorial series
(\ref{factorial_series_Laplace_transform_in_alpha_en_beta}) of a non-entire
function may converge, whereas the Laurent series
(\ref{Laplace_transform_inverse_powers_z}) does not. An example is the Binet
function $\mu\left(  z\right)  $, whose Laurent series
(\ref{Laplace_transform_inverse_powers_z}) is Stirling's famous asymptotic,
but divergent series (\ref{Binet_asymptotic_Stirling_expansion}). Hence, the
question arises: \textquotedblleft Given the Taylor series $f\left(  t\right)
=\sum_{k=0}^{\infty}f_{k}t^{k}$ with radius of convergence $R_{f}$, when does
a factorial series (\ref{factorial_series_Laplace_transform_in_alpha_en_beta})
of the Laplace transform $\varphi(z)$ converge?\textquotedblright

We can only give a partial insight. The recipe in 5 steps for a factorial
series (\ref{factorial_series_Laplace_transform_in_alpha_en_beta}) of the
Laplace transform $\varphi(z)$ requires that the Taylor series
(\ref{Taylor_series_f(-beta*log(1-u))_rond_u0=0}) of $g\left(  u;0,\beta
\right)  =f(-\beta\log\left(  1-u\right)  )$ around the origin $u_{0}=0$
converges within the unit circle. Its corresponding Taylor coefficient is
written in terms of $f_{k}=\frac{1}{k!}\left.  \frac{d^{k}f(t)}{dt^{k}%
}\right\vert _{t=0}$ as $g_{0}=f_{0}$ and%
\begin{equation}
g_{m}\left(  \beta\right)  =\frac{1}{m!}\sum_{k=1}^{m}k!f_{k}\beta^{k}\left(
-1\right)  ^{m-k}S_{m}^{(k)}\,\hspace{1cm}\text{for }m>0
\label{Taylor_coefficient_gm_as_fk}%
\end{equation}
from which it follows that, for $m>0$, the Taylor coefficient $g_{m}\left(
\beta\right)  $ is independent of $f_{0}=f\left(  0\right)  $ (like any
characteristic coefficient (\ref{def_s})). The inverse transform of
(\ref{Taylor_coefficient_gm_as_fk})%
\begin{equation}
\beta^{m}f_{m}=\frac{1}{m!}\sum_{k=1}^{m}k!g_{k}\left(  \beta\right)  \left(
-1\right)  ^{m-k}\mathcal{S}_{m}^{(k)}\hspace{1cm}\text{for }m>0
\label{Taylor_coefficient_fj_as_gm}%
\end{equation}
where $\mathcal{S}_{m}^{(k)}$ is the Stirling number of the Second Kind, is
deduced as in the proof of Property
\ref{property_generalized_Binet_Bernoulli_polynomials}.

In contrast to $S_{m}^{(k)}$, the Stirling numbers $\mathcal{S}_{m}^{(k)}$ are
non-negative. Thus, $\left(  -1\right)  ^{m-k}S_{m}^{(k)}>0$ in
(\ref{Taylor_coefficient_gm_as_fk}), whereas $\left(  -1\right)
^{k}\mathcal{S}_{m}^{(k)}$ in (\ref{Taylor_coefficient_fj_as_gm}) is
alternating with $k$. If $f_{k}$ is non-negative, then
(\ref{Taylor_coefficient_gm_as_fk}) indicates that also $g_{m}\left(
\beta\right)  $ is non-negative. Moreover, (\ref{Taylor_coefficient_gm_as_fk})
written as $g_{m}\left(  \beta\right)  =f_{m}+\frac{1}{m!}\sum_{k=1}%
^{m-1}k!f_{k}\left(  -1\right)  ^{m-k}S_{m}^{(k)}$ then shows that
$g_{m}\left(  \beta\right)  \geq f_{m}$. Consequently, the radius $R_{g}$ of
convergence of $g\left(  u;0,\beta\right)  =f(-\beta\log\left(  1-u\right)
)=\sum_{m=0}^{\infty}g_{m}\left(  \beta\right)  u^{m}$ is not larger than the
radius $R_{f}$ of convergence of $f\left(  t\right)  =\sum_{k=0}^{\infty}%
f_{k}t^{k}$. However, if $f_{k}$ is non-negative, then the non-entire function
$f\left(  t\right)  $ has a pole at a finite, real $t=R_{f}$ and its Laplace
integral $\varphi(z)$ in (\ref{def_pgf_continuousrv}) does not exist. On the
other hand, if $f_{k}=\left(  -1\right)  ^{k}\left\vert f_{k}\right\vert $ is
alternating, then $f\left(  -t\right)  $ has non-negative Taylor coefficients
so that $f\left(  t\right)  $ is decreasing in $t$. The Laplace integral
$\varphi(z)$ exists for decreasing functions $f\left(  t\right)  $. If
$g_{m}\left(  \beta\right)  =\left(  -1\right)  ^{m}\left\vert g_{m}\left(
\beta\right)  \right\vert $ is alternating, then
(\ref{Taylor_coefficient_fj_as_gm}) shows that $f_{m}=\left(  -1\right)
^{m}\left(  \left\vert g_{m}\left(  \beta\right)  \right\vert +\frac{1}%
{m!}\sum_{k=1}^{m-1}k!\left\vert g_{k}\left(  \beta\right)  \right\vert
\mathcal{S}_{m}^{(k)}\right)  $ is alternating and $\left\vert f_{m}%
\right\vert >\left\vert g_{m}\left(  \beta\right)  \right\vert $, implying
that the radius of convergence $R_{g}\geq R_{f}$.

In summary, if $f_{k}=\left(  -1\right)  ^{k}\left\vert f_{k}\right\vert $ is
alternating, then the factorial series
(\ref{factorial_series_Laplace_transform_in_alpha_en_beta}) has higher
probability of convergence than its corresponding Laurent series
(\ref{Laplace_transform_inverse_powers_z}).

\medskip\textbf{Acknowledgement} I am very grateful to the late R. B. Paris
and R. E. Kooij for pointing to a few misprints in an early version and to A.
Olde Daalhuis for pointing me to Nemes' paper \cite{Nemes_2013}. After
completing this paper, G. Nemes has informed me about the work of Shi \emph{et
al.} \cite{ShiX2006}, whose Theorem 1 is essentially the same as our Theorem
\ref{theorem_mu_infinite_convergent_factorial_expansions}, though their
coefficient for $b_{m}\left(  \alpha\right)  $ is only mentioned in integral
form, different from (\ref{generalized_Binet_numbers_alpha_integral}), and our
validity range of $\alpha$ is slightly broader. Moreover, our proof(s) and
approach are entirely different. A. Olde Daalhuis and G. Nemes pointed me to
the book of Lauwerier \cite{Lauwerier1974} and to \cite{Delabaere2007}. N.
Temme has sent me his earlier report \cite{Temme1967} of which parts appeared
in Lauwerier's book.

{\footnotesize
\bibliographystyle{plain}
\bibliography{cac,math,misc,net,pvm,qth,tel}

\begin{thebibliography}{10}

\bibitem{Abramowitz}
M.~Abramowitz and I.~A. Stegun.
\newblock {\em Handbook of Mathematical Functions}.
\newblock Dover Publications, Inc., New York, 1968.

\bibitem{Binet1839}
J.~P.~M. Binet.
\newblock M{\'e}moire sur les int{\'e}grales d{\'e}finites {E}ul{\'e}riennes et
  sur leur application {\`a} la th{\'e}orie des suites ainsi qu` {\`a} l\`
  {\'e}valuation des functions des grands nombres.
\newblock {\em Journal de l`\'{E}cole Polytechnique}, XVI:123--343, July 1839.

\bibitem{Blagouchine_2016}
I.~V. Blagouchine.
\newblock Two series expansions for the logarithm of the gamma function
  involving {S}tirling numbers and containing only rational coefficients for
  certain arguments related to $\pi^{-1}$.
\newblock {\em arXiv1408.3902v9}, May 2016.

\bibitem{Comtet}
L.~Comtet.
\newblock {\em Advanced Combinatorics}.
\newblock D. Riedel Publishing Company, Dordrecht, Holland, revised and
  enlarged edition, 1974.

\bibitem{Delabaere2007}
E.~Delabaere and J.-M. Rasoamanana.
\newblock Sommation effective d'une somme de {B}orel par s{\'e}ries de
  factorielles.
\newblock {\em Annales de l`institut Fourier}, 57(2):421--456, 2007.

\bibitem{Erdelyi_v1}
A.~Erd{\'e}lyi, W.~Magnus, F.~Oberhettinger, and F.~G. Tricomi.
\newblock {\em Higher Transcendental Functions}, volume~1 of {\em California
  Institute of Technology Bateman Manuscript Project}.
\newblock McGraw-Hill Book Company, New York, 1953.

\bibitem{Evgrafov_1965}
M.~A. Evgrafov.
\newblock {\em Analytic Functions}.
\newblock W. B. Saunders Company, 1966; Reprinted by Dover Publications, Inc.,
  New York, dover 2019 edition, 2019.

\bibitem{Gilbert1873}
Ph. Gilbert.
\newblock Recherches sur de d{\'e}veloppement de la function $\gamma$ et sur
  certaines int{\'e}grals d{\'e}finies qui en d{\'e}pendent.
\newblock {\em M{\'e}moires de l'Acad{\'e}mie royale des Sciences, des Lettres
  et des Beaux-Arts de Belgique}, 41:1--60, 1876.

\bibitem{Gilbert1886}
Ph. Gilbert.
\newblock Sur les produits compos{\'e}s d'un grand nombre de facteurs et sur le
  reste de la s{\'e}rie de binet.
\newblock {\em Annales de la Soci{\'e}t{\'e} scientifique de Bruxelles},
  10:191--200, 1886.

\bibitem{Gorenflo_2020}
R.~Gorenflo, A.~A. Kilbas, F.~Mainardi, and S.~V. Rogosin.
\newblock {\em Mittag-Leffler Functions, Related Topics and Applications}.
\newblock Springer, second edition, 2020.

\bibitem{Gudermann1845}
C.~Gudermann.
\newblock Additamentum ad functionis {$\Gamma$}$(a)=\int_{0}^{\infty}e^{-x}
  x^{a-1} dx$ theoriam.
\newblock {\em Journal f{\"u}r die reine und angewandte Mathematik},
  29:209--212, 1845.

\bibitem{Hardy_div}
G.~H. Hardy.
\newblock {\em Divergent Series}.
\newblock Oxford University Press, London, 1948.

\bibitem{Hardy_pure_math}
G.~H. Hardy.
\newblock {\em A Course of Pure Mathematics}.
\newblock Cambridge University Press, 10nth edition, 2006.

\bibitem{Hermite1895}
Ch. Hermite.
\newblock Sur la function {$\log\Gamma(a)$}.
\newblock {\em Journal f{\"u}r die reine und angewandte Mathematik},
  115:201--208, 1895.

\bibitem{Jensen_Gronwall1916}
J.~L. W.~V. Jensen and T.~H. Gronwall.
\newblock An elementary exposition of the theory of the {G}amma function.
\newblock {\em Annals of Mathematics}, 17(3):124--166, March 1916.

\bibitem{Lauwerier1974}
H.~A. Lauwerier.
\newblock {\em Aymptotic Analysis}.
\newblock Mathematical Centre Tracts. Mathematisch Centrum, Amsterdam, 1974.

\bibitem{Nemes_2013}
G.~Nemes.
\newblock Generalization of {B}inet's {G}amma function formulas.
\newblock {\em Integral Transforms and Special Functions}, 24(8):595--606,
  2013.

\bibitem{NielsenChelsea}
N.~Nielsen.
\newblock {\em Die Gammafunktion: Band I. Handbuch der Theorie der
  Gammafunktion und Band II. Theorie des Integrallogarithmus und verwandter
  Transzendenten}.
\newblock B. G. Teubner, Leipzig 1906; republished by Chelsea, New York, 1956.

\bibitem{Olver}
F.~W.~J. Olver, D.~W Lozier, R.~F. Boisvert, and C.~W. Clark.
\newblock {\em {NIST} {H}andbook of {M}athematical {F}unctions}.
\newblock Cambridge University Press, New York, 2010.

\bibitem{Paris_Kaminski}
R.~B. Paris and D.~Kaminski.
\newblock {\em Asymptotics and Mellin-Barnes Integrals}, volume~85 of {\em
  Encyclopedia of Mathematics and its Applications}.
\newblock Cambridge University Press, U. K., 2001.

\bibitem{Riordan}
J.~Riordan.
\newblock {\em Combinatorial Identities}.
\newblock John Wiley \& Sons, New York, 1968.

\bibitem{Sansone}
G.~Sansone and J.~Gerretsen.
\newblock {\em Lectures on the Theory of Functions of a Complex Variable},
  volume 1 and 2.
\newblock P. Noordhoff, Groningen, 1960.

\bibitem{ShiX2006}
X.~Shi, F.~Liu, and M.~Hu.
\newblock A new asymptotic series for the {G}amma function.
\newblock {\em Journal of Computational and Applied Mathematics}, 195:134--154,
  2006.

\bibitem{Temme1967}
N.~M. Temme.
\newblock Asymptotische ontwikkelingen van {F}akulteitsreeksen.
\newblock {\em Technical Report TN 48}, Mathematisch Centrum, Amsterdam, March
  1967.

\bibitem{Titchmarshfourier}
E.~C. Titchmarsh.
\newblock {\em Introduction to the Theory of Fourier Integrals}.
\newblock Oxford University Press, Ely House, London W. I, 2nd edition, 1948.

\bibitem{Titchmarshfunctions}
E.~C. Titchmarsh.
\newblock {\em The Theory of Functions}.
\newblock Oxford University Press, Amen House, London, 1964.

\bibitem{Titchmarshzeta}
E.~C. Titchmarsh and D.~R. Heath-Brown.
\newblock {\em The Theory of the Zeta-function}.
\newblock Oxford Science Publications, Oxford, 2nd edition, 1986.

\bibitem{PVM_ASYM}
P.~Van~Mieghem.
\newblock The asymptotic behaviour of queueing systems: Large deviations theory
  and dominant pole approximation.
\newblock {\em Queueing Systems}, 23:27--55, 1996.

\bibitem{PVM_graphspectra}
P.~Van~Mieghem.
\newblock {\em Graph Spectra for Complex Networks}.
\newblock Cambridge University Press, Cambridge, U.K., 2011.

\bibitem{PVM_PAComplexNetsCUP}
P.~Van~Mieghem.
\newblock {\em Performance Analysis of Complex Networks and Systems}.
\newblock Cambridge University Press, Cambridge, U.K., 2014.

\bibitem{PVM_Mittag_Leffler}
P.~Van~Mieghem.
\newblock The {M}ittag-{L}effler function.
\newblock {\em arXiv:2005.13330}, 2020.

\bibitem{Whittaker_Watson}
E.~T. Whittaker and G.~N. Watson.
\newblock {\em A Course of Modern Analysis}.
\newblock Cambridge University Press, Cambridge, UK, cambridge mathematical
  library edition, 1996.

\bibitem{Widder}
D.~V. Widder.
\newblock {\em The Laplace transform}.
\newblock Princeton University Press, Princeton, 1946.

\end{thebibliography}
}

\appendix

\section{Complex integral for Binet's function $\mu\left(  z\right)  $}

\label{sec_complex_integral_mu}

\subsection{Derivation of the complex integral in
(\ref{Binet_complex_integral_Zeta})}

From Weierstrass's product (\ref{Weierstrass_product_Gamma_function}) of the
Gamma function, Whittaker and Watson \cite[p. 277]{Whittaker_Watson} deduce
the formula, valid for all $a$ and $z$,%
\[
\log\frac{\Gamma\left(  a\right)  }{\Gamma\left(  z+a\right)  }=-z\frac
{\Gamma^{\prime}\left(  a\right)  }{\Gamma\left(  a\right)  }+\frac{1}{2\pi
i}\int_{q-i\infty}^{q+i\infty}\frac{\pi}{\sin\pi s}\frac{\zeta\left(
s,a\right)  }{s}z^{s}ds\hspace{1cm}\text{with }1<q<2
\]
where $\zeta\left(  s,a\right)  =\sum_{n=0}^{\infty}\frac{1}{\left(
n+a\right)  ^{s}}$ is the Hurwitz Zeta-function, which reduces for $a=1$ to
the Riemann Zeta-function $\zeta\left(  s\right)  $. Thus, for $a=1$ and
$\frac{\Gamma^{\prime}\left(  1\right)  }{\Gamma\left(  1\right)  }=-\gamma$,
where $\gamma$ is the Euler constant, we have%
\[
\log\Gamma\left(  z+1\right)  =-\gamma z-\frac{1}{2\pi i}\int_{q-i\infty
}^{q+i\infty}\frac{\pi}{\sin\pi s}\frac{\zeta\left(  s\right)  }{s}%
z^{s}ds\hspace{1cm}\text{with }1<q<2
\]
If we move the line of integration to $0<c=\operatorname{Re}\left(  s\right)
<1$, we encounter a double pole at $s=1$, because $\zeta\left(  s\right)
=\frac{1}{s-1}+\gamma+O\left(  s-1\right)  $ around $s=1$ and a zero of
$\sin\pi s$. The residue at $s=1$ follows from Cauchy's integral theorem
$\frac{1}{k!}\;\left.  \frac{d^{k}f(z)}{dz^{k}}\right\vert _{z=z_{0}}=\frac
{1}{2\pi i}\int_{C(z_{0})}\frac{f(\omega)\;d\omega}{(\omega-z_{0})^{k+1}}$,
where $f\left(  z\right)  $ is analytic within the contour $C\left(
z_{0}\right)  $ that encloses the point $z_{0}$ and we obtain%
\[
\frac{1}{2\pi i}\int_{q-i\infty}^{q+i\infty}\frac{\pi}{\sin\pi s}\frac
{\zeta\left(  s\right)  }{s}z^{s}ds=\frac{1}{2\pi i}\int_{c-i\infty
}^{c+i\infty}\frac{\pi}{\sin\pi s}\frac{\zeta\left(  s\right)  }{s}%
z^{s}ds+\lim_{s\rightarrow1}\frac{d}{ds}\left(  \frac{z^{s}}{s}\frac{\pi
\zeta\left(  s\right)  \left(  s-1\right)  ^{2}}{\sin\pi s}\right)
\]
because the function between brackets is analytic at $s=1$. Executing the
derivative,%
\[
\frac{d}{ds}\left(  \frac{z^{s}}{s}\frac{\pi\zeta\left(  s\right)  \left(
s-1\right)  ^{2}}{\sin\pi s}\right)  =\frac{z^{s}\pi\left(  s-1\right)
}{s\sin\pi s}\left(  \left(  \log z-\frac{1}{s}\right)  \zeta\left(  s\right)
\left(  s-1\right)  +\zeta^{\prime}\left(  s\right)  \left(  s-1\right)
+2\zeta\left(  s\right)  -\frac{\pi\zeta\left(  s\right)  \left(  s-1\right)
\cos\pi s}{\sin\pi s}\right)
\]
and using the Taylor expansions of $\left(  s-1\right)  \zeta\left(  s\right)
$ around $s=1$ gives us%
\[
\lim_{s\rightarrow1}\frac{d}{ds}\left(  \frac{z^{s}}{s}\frac{\pi\zeta\left(
s\right)  \left(  s-1\right)  ^{2}}{\sin\pi s}\right)  =-z\left(  \log
z-1+\gamma\right)
\]
and we obtain, for $0<c<1$,%
\[
\log\Gamma\left(  z+1\right)  =z\log z-z-\frac{1}{2\pi i}\int_{c-i\infty
}^{c+i\infty}\frac{\pi}{\sin\pi s}\frac{\zeta\left(  s\right)  }{s}z^{s}ds
\]
Moving the line of integration over the double pole at $s=0$ to the left
yields, for $-1<c^{\prime}<0$,%
\[
\log\Gamma\left(  z+1\right)  =z\log z-z-\frac{1}{2\pi i}\int_{c^{\prime
}-i\infty}^{c^{\prime}+i\infty}\frac{\pi}{\sin\pi s}\frac{\zeta\left(
s\right)  }{s}z^{s}ds-\lim_{s\rightarrow0}\frac{d}{ds}\left(  \frac{\pi
s\zeta\left(  s\right)  }{\sin\pi s}z^{s}\right)
\]
The derivative is%
\[
\frac{d}{ds}\left(  \frac{s\zeta\left(  s\right)  }{\sin\pi s}z^{s}\right)
=\frac{z^{s}\zeta\left(  s\right)  s}{\sin\pi s}\left\{  \frac{1}{s}-\frac
{\pi\cos\pi s}{\sin\pi s}+\log z+\frac{\zeta^{\prime}\left(  s\right)  }%
{\zeta\left(  s\right)  }\right\}
\]
Since $\pi\cot\left(  \pi x\right)  =\frac{1}{x}-2\sum_{n=1}^{\infty}%
\zeta\left(  2n\right)  \,x^{2n-1}$, we find that $\lim_{s\rightarrow0}%
\frac{1}{s}-\pi\cot\left(  \pi s\right)  =0$ and%
\[
\lim_{s\rightarrow0}\frac{ds}{ds}\left(  \frac{\pi s\zeta\left(  s\right)
}{\sin\pi s}z^{s}\right)  =\zeta\left(  0\right)  \left\{  \log z+\frac
{\zeta^{\prime}\left(  0\right)  }{\zeta\left(  0\right)  }\right\}
\]
With $\zeta\left(  0\right)  =-\frac{1}{2}$ and $\zeta^{\prime}\left(
0\right)  =-\frac{1}{2}\log(2\pi)$, we arrive at $\lim_{s\rightarrow0}%
\frac{ds}{ds}\left(  \frac{\pi s\zeta\left(  s\right)  }{\sin\pi s}%
z^{s}\right)  =-\frac{1}{2}\left\{  \log z+\log(2\pi)\right\}  $ and%
\[
\log\Gamma\left(  z+1\right)  =\left(  z+\frac{1}{2}\right)  \log z-z+\frac
{1}{2}\log(2\pi)-\frac{1}{2\pi i}\int_{c^{\prime}-i\infty}^{c^{\prime}%
+i\infty}\frac{\pi}{\sin\pi s}\frac{\zeta\left(  s\right)  }{s}z^{s}%
ds\hspace{1cm}\text{with }-1<c^{\prime}<0
\]
From the definition $\log\Gamma\left(  z\right)  =\left(  z-\frac{1}%
{2}\right)  \log z-z+\frac{1}{2}\log\left(  2\pi\right)  +\mu\left(  z\right)
$, we find (\ref{Binet_complex_integral_Zeta}).

We present a second, shorter derivation of (\ref{Binet_complex_integral_Zeta})
by employing the inverse Mellin transform
\[
\frac{1}{e^{2\pi t}-1}=\frac{1}{2\pi i}\int_{c-i\infty}^{c+i\infty}%
\Gamma\left(  s\right)  \zeta\left(  s\right)  \left(  2\pi t\right)
^{-s}ds\hfill\hspace{1cm}\text{for }c>1
\]
Substitution into Binet's integral (\ref{Binet_mu_integral_arctan}) and
reversing the integrals gives%
\[
\mu\left(  z\right)  =2\int_{0}^{\infty}\frac{\arctan\left(  \frac{t}%
{z}\right)  }{e^{2\pi t}-1}dt=\frac{1}{\pi i}\int_{c-i\infty}^{c+i\infty
}\Gamma\left(  s\right)  \zeta\left(  s\right)  \left(  2\pi\right)
^{-s}\left(  \int_{0}^{\infty}\arctan\left(  \frac{t}{z}\right)
t^{-s}dt\right)  ds
\]
Partial integration, followed by a substitution $u=\left(  \frac{t}{z}\right)
^{2}$ and the use of the Beta integral and the Gamma reflection formula
results in%
\[
-\int_{0}^{\infty}\arctan\left(  \frac{t}{z}\right)  t^{-s}dt=\frac{z^{1-s}%
\pi}{2\left(  1-s\right)  \sin\frac{\pi s}{2}}%
\]
and%
\[
\mu\left(  z\right)  =-\frac{z}{2i}\int_{c-i\infty}^{c+i\infty}\zeta\left(
s\right)  \frac{\Gamma\left(  s\right)  \left(  2\pi z\right)  ^{-s}}{\left(
1-s\right)  \sin\frac{\pi s}{2}}ds
\]
Using the functional equation $\zeta(s)=2(2\pi)^{s-1}\sin\frac{\pi s}{2}%
\Gamma(1-s)\zeta(1-s)$ yields%
\[
\mu\left(  z\right)  =-\frac{1}{2\pi i}\int_{c-i\infty}^{c+i\infty}\frac{\pi
}{\sin\pi s}\frac{\zeta(1-s)}{\left(  1-s\right)  }z^{1-s}ds\hspace
{1cm}\text{with }1<c<2
\]
and a change of variable $w=1-s$ then returns again the complex integral in
(\ref{Binet_complex_integral_Zeta}).

\subsection{Derivation of the convergent series
(\ref{mu_series_TaylorZeta_1_op_z})}

\label{sec_convergent_mu_series_Stirling}Substituting the Taylor series
$\left(  s-1\right)  \zeta\left(  s\right)  =\sum_{m=0}^{\infty}g_{m}\left(
1\right)  \left(  s-1\right)  ^{m}$ into the integral
(\ref{Binet_complex_integral_Zeta}) yields
\begin{align*}
\mu\left(  z\right)   &  =-\frac{1}{2\pi i}\int_{c-i\infty}^{c+i\infty}%
\frac{\pi}{\sin\pi s}\frac{\zeta\left(  s\right)  }{s}z^{s}ds\hspace
{1cm}\text{with }-1<c<0\\
&  =-\frac{1}{2\pi i}\int_{c-i\infty}^{c+i\infty}\frac{\pi}{\sin\pi s}\frac
{1}{s\left(  s-1\right)  }\sum_{m=0}^{\infty}g_{m}\left(  1\right)  \left(
s-1\right)  ^{m}z^{s}ds
\end{align*}
Integration and summation can be reversed, because the Taylor series converges
for all complex $s$ and within the radius of convergence, a Taylor series
represents an analytic function that can be integrated and differentiated
\cite[p. 97]{Titchmarshfunctions},%
\begin{align*}
\mu\left(  z\right)   &  =-\sum_{m=0}^{\infty}g_{m}\left(  1\right)  \frac
{1}{2\pi i}\int_{c-i\infty}^{c+i\infty}\frac{\pi}{\sin\pi s}\frac{\left(
s-1\right)  ^{m}z^{s}}{s\left(  s-1\right)  }ds\\
&  =-\frac{1}{2\pi i}\int_{c-i\infty}^{c+i\infty}\frac{\pi}{\sin\pi s}%
\frac{z^{s}}{s\left(  s-1\right)  }ds-\sum_{m=1}^{\infty}g_{m}\left(
1\right)  \frac{1}{2\pi i}\int_{c-i\infty}^{c+i\infty}\frac{\pi}{\sin\pi
s}\frac{\left(  s-1\right)  ^{m-1}z^{s}}{s}ds
\end{align*}
We evaluate the first integral. If $\left\vert z\right\vert \leq1$, then we
close the contour over the positive $\operatorname{Re}\left(  s\right)
$-plane (where the integral over semi-circle at infinity vanishes). Cauchy's
residu theorem tells us that%
\begin{align*}
-\frac{1}{2\pi i}\int_{c-i\infty}^{c+i\infty}\frac{\pi}{\sin\pi s}\frac{z^{s}%
}{s\left(  s-1\right)  }ds  &  =\pi\sum_{n=2}^{\infty}\lim_{s\rightarrow
n}\frac{\left(  s-n\right)  }{\sin\pi s}\frac{z^{s}}{s\left(  s-1\right)
}+\pi\lim_{s\rightarrow0}\frac{d}{ds}\frac{s}{\sin\pi s}\frac{z^{s}}{\left(
s-1\right)  }\\
&  +\pi\lim_{s\rightarrow1}\frac{d}{ds}\frac{\left(  s-1\right)  }{\sin\pi
s}\frac{z^{s}}{s}%
\end{align*}
With%
\[
\frac{d}{ds}\frac{s}{\sin\pi s}\frac{z^{s}}{\left(  s-1\right)  }=\frac
{s}{\sin\pi s}\frac{z^{s}}{\left(  s-1\right)  }\left\{  \frac{1}{s}-\pi
\cot\pi s+\log z-\frac{1}{\left(  s-1\right)  }\right\}
\]
Since $\lim_{s\rightarrow0}\frac{1}{s}-\pi\cot\pi s=0$, we find that%
\[
\pi\lim_{s\rightarrow0}\frac{d}{ds}\frac{s}{\sin\pi s}\frac{z^{s}}{\left(
s-1\right)  }=-\left(  \log z+1\right)  \lim_{s\rightarrow0}\frac{\pi s}%
{\sin\pi s}=-\left(  \log z+1\right)
\]
and, similarly, that%
\[
\pi\lim_{s\rightarrow1}\frac{d}{ds}\frac{\left(  s-1\right)  }{\sin\pi s}%
\frac{z^{s}}{s}=-z\left(  \log z-1\right)
\]
Hence, for $\left\vert z\right\vert \leq1$,%
\[
-\frac{1}{2\pi i}\int_{c-i\infty}^{c+i\infty}\frac{\pi}{\sin\pi s}\frac{z^{s}%
}{s\left(  s-1\right)  }ds=\sum_{n=2}^{\infty}\frac{\left(  -1\right)
^{n}z^{n}}{n\left(  n-1\right)  }-\left(  \log z+1\right)  -z\left(  \log
z-1\right)
\]
With $\frac{1}{n-1}-\frac{1}{n}=\frac{1}{n\left(  n-1\right)  }$, we have%
\begin{align*}
\sum_{n=2}^{\infty}\frac{\left(  -1\right)  ^{n}z^{n}}{n\left(  n-1\right)  }
&  =\sum_{n=2}^{\infty}\frac{\left(  -z\right)  ^{n}}{n-1}-\sum_{n=2}^{\infty
}\frac{\left(  -z\right)  ^{n}}{n}=\left(  -z-1\right)  \sum_{n=1}^{\infty
}\frac{\left(  -z\right)  ^{n}}{n}+\left(  -z\right) \\
&  =\left(  z+1\right)  \log\left(  1+z\right)  -z
\end{align*}
Thus, for $\left\vert z\right\vert \leq1$, we find%
\[
-\frac{1}{2\pi i}\int_{c-i\infty}^{c+i\infty}\frac{\pi}{\sin\pi s}\frac{z^{s}%
}{s\left(  s-1\right)  }ds=\left(  z+1\right)  \log\left(  1+\frac{1}%
{z}\right)  -1
\]
For $\left\vert z\right\vert \geq1$, we close the contour over the negative
$\operatorname{Re}\left(  s\right)  $-plane,%
\begin{align*}
-\frac{1}{2\pi i}\int_{c-i\infty}^{c+i\infty}\frac{\pi}{\sin\pi s}\frac{z^{s}%
}{s\left(  s-1\right)  }ds  &  =-\sum_{n=1}^{\infty}\lim_{s\rightarrow-n}%
\frac{\pi\left(  s+n\right)  }{\sin\pi s}\frac{z^{s}}{s\left(  s-1\right)
}=-\sum_{n=1}^{\infty}\frac{\left(  -1\right)  ^{n}z^{-n}}{n\left(
n+1\right)  }\\
&  =\left(  z+1\right)  \log\left(  1+\frac{1}{z}\right)  -1
\end{align*}
In summary, the first term equals%
\[
-\frac{1}{2\pi i}\int_{c-i\infty}^{c+i\infty}\frac{\pi}{\sin\pi s}\frac{z^{s}%
}{s\left(  s-1\right)  }ds=\left(  z+1\right)  \log\left(  1+\frac{1}%
{z}\right)  -1
\]
and%
\[
\mu\left(  z\right)  =\left(  z+1\right)  \log\left(  1+\frac{1}{z}\right)
-1-\sum_{m=1}^{\infty}g_{m}\left(  1\right)  \frac{1}{2\pi i}\int_{c-i\infty
}^{c+i\infty}\frac{\pi}{\sin\pi s}\frac{\left(  s-1\right)  ^{m-1}z^{s}}{s}ds
\]

The remaining integral is evaluated similarly. For $\left\vert z\right\vert
>1$, we close the contour over negative $\operatorname{Re}\left(  s\right)
$-plane and obtain%
\[
\frac{1}{2\pi i}\int_{c-i\infty}^{c+i\infty}\frac{\pi}{\sin\pi s}\frac{\left(
s-1\right)  ^{m-1}z^{s}}{s}ds=\sum_{n=1}^{\infty}\lim_{s\rightarrow-n}%
\frac{\pi\left(  s+n\right)  }{\sin\pi s}\frac{\left(  s-1\right)  ^{m-1}%
z^{s}}{s}=\sum_{n=1}^{\infty}\left(  -1\right)  ^{n}\frac{\left(  -n-1\right)
^{m-1}z^{-n}}{-n}%
\]
and%
\[
\frac{1}{2\pi i}\int_{c-i\infty}^{c+i\infty}\frac{\pi}{\sin\pi s}\frac{\left(
s-1\right)  ^{m-1}z^{s}}{s}ds=\left(  -1\right)  ^{m}\sum_{n=1}^{\infty}%
\frac{\left(  n+1\right)  ^{m-1}}{n}\frac{1}{\left(  -z\right)  ^{n}}%
\]
From $\log\left(  1+x\right)  =-\sum_{n=1}^{\infty}\frac{\left(  -x\right)
^{n}}{n}$, we have that%
\[
\left(  -1\right)  ^{m}\sum_{n=1}^{\infty}\frac{\left(  n+1\right)  ^{m-1}}%
{n}\frac{1}{\left(  -z\right)  ^{n}}=\left(  -1\right)  ^{m-1}\left.
e^{-y}\frac{d^{m-1}}{dy^{m-1}}\left(  e^{y}\log\left(  1+e^{y}\right)
\right)  \right\vert _{z=e^{-y}}%
\]
Leibniz' rule gives%
\begin{align*}
\frac{d^{m}}{dy^{m}}\left(  e^{y}\log\left(  1+e^{y}\right)  \right)   &
=\sum_{l=0}^{m}\binom{m}{l}\frac{d^{m-l}}{dy^{m-l}}\left(  e^{y}\right)
\frac{d^{l}}{dy^{l}}\left(  \log\left(  1+e^{y}\right)  \right) \\
&  =e^{y}\left(  \log\left(  1+e^{y}\right)  \right)  +e^{y}\sum_{l=1}%
^{m}\binom{m}{l}\frac{d^{l-1}}{dy^{l-1}}\left(  1+e^{-y}\right)  ^{-1}%
\end{align*}
For $k>0$, it holds\footnote{In the theory of the Fermi-Dirac integral
$F_{p}(z)=\frac{1}{\Gamma(p+1)}\int_{0}^{\infty}\frac{x^{p}}{1+e^{x-z}}\;dx$
for complex $p$ and $z$, the functional equation $\frac{dF_{p}(y)}{dy}%
=F_{p-1}(y)$ leads to (\ref{deriv_kf_neg}).} that
\begin{equation}
F_{-k}(y)=\frac{d^{k-1}}{dy^{k-1}}\left(  \frac{1}{1+e^{-y}}\right)
=\sum_{m=1}^{k}(m-1)!(-1)^{m-1}\mathcal{S}_{k}^{(m)}\,\left(  \frac
{1}{1+e^{-y}}\right)  ^{m} \label{deriv_kf_neg}%
\end{equation}
Thus, we find%
\begin{align*}
\frac{1}{2\pi i}\int_{c-i\infty}^{c+i\infty}\frac{\pi}{\sin\pi s}\frac{\left(
s-1\right)  ^{m-1}z^{s}}{s}ds  &  =\left(  -1\right)  ^{m-1}\left.
e^{-y}\frac{d^{m-1}}{dy^{m-1}}\left(  e^{y}\log\left(  1+e^{y}\right)
\right)  \right\vert _{z=e^{-y}}\\
&  =\left(  -1\right)  ^{m-1}\log\left(  1+\frac{1}{z}\right)  +\sum
_{l=1}^{m-1}\binom{m-1}{l}\sum_{v=1}^{l}(v-1)!(-1)^{m-v}\mathcal{S}_{l}%
^{(v)}\,\left(  \frac{1}{1+z}\right)  ^{v}%
\end{align*}
Reversal of the last double sum and with $\mathcal{S}_{m+1}^{(v+1)}=\sum
_{l=v}^{m}\binom{m}{l}\mathcal{S}_{l}^{(v)}$, we have%
\[
\frac{1}{2\pi i}\int_{c-i\infty}^{c+i\infty}\frac{\pi}{\sin\pi s}\frac{\left(
s-1\right)  ^{m-1}z^{s}}{s}ds=\left(  -1\right)  ^{m-1}\log\left(  1+\frac
{1}{z}\right)  +\sum_{v=1}^{m-1}(v-1)!(-1)^{m-v}\mathcal{S}_{m}^{(v+1)}%
\,\left(  \frac{1}{1+z}\right)  ^{v}%
\]
and%
\begin{align*}
\mu\left(  z\right)   &  =\left(  z+1\right)  \log\left(  1+\frac{1}%
{z}\right)  -1-\sum_{m=1}^{\infty}g_{m}\left(  1\right)  \frac{1}{2\pi i}%
\int_{c-i\infty}^{c+i\infty}\frac{\pi}{\sin\pi s}\frac{\left(  s-1\right)
^{m-1}z^{s}}{s}ds\\
&  =\left(  z+1\right)  \log\left(  1+\frac{1}{z}\right)  -1+\log\left(
1+\frac{1}{z}\right)  \sum_{m=1}^{\infty}g_{m}\left(  1\right)  \left(
-1\right)  ^{m}\\
&  \hspace{0.5cm}-\sum_{m=1}^{\infty}g_{m}\left(  1\right)  \sum_{v=1}%
^{m-1}(v-1)!(-1)^{m-v}\mathcal{S}_{m}^{(v+1)}\,\left(  \frac{1}{1+z}\right)
^{v}%
\end{align*}
Further, with $\left(  -1\right)  \zeta\left(  0\right)  =\sum_{m=0}^{\infty
}g_{m}\left(  1\right)  \left(  -1\right)  ^{m}=\frac{1}{2}$, we obtain%
\[
\mu\left(  z\right)  =\left(  z+\frac{1}{2}\right)  \log\left(  1+\frac{1}%
{z}\right)  -1-\sum_{m=1}^{\infty}g_{m}\left(  1\right)  \sum_{v=1}%
^{m-1}(v-1)!(-1)^{m-v}\mathcal{S}_{m}^{(v+1)}\,\left(  \frac{1}{1+z}\right)
^{v}%
\]
The forward difference formula (\ref{forwards_difference_mu}) shows that the
first terms are equal to $\mu\left(  z\right)  -\mu\left(  z+1\right)
=-\left(  z+\frac{1}{2}\right)  \log\frac{z}{z+1}-1$ and that%
\[
\mu\left(  z+1\right)  =\sum_{m=1}^{\infty}g_{m}\left(  1\right)  \sum
_{v=1}^{m-1}(v-1)!(-1)^{m-1-v}\mathcal{S}_{m}^{(v+1)}\,\left(  \frac{1}%
{1+z}\right)  ^{v}%
\]
which is (\ref{mu_series_TaylorZeta_1_op_z}), after replacing $z+1\rightarrow
z$.

\subsection{Taylor coefficients $g_{m}\left(  1\right)  $ of the Riemann Zeta
function}

\label{sec_Taylor_coefficient_RZeta_around1}The convergent Dirichlet series of
the Eta function $\eta(s)=\sum_{n=1}^{\infty}\frac{(-1)^{n+1}}{n^{s}}$ for
$\operatorname{Re}\left(  s\right)  \geq0$ immediately leads to the Taylor
expansion
\begin{equation}
\eta(s)=\sum_{k=0}^{\infty}\frac{\eta^{(k)}(s)}{k!}(s-x)^{k}
\label{def_Taylor_expansion_Etafunction_rond_x}%
\end{equation}
with
\begin{equation}
\eta^{(k)}(s)=(-1)^{k}\;\sum_{n=1}^{\infty}\frac{(-1)^{n+1}}{n^{s}}\;\log^{k}n
\label{afgeleiden_eta(s)_Dirichlet}%
\end{equation}
However, the Dirichlet series of $\eta^{(k)}(s)$ converges too slowly to be of
any practical use. Fortunately, fast converging series are obtained for real
$s\geq0$ by the Euler transform \cite{Hardy_div},
\begin{equation}
\eta^{(k)}\left(  s\right)  =\left(  -1\right)  ^{k+1}\sum_{m=1}^{\infty
}\left[  \sum_{j=1}^{m}{\binom{m-1}{j-1}}\frac{\left(  -1\right)  ^{j}\ln
^{k}j}{j^{s}}\,\right]  \left(  \frac{1}{2}\right)  ^{m}
\label{afgeleiden_eta(s)_Euler}%
\end{equation}

Invoking the relation $\zeta(s)=$ $\frac{\eta(s)}{1-2^{1-s}}$ and using the
generating function (\ref{gf_Bernoullinumbers}) of the Bernoulli numbers, we
have
\[
\frac{1}{1-2^{1-s}}=-\frac{1}{e^{-(s-1)\log2}-1}=\frac{1}{(s-1)\log2}%
\sum_{n=0}^{\infty}B_{n}\,\frac{\left(  -\log2\right)  ^{n}}{n!}(s-1)^{n}%
\]
After executing the Cauchy product of the Taylor series for $\frac
{1}{1-2^{1-s}}$ and that of the Eta function in
(\ref{def_Taylor_expansion_Etafunction_rond_x}), we obtain
\begin{align*}
\zeta(s)  &  =\frac{1}{(s-1)\log2}\sum_{k=0}^{\infty}\left[  \sum_{j=0}%
^{k}B_{j}\,\frac{\left(  -\log2\right)  ^{j}}{j!}\frac{\eta^{(k-j)}%
(1)}{(k-j)!}\right]  (s-1)^{k}\\
&  =\frac{1}{s-1}+\sum_{k=1}^{\infty}\left[  \sum_{j=0}^{k}B_{j}%
\,\frac{(-1)^{j}\log^{j-1}2}{j!}\frac{\eta^{(k-j)}(1)}{(k-j)!}\right]
(s-1)^{k-1}%
\end{align*}
where we have used that $\eta(1)=\log2$. Equating corresponding powers in
$\left(  s-1\right)  $ in both Taylor series of $\left(  s-1\right)
\zeta\left(  s\right)  $ yields, with $g_{0}(1)=1$ and for $k>0$,
\begin{equation}
g_{k}(1)=\frac{1}{k!}\sum_{j=0}^{k}\binom{k}{j}B_{j}\,(-1)^{j}\eta
^{(k-j)}(1)\log^{j-1}2 \label{Taylor_zeta_rond1}%
\end{equation}

The Taylor coefficient of $(s-1)\zeta(s)=\sum_{m=0}^{\infty}$ $g_{m}(1)\left(
s-1\right)  ^{m}$ around $s_{0}=1$ follow from (\ref{Taylor_zeta_rond1})
as{\small
\[%
\begin{array}
[c]{cclcccl}%
g_{0}(1) & = & 1.0 &  & g_{1}(1) & = & \gamma=0.5772156649015328606\\
g_{2}(1) & = & 0.07281584548367672486 &  & g_{3}(1) & = &
-0.004845181596436159243\\
g_{4}(1) & = & -0.000342305736717224311 &  & g_{5}(1) & = &
0.00009689041939447083573\\
g_{6}(1) & = & -6.611031810842189181\;10^{-6} &  & g_{7}(1) & = &
-3.31624090875277236\;10^{-7}\\
g_{8}(1) & = & 1.0462094584479187422\;10^{-7} &  & g_{9}(1) & = &
-8.733218100273797361\;10^{-9}\\
g_{10}(1) & = & 9.478277782762358956\;10^{-11} &  & g_{11}(1) & = &
5.658421927608707966\;10^{-11}\\
g_{12}(1) & = & -6.768689863513696656\;10^{-12} &  & g_{13}(1) & = &
3.492115936672031855\;10^{-13}\\
g_{14}(1) & = & 4.41042474175775338\;10^{-15} &  & g_{15}(1) & = &
-2.3997862217709991766\;10^{-15}\\
g_{16}(1) & = & 2.167731220072682855\;10^{-16} &  & g_{17}(1) & = &
-9.54446607636696516\;10^{-18}\\
g_{18}(1) & = & -7.387676660538636498\;10^{-20} &  & g_{19}(1) & = &
4.800850782488065211\;10^{-20}\\
g_{20}(1) & = & -4.139956737713305639\;10^{-21} &  & g_{21}(1) & = &
1.19168201593979951\;10^{-22}%
\end{array}
\]
}

\section{Taylor series of $\frac{x}{\log^{n}\left(  1+x\right)  }$ for integer
$n$}

\label{sec_Taylor_series_x/log(1+x)}Integrating the double generating
function
\[
\left(  1+x\right)  ^{u}=e^{u\log\left(  1+x\right)  }=\sum_{m=0}^{\infty}%
\sum_{k=0}^{m}\frac{S_{m}^{(k)}}{m!}u^{k}x^{m}%
\]
of the Stirling numbers $S_{m}^{(k)}$ of the First Kind \cite[Sec. 24.1.3 and
24.1.4]{Abramowitz} with respect to $u$ results, for $\left\vert x\right\vert
<1$, in%
\[
\frac{e^{b\log\left(  1+x\right)  }-e^{a\log\left(  1+x\right)  }}{\log\left(
1+x\right)  }=\sum_{m=0}^{\infty}\sum_{k=0}^{m}S_{m}^{(k)}\frac{b^{k+1}%
-a^{k+1}}{k+1}\frac{x^{m}}{m!}%
\]
In particular, for $b=1$ and $a=0$, we obtain the Taylor series, valid for
$\left\vert x\right\vert <1$,%
\begin{equation}
\frac{x}{\log\left(  1+x\right)  }=\sum_{m=0}^{\infty}\left(  \sum_{k=0}%
^{m}\frac{S_{m}^{(k)}}{k+1}\right)  \frac{x^{m}}{m!}=1+\sum_{m=1}^{\infty
}\left(  \sum_{k=1}^{m}\frac{S_{m}^{(k)}}{k+1}\right)  \frac{x^{m}}{m!}
\label{Taylor_series_x/Log(1+x)}%
\end{equation}

We generalize the above. The $n$-fold integral of $e^{u\lambda}$ equals%
\[
\int_{a}^{b}du_{1}\int_{a}^{u_{1}}du_{2}\ldots\int_{a}^{u_{n-1}}du_{n}%
e^{u_{n}\lambda}=\frac{1}{\left(  n-1\right)  !}\int_{a}^{b}\left(
b-u\right)  ^{n-1}e^{u\lambda}du
\]
Let $t=b-u$, followed by $y=\lambda t$, then%
\[
\int_{a}^{b}\left(  b-u\right)  ^{n-1}e^{u\lambda}du=\frac{e^{\lambda b}%
}{\lambda^{n}}\int_{0}^{\lambda\left(  b-a\right)  }y^{n-1}e^{-y}dy
\]
and the integral can be executed leading to%
\[
\frac{1}{\left(  n-1\right)  !}\int_{a}^{b}\left(  b-u\right)  ^{n-1}%
e^{u\lambda}du=\frac{e^{\lambda b}}{\lambda^{n}}\left(  1-e^{-\lambda\left(
b-a\right)  }\sum_{k=0}^{m}\frac{\left(  \lambda\left(  b-a\right)  \right)
^{k}}{k!}\right)
\]
On the other hand, $e^{u\lambda}=\sum_{k=0}^{\infty}\frac{\lambda^{k}u^{k}%
}{k!}$ and the $n$-fold integration of $u^{k}$ is
\begin{equation}
\int_{a}^{b}du_{1}\int_{a}^{u_{1}}du_{2}\ldots\int_{a}^{u_{n-1}}du_{n}%
u_{n}^{k}=\frac{1}{\left(  n-1\right)  !}\int_{a}^{b}\left(  b-u\right)
^{n-1}u^{k}du \label{n-fold_integration_power_u}%
\end{equation}
and%
\[
\int_{a}^{b}\left(  b-u\right)  ^{n-1}u^{k}du=b^{n}\int_{a}^{b}\left(
1-\frac{u}{b}\right)  ^{n-1}u^{k}du=b^{n+k}\int_{\frac{a}{b}}^{1}\left(
1-w\right)  ^{n-1}w^{k}dw
\]
which simplifies considerably if $a=0$, due to the Beta integral $\int_{0}%
^{1}\left(  1-w\right)  ^{n-1}w^{k}dw=\frac{\Gamma\left(  n\right)
\Gamma\left(  k+1\right)  }{\Gamma\left(  n+k+1\right)  }$. Thus, choosing
$a=0$ leads to%
\[
\frac{1}{\lambda^{n}}\left(  e^{\lambda b}-\sum_{k=0}^{m}\frac{\left(  \lambda
b\right)  ^{k}}{k!}\right)  =b^{n}\sum_{k=0}^{\infty}\frac{\lambda^{k}b^{k}%
}{\left(  n+k\right)  !}%
\]
Further, with $\lambda=\log\left(  1+x\right)  $ and introducing the Taylor
series $\lambda^{k}=\log^{k}(1+x)=k!\sum_{m=k}^{\infty}S_{m}^{(k)}%
\;\frac{x^{m}}{m!}$, valid for $x|<1$, yields, after reversal of the $k$- and
$m$-sum,%
\[
\frac{1}{\log^{n}\left(  1+x\right)  }\left(  e^{\log\left(  1+x\right)
b}-\sum_{k=0}^{n-1}\frac{\left(  \log\left(  1+x\right)  b\right)  ^{k}}%
{k!}\right)  =b^{n}\sum_{m=0}^{\infty}\left(  \sum_{k=0}^{m}\frac
{k!S_{m}^{(k)}b^{k}}{\left(  k+n\right)  !}\right)  \frac{x^{m}}{m!}%
\]
valid for $\left\vert x\right\vert <1$ and which simplifies for $b=1$ to,
\begin{equation}
\frac{x}{\log^{n}\left(  1+x\right)  }=\frac{1}{n!}+\sum_{k=1}^{n-1}\frac
{1}{k!\log^{n-k}\left(  1+x\right)  }+\sum_{m=1}^{\infty}\left(  \sum
_{k=1}^{m}\frac{k!S_{m}^{(k)}}{\left(  k+n\right)  !}\right)  \frac{x^{m}}{m!}
\label{Taylor_series_x/Log^n(1+x)}%
\end{equation}
For $n=1$ in (\ref{Taylor_series_x/Log^n(1+x)}), we find again
(\ref{Taylor_series_x/Log(1+x)}).

Applying $n$-fold integration to the generating function
(\ref{genfunc_stirling})%
\[
\int_{0}^{b}du_{1}\int_{0}^{u_{1}}du_{2}\ldots\int_{0}^{u_{n-1}}du_{n}%
\prod_{k=0}^{m-1}(u_{n}-k)=\sum_{k=0}^{m}S_{m}^{(k)}\frac{k!}{\left(
n+k\right)  !}b^{k+n}%
\]
and executing the left-hand side (via partial integration) yields
\[
\frac{1}{\left(  n-1\right)  !}\int_{0}^{b}\left(  b-u\right)  ^{n-1}%
\prod_{k=0}^{m-1}(u-k)du=\sum_{k=0}^{m}S_{m}^{(k)}\frac{k!}{\left(
n+k\right)  !}b^{k+n}%
\]
which links Stirling numbers to the general integral form used by Binet
\cite[p. 339]{Binet1839} in the series expansion of his Binet function
$\mu\left(  z\right)  $.

\section{The Taylor series (\ref{Taylor_series_g_alpha(u)})}

\label{sec_Taylor_series_c[m]}

Inspired by Nemes \cite[Section 3]{Nemes_2013} and using the integral
(\ref{generalized_Binet_numbers_alpha_integral}), we compute the
\textquotedblleft exponential\textquotedblright\ generating function of the
Binet polynomials $b_{m}\left(  \alpha\right)  $,%
\[
\sum_{m=1}^{\infty}\frac{b_{m}\left(  \alpha\right)  }{\left(  m-1\right)
!}u^{m}=\int_{\alpha-1}^{\alpha}\left(  x+\left(  \frac{1}{2}-\alpha\right)
\right)  \sum_{m=1}^{\infty}\frac{u^{m}}{m!}\prod_{k=0}^{m-1}(k+x)dx
\]
From (\ref{genfunc_stirling}), it follows that $\prod_{k=0}^{m-1}%
(k+x)=\frac{\left(  -1\right)  ^{m}\Gamma(1-x)}{\Gamma(1-x-m)}$ and
$\frac{\Gamma(1-x)}{\Gamma(1-x-m)m!}=\binom{-x}{m}$. Provided that $\left\vert
u\right\vert <1$, the binomial sum equals%
\[
\sum_{m=1}^{\infty}\frac{u^{m}}{m!}\prod_{k=0}^{m-1}(k+x)=\sum_{m=1}^{\infty
}\binom{-x}{m}\left(  -1\right)  ^{m}u^{m}=\left(  1-u\right)  ^{-x}-1
\]
Hence, we obtain, for $\left\vert u\right\vert <1$,%
\[
\sum_{m=1}^{\infty}\frac{b_{m}\left(  \alpha\right)  }{\left(  m-1\right)
!}u^{m}=\int_{\alpha-1}^{\alpha}\left(  x+\left(  \frac{1}{2}-\alpha\right)
\right)  e^{-x\log\left(  1-u\right)  }dx-\int_{\alpha-1}^{\alpha}\left(
x+\left(  \frac{1}{2}-\alpha\right)  \right)  dx
\]
and%
\begin{equation}
g_{\alpha}\left(  u\right)  =-\frac{1}{2}\frac{\left(  1-u\right)  ^{-\alpha
}+\left(  1-u\right)  ^{1-\alpha}}{\log\left(  1-u\right)  }+\frac{\left(
1-u\right)  ^{1-\alpha}-\left(  1-u\right)  ^{-\alpha}}{\log^{2}\left(
1-u\right)  }=\sum_{m=1}^{\infty}\frac{b_{m}\left(  \alpha\right)  }{\left(
m-1\right)  !}u^{m} \label{exponential_generating_function_b_m(alpha)}%
\end{equation}
which reduces, for $\alpha=0$, to the Taylor series
(\ref{Taylor_series_bijna_cm_coefficienten}) in Binet's second derivation
(\ref{Binet_function_integral}) in Theorem
\ref{theorem_Binet_convergent_expansion}. Moreover,%
\[
g_{\alpha}\left(  u\right)  =\left(  1-u\right)  ^{-\alpha}\left(
\frac{\left(  \frac{u}{2}-1\right)  }{\log\left(  1-u\right)  }-\frac{u}%
{\log^{2}\left(  1-u\right)  }\right)  =\left(  1-u\right)  ^{-\alpha}%
g_{0}\left(  u\right)
\]
where $g_{0}\left(  u\right)  =\sum_{m=1}^{\infty}\frac{c_{m}}{\left(
m-1\right)  !}u^{m}$ was prominent in Binet's proof of Theorem
\ref{theorem_Binet_convergent_expansion}.

We make Binet's substitution $u=1-e^{-t}$ in
(\ref{exponential_generating_function_b_m(alpha)}) and obtain%
\[
\frac{e^{\alpha t}+e^{\left(  \alpha-1\right)  t}}{2t}+\frac{e^{\left(
\alpha-1\right)  t}-e^{\alpha t}}{t^{2}}=\sum_{m=1}^{\infty}\frac{b_{m}\left(
\alpha\right)  }{\left(  m-1\right)  !}\left(  1-e^{-t}\right)  ^{m}%
\]
The Taylor series around $t_{0}=0$ of the left-hand side is%
\[
\frac{1}{t}\left(  \frac{e^{\alpha t}+e^{\left(  \alpha-1\right)  t}}{2}%
+\frac{e^{\left(  \alpha-1\right)  t}-e^{\alpha t}}{t}\right)  =\sum
_{k=0}^{\infty}\left\{  \frac{\frac{k+2}{2}\left(  \alpha^{k+1}+\left(
\alpha-1\right)  ^{k+1}\right)  +\left(  \left(  \alpha-1\right)
^{k+2}-\alpha^{k+2}\right)  }{\left(  k+2\right)  !}\right\}  t^{k}%
\]
while the right-hand side is%
\[
\sum_{m=1}^{\infty}\frac{b_{m}\left(  \alpha\right)  }{\left(  m-1\right)
!}\left(  1-e^{-t}\right)  ^{m}=\sum_{m=1}^{\infty}\frac{b_{m}\left(
\alpha\right)  }{\left(  m-1\right)  !}\sum_{j=0}^{m}\binom{m}{j}\left(
-1\right)  ^{j}e^{-jt}%
\]
but the reversal of the $m$- and $k$-sum is not allowed\footnote{Indeed,
$\sum_{j=0}^{\infty}\left(  \sum_{m=j}^{\infty}\frac{mb_{m}\left(
\alpha\right)  }{\left(  m-j\right)  !}\right)  \frac{\left(  -1\right)
^{j}e^{-jt}}{j!}$ diverges any $j>0$, because
\[
\sum_{m=j}^{\infty}\frac{mb_{m}\left(  \alpha\right)  }{\left(  m-j\right)
!}=\lim_{u\rightarrow1}\frac{d^{j}}{du^{j}}g_{\alpha}\left(  u\right)
=\lim_{u\rightarrow1}\frac{d^{j}}{du^{j}}\left(  1-u\right)  ^{-\alpha}%
g_{0}\left(  u\right)
\]
and $\sum_{m=j}^{\infty}\frac{mc_{m}}{\left(  m-j\right)  !}=\lim
_{u\rightarrow1}\frac{d^{j}}{du^{j}}g_{0}\left(  u\right)  =-\infty$ for
$j>0$, but $g_{0}\left(  1\right)  =\sum_{m=1}^{\infty}\frac{c_{m}}{\left(
m-1\right)  !}=0$, that already appeared in Section \ref{sec_growth_cm_with_m}
in the determination of $\lim_{z\rightarrow0}z\mu\left(  z\right)  =0$.}.
After Taylor expansion of $e^{-jt}=\sum_{k=0}^{\infty}j^{k}\frac{\left(
-1\right)  ^{k}t^{k}}{k!}$ around $t_{0}=0$,%
\[
\sum_{m=1}^{\infty}\frac{b_{m}\left(  \alpha\right)  }{\left(  m-1\right)
!}\left(  1-e^{-t}\right)  ^{m}=\sum_{m=1}^{\infty}\frac{b_{m}\left(
\alpha\right)  }{\left(  m-1\right)  !}\sum_{k=0}^{\infty}\left(  \sum
_{j=0}^{m}\binom{m}{j}\left(  -1\right)  ^{j}j^{k}\right)  \frac{\left(
-1\right)  ^{k}t^{k}}{k!}%
\]
and recognizing the closed form expression \cite[sec. 24.1.4.C]{Abramowitz} of
the Stirling number of the Second Kind
\begin{equation}
\mathcal{S}_{k}^{(m)}=\frac{1}{m!}\sum_{j=0}^{m}(-1)^{m-j}{\binom{m}{j}}j^{k}
\label{stirling2closed}%
\end{equation}
we have%
\[
\sum_{m=1}^{\infty}\frac{b_{m}\left(  \alpha\right)  }{\left(  m-1\right)
!}\left(  1-e^{-t}\right)  ^{m}=\sum_{k=0}^{\infty}\left(  \sum_{m=1}^{\infty
}mb_{m}\left(  \alpha\right)  \mathcal{S}_{k}^{(m)}\left(  -1\right)
^{m-k}\right)  \frac{t^{k}}{k!}%
\]
Since the positive integer $\mathcal{S}_{k}^{(m)}=0$ for $m>k$, we finally
arrive at%
\[
\sum_{k=0}^{\infty}\left\{  \frac{\frac{k+2}{2}\left(  \alpha^{k+1}+\left(
\alpha-1\right)  ^{k+1}\right)  +\left(  \left(  \alpha-1\right)
^{k+2}-\alpha^{k+2}\right)  }{\left(  k+2\right)  !}\right\}  t^{k}=\sum
_{k=0}^{\infty}\left(  \sum_{m=1}^{k}mb_{m}\left(  \alpha\right)
\mathcal{S}_{k}^{(m)}\left(  -1\right)  ^{m-k}\right)  \frac{t^{k}}{k!}%
\]
which is, after equating corresponding powers in $t$, again
(\ref{recursion_general_coeff_b_m(alpha)_StirlingS2}).

\section{Integral for the Binet coefficient $c_{m}$}

\label{sec_integral_Binet_coefficient}The Taylor series
(\ref{Taylor_series_bijna_cm_coefficienten}), which is a special case of
(\ref{Taylor_series_g_alpha(u)}) for $\alpha=0$, is written in terms of the
Binet coefficients $c_{m}$ in (\ref{def_Binet_coefficients}) as%
\[
g_{0}\left(  -u\right)  =\frac{u}{\log^{2}\left(  1+u\right)  }-\frac
{1+\frac{1}{2}u}{\log\left(  1+u\right)  }=\sum_{m=1}^{\infty}\frac{\left(
-1\right)  ^{m}c_{m}}{\left(  m-1\right)  !}u^{m}%
\]
The corresponding integral form for the Taylor coefficient
\cite{Titchmarshfunctions} is%
\[
\frac{\left(  -1\right)  ^{m}c_{m}}{\left(  m-1\right)  !}=\frac{1}{2\pi
i}\int_{C(0)}\frac{d\omega\;}{\omega^{m+1}}\left(  \frac{\omega}{\log
^{2}\left(  1+\omega\right)  }-\frac{1+\frac{1}{2}\omega}{\log\left(
1+\omega\right)  }\right)
\]
where $C\left(  0\right)  $ is a contour that encloses the point $u_{0}=0$ in
counter-clockwise sense. A straightforward execution of the contour $C\left(
0\right)  $ is to choose a circle around the origin, with radius $0<r\leq1$,
due to the branch cut at the negative real axis for $\operatorname{Re}\left(
\omega\right)  <-1$. The resulting integral is numerically not stable. An
alternative way is to deform the contour to enclose the entire complex plane
(except for the point $\omega=0$ and avoiding the branch cut) in clockwise
sense. The integrand vanishes at $\omega=-1$. For $\left\vert \omega
\right\vert \rightarrow\infty$, the integrand vanishes for $m>1$ and we only
maintain the path around the branch cut. In particular, we construct a path
that travels from infinity to the point $-1\leq q<0$ under an angle $-\theta$,
where $\theta\in\left(  0,\pi\right)  $ and returns from the point $q$ along a
straight line under angle $\theta$ to infinity. We thus obtain%
\begin{align*}
\frac{\left(  -1\right)  ^{m}c_{m}}{\left(  m-1\right)  !}  &  =-\frac{1}{2\pi
i}\int_{\infty}^{0}\frac{d\left(  q+xe^{-i\theta}\right)  \;}{\left(
q+xe^{-i\theta}\right)  ^{m+1}}\left(  \frac{q+xe^{-i\theta}}{\log^{2}\left(
1+q+xe^{-i\theta}\right)  }-\frac{1+\frac{1}{2}q+\frac{1}{2}xe^{-i\theta}%
}{\log\left(  1+q+xe^{-i\theta}\right)  }\right) \\
&  -\frac{1}{2\pi i}\int_{0}^{\infty}\frac{d\left(  q+xe^{i\theta}\right)
\;}{\left(  q+xe^{i\theta}\right)  ^{m+1}}\left(  \frac{q+xe^{i\theta}}%
{\log^{2}\left(  1+q+xe^{i\theta}\right)  }-\frac{1+\frac{1}{2}q+\frac{1}%
{2}xe^{i\theta}}{\log\left(  1+q+xe^{i\theta}\right)  }\right)
\end{align*}
The computation simplifies if we choose $q=-1$,%
\[
\frac{\left(  -1\right)  ^{m-1}c_{m}}{\left(  m-1\right)  !}=\frac{1}{\pi}%
\int_{0}^{\infty}\operatorname{Im}\left(  \frac{e^{i\theta}}{\left(
-1+xe^{i\theta}\right)  ^{m+1}}\left(  \frac{-1+xe^{i\theta}}{\left(  \log
x+i\theta\right)  ^{2}}-\frac{1}{2}\frac{1+xe^{i\theta}}{\log x+i\theta
}\right)  \right)  \;dx
\]
We evaluate the integrand. Denoting%
\[
\Psi\left(  x\right)  =\theta-\left(  m+1\right)  \arccos\frac{x\cos\theta
-1}{\sqrt{x^{2}-2x\cos\theta+1}}+\arccos\frac{x\cos\theta+1}{\sqrt
{x^{2}+2x\cos\theta+1}}-\arctan\frac{\theta}{\log x}%
\]
we obtain%
\begin{align*}
\frac{\left(  -1\right)  ^{m-1}c_{m}}{\left(  m-1\right)  !}  &  =\frac{1}%
{\pi}\int_{0}^{\infty}\frac{dx}{\log^{2}x+\theta^{2}}\frac{\sin\left(
\theta-m\arccos\frac{x\cos\theta-1}{\sqrt{x^{2}-2x\cos\theta+1}}-2\arctan
\frac{\theta}{\log x}\right)  }{\left(  x^{2}-2x\cos\theta+1\right)
^{\frac{m}{2}}}\\
&  -\frac{1}{2\pi}\int_{0}^{\infty}\frac{\sqrt{\left(  1+2x\cos\theta
+x^{2}\right)  }}{\sqrt{\log^{2}x+\theta^{2}}}\frac{\sin\Psi\left(  x\right)
}{\left(  x^{2}-2x\cos\theta+1\right)  ^{\frac{m+1}{2}}}dx
\end{align*}
This form simplifies substantially if we choose $\theta=\frac{\pi}{2}$. After
simplifying the sines, we arrive at%
\begin{equation}
\frac{\left(  -1\right)  ^{m-1}c_{m}}{\left(  m-1\right)  !}=\frac{1}{\pi}%
\int_{0}^{\infty}\frac{\frac{\cos\left(  m\arccos\frac{-1}{\sqrt{x^{2}+1}%
}+2\arctan\frac{\theta}{\log x}\right)  }{\log^{2}x+\left(  \frac{\pi}%
{2}\right)  ^{2}}+\frac{\cos\left(  \left(  m+2\right)  \arccos\frac{-1}%
{\sqrt{x^{2}+1}}+\arctan\frac{\theta}{\log x}\right)  }{2\sqrt{\log
^{2}x+\left(  \frac{\pi}{2}\right)  ^{2}}}}{\left(  x^{2}+1\right)  ^{\frac
{m}{2}}}dx \label{Integral_Binet_coefficient_c[m]_cosines}%
\end{equation}
However, the numerical evaluation of the integral
(\ref{Integral_Binet_coefficient_c[m]_cosines}) is remarkably inaccurate.
Therefore, we simplify the cosines. After some manipulations, we arrive at an
integral for the Binet coefficient $c_{m}$ for $m>1$,%
\begin{align}
\frac{\left(  -1\right)  ^{m-1}c_{m}}{\left(  m-1\right)  !}  &  =\frac{1}%
{\pi}\int_{0}^{\infty}\frac{\cos\left(  m\arccos\frac{-1}{\sqrt{x^{2}+1}%
}\right)  \left\{  \frac{\log^{2}x-\left(  \frac{\pi}{2}\right)  ^{2}}%
{\log^{2}x+\left(  \frac{\pi}{2}\right)  ^{2}}+\frac{\left(  1-x^{2}\right)
\frac{\log x}{2}+\frac{\pi}{2}x}{\left(  1+x^{2}\right)  }\right\}  }{\left(
x^{2}+1\right)  ^{\frac{m}{2}}\left(  \log^{2}x+\left(  \frac{\pi}{2}\right)
^{2}\right)  }dx\nonumber\\
&  +\frac{1}{\pi}\int_{0}^{\infty}\frac{\sin\left(  m\arccos\frac{-1}%
{\sqrt{x^{2}+1}}\right)  \left\{  \frac{-\pi\log x}{\log^{2}x+\left(
\frac{\pi}{2}\right)  ^{2}}+\frac{x\log x-\frac{\pi}{4}\left(  1-x^{2}\right)
}{\left(  1+x^{2}\right)  }\right\}  }{\left(  x^{2}+1\right)  ^{\frac{m}{2}%
}\left(  \log^{2}x+\left(  \frac{\pi}{2}\right)  ^{2}\right)  }dx
\label{Integral_Binet_coefficient_c[m]}%
\end{align}
where $\cos\left(  m\arccos y\right)  $ is the Chebyshev orthogonal polynomial
in $y$. The integral (\ref{Integral_Binet_coefficient_c[m]}) can be evaluated accurately.

Upper bounding the cosines in (\ref{Integral_Binet_coefficient_c[m]_cosines})
leads to%
\[
\left\vert \frac{\left(  -1\right)  ^{m-1}c_{m}}{\left(  m-1\right)
!}\right\vert \leq\frac{1}{\pi}\int_{0}^{\infty}\frac{dx}{\left(
x^{2}+1\right)  ^{\frac{m}{2}}}\left\{  \frac{1}{\log^{2}x+\left(  \frac{\pi
}{2}\right)  ^{2}}+\frac{1}{2\sqrt{\log^{2}x+\left(  \frac{\pi}{2}\right)
^{2}}}\right\}
\]
The function in between brackets $\left\{  .\right\}  $ is maximal at $x=1$,
where it equals $\frac{4}{\pi^{2}}+\frac{1}{\pi}$. Thus,%
\[
\left\vert \frac{\left(  -1\right)  ^{m-1}c_{m}}{\left(  m-1\right)
!}\right\vert \leq\frac{1}{\pi}\left(  \frac{4}{\pi^{2}}+\frac{1}{\pi}\right)
\int_{0}^{\infty}\frac{dx}{\left(  x^{2}+1\right)  ^{\frac{m}{2}}}=\frac
{1}{2\sqrt{\pi}}\left(  \frac{4}{\pi^{2}}+\frac{1}{\pi}\right)  \frac
{\Gamma\left(  \frac{m-1}{2}\right)  }{\Gamma\left(  \frac{m}{2}\right)
}=O\left(  \frac{1}{\sqrt{m}}\right)
\]
but this upper bound is rather weak. In particular, since all coefficients
$c_{m}$ for $m>2$ have the same sign by Theorem
\ref{theorem_Binet_coefficient_integral}, the convergence of $\sum
_{m=1}^{\infty}\frac{c_{m}}{\left(  m-1\right)  !}$ implies that
$\frac{\left\vert c_{m}\right\vert }{\left(  m-1\right)  !}=O\left(  \frac
{1}{m^{1+\varepsilon}}\right)  $ for $\varepsilon>0$.

\section{Asymptotic expansion for $b_{m}(\frac{1}{2}+\alpha)$}

\label{sec_asymptotic_bm(alpha)_large_m}We start from the integral in
(\ref{generalized_Binet_numbers_alpha_integral_alpha+1/2}),
\begin{align*}
\left(  m+1\right)  b_{m+1}\left(  \frac{1}{2}+\alpha\right)   &
=\int_{-\frac{1}{2}}^{\frac{1}{2}}u\left(  u+\alpha\right)  \prod_{k=1}%
^{m}(k+u+\alpha)du\\
&  =\prod_{k=1}^{m}\left(  k+\alpha\right)  \int_{-\frac{1}{2}}^{\frac{1}{2}%
}u\left(  u+\alpha\right)  \prod_{k=1}^{m}\left(  1+\frac{u}{k+\alpha}\right)
du
\end{align*}
Provided $\frac{\left\vert u\right\vert }{\left\vert k+\alpha\right\vert }<1$,
the product can be expanded around $u_{0}=0$ as%
\begin{align*}
\prod_{k=n}^{m}\left(  1+\frac{u}{k+\alpha}\right)   &  =\exp\left(
\sum_{k=n}^{m}\log\left(  1+\frac{u}{k+\alpha}\right)  \right)  =\exp\left(
\sum_{j=1}^{\infty}\frac{\left(  -1\right)  ^{j-1}}{j}\sum_{k=n}^{m}\frac
{1}{\left(  k+\alpha\right)  ^{j}}u^{j}\right) \\
&  =\exp\left(  u\sum_{k=n}^{m}\frac{1}{\left(  k+\alpha\right)  }\right)
\exp\left(  \sum_{j=2}^{\infty}\frac{\left(  -1\right)  ^{j-1}}{j}\sum
_{k=n}^{m}\frac{1}{\left(  k+\alpha\right)  ^{j}}u^{j}\right)
\end{align*}
The convergence requirement indicates that $\left\vert k+\alpha\right\vert
>\frac{1}{2}$ for $k\geq n\geq1$, which means that $\frac{1}{2}-n<\alpha$. We
limit ourselves here to $n=1$, implying that the analysis is valid for real
$\alpha>-\frac{1}{2}$ and, thus, after translating $b_{m}(\frac{1}{2}+\alpha)$
to $b_{m}\left(  \alpha^{\prime}\right)  $ for $\alpha^{\prime}>0$. If that
range must be larger, then we can increase $n\geq2$, so that $\alpha>\frac
{3}{2}$ and so on; the only effect is that the integral $I_{j}$ below is a
little more involved, but still analytically computable. For $j\geq2$, the sum
$\sum_{k=n}^{m}\frac{1}{\left(  k+\alpha\right)  ^{j}}$ converges for all $m$,
whereas $\sum_{k=n}^{m}\frac{1}{\left(  k+\alpha\right)  }$ diverges when
$m\rightarrow\infty$, which justifies the split-off of the $j=1$ term. The
limit $m\rightarrow\infty$ case can be expressed in terms of the Hurwitz
Zeta-function $\zeta\left(  s,\alpha\right)  =\sum_{k=1}^{\infty}\frac
{1}{\left(  k+\alpha\right)  ^{s}}$ (see Appendix
\ref{sec_complex_integral_mu}). The remaining $j$-series is alternating with
decreasing coefficients and can thus be bounded as%
\[
-\frac{u^{2}}{2}\sum_{k=n}^{m}\frac{1}{\left(  k+\alpha\right)  ^{2}}%
<\sum_{j=2}^{\infty}\frac{\left(  -1\right)  ^{j-1}}{j}\sum_{k=n}^{m}\frac
{1}{\left(  k+\alpha\right)  ^{j}}u^{j}<-\frac{u^{2}}{2}\sum_{k=n}^{m}\frac
{1}{\left(  k+\alpha\right)  ^{2}}+\frac{u^{3}}{3}\sum_{k=n}^{m}\frac
{1}{\left(  k+\alpha\right)  ^{3}}%
\]
Rather than continuing with these bounds, we proceed with an exact computation
using our characteristic coefficients \cite[Appendix]{PVM_ASYM}, that enables
us to expand $\exp\left(  \sum_{j=2}^{\infty}\frac{\left(  -1\right)  ^{j-1}%
}{j}\sum_{k=n}^{m}\frac{1}{\left(  k+\alpha\right)  ^{j}}u^{j}\right)
=\sum_{l=0}^{\infty}\phi_{l}u^{l}$ in a Taylor series around $u_{0}=0$. The
Taylor series of a function $G\left(  z\right)  $ of a function $f\left(
z\right)  $ is
\begin{equation}
G(f(z))=\sum_{m=0}^{\infty}\left(  \sum_{k=0}^{m}\frac{1}{k!}\left.
\frac{d^{k}G(f)}{df^{k}}\right\vert _{f=f\left(  z_{0}\right)  }%
s[k,m]_{f\left(  z\right)  }(z_{0})\,\right)  (z-z_{0})^{m} \label{Gf_rond_z0}%
\end{equation}
The characteristic coefficients of a complex function $f\left(  z\right)  $
with Taylor series $f\left(  z\right)  =\sum_{k=0}^{\infty}f_{k}\left(
z_{0}\right)  \left(  z-z_{0}\right)  ^{k}$, defined by $\left.
s[k,m]\right\vert _{f}\left(  z_{0}\right)  =\frac{1}{m!}\left.  \frac{d^{m}%
}{dz^{m}}\left(  f\left(  z\right)  -f\left(  z_{0}\right)  ^{k}\right)
\right\vert _{z=z_{0}}$, possesses a general form
\begin{equation}
\left.  s[k,m]\right\vert _{f}\left(  z_{0}\right)  =\sum_{\sum_{i=1}^{k}%
j_{i}=m;j_{i}>0}\prod_{i=1}^{k}f_{j_{i}}\left(  z_{0}\right)  \label{def_s}%
\end{equation}
and obeys $\left.  s[k,m]\right\vert _{f}\left(  z_{0}\right)  =0$ if $k<0$
and $k>m$. Moreover, $\left.  s[k,m]\right\vert _{f}\left(  z_{0}\right)  $
possesses a recursion and the coefficients $\phi_{l}$ can be computed up to
any desired order. The function $f\left(  u\right)  =\sum_{j=2}^{\infty
}\left(  \frac{\left(  -1\right)  ^{j-1}}{j}\sum_{k=n}^{m}\frac{1}{\left(
k+\alpha\right)  ^{j}}\right)  u^{j}$ has clearly two vanishing Taylor
coefficients, $f_{0}=f_{1}=0$, while $f_{j}=\frac{\left(  -1\right)  ^{j-1}%
}{j}\sum_{k=n}^{m}\frac{1}{\left(  k+\alpha\right)  ^{j}}$. Invoking
(\ref{Gf_rond_z0})%
\[
e^{f(u)}=1+\sum_{m=1}^{\infty}\left[  \sum_{k=1}^{m}\frac{1}{k!}%
\,s[k,m]\right]  \,u^{m}%
\]
indicates that $\phi_{0}=1$ and $\phi_{l}=\sum_{k=1}^{l}\frac{1}{k!}\,s[k,l]$.
Because $f_{1}=0$, it holds that $\phi_{1}=s[1,1]=0$. We list the first Taylor
coefficients $\phi_{l}$,%
\begin{align*}
\phi_{2}  &  =-\frac{1}{2}\sum_{k=n}^{m}\frac{1}{\left(  k+\alpha\right)
^{2}}\text{ and }\phi_{3}=\frac{1}{3}\sum_{k=n}^{m}\frac{1}{\left(
k+\alpha\right)  ^{3}}\\
\phi_{4}  &  =\frac{1}{8}\left(  \sum_{k=n}^{m}\frac{1}{\left(  k+\alpha
\right)  ^{2}}\right)  ^{2}-\frac{1}{4}\sum_{k=n}^{m}\frac{1}{\left(
k+\alpha\right)  ^{4}}\\
\phi_{5}  &  =-\frac{1}{6}\sum_{k=n}^{m}\frac{1}{\left(  k+\alpha\right)
^{2}}\sum_{k=n}^{m}\frac{1}{\left(  k+\alpha\right)  ^{3}}+\frac{1}{5}%
\sum_{k=n}^{m}\frac{1}{\left(  k+\alpha\right)  ^{5}}\\
\phi_{6}  &  =-\frac{1}{48}\left(  \sum_{k=n}^{m}\frac{1}{\left(
k+\alpha\right)  ^{2}}\right)  ^{3}+\frac{1}{18}\left(  \sum_{k=n}^{m}\frac
{1}{\left(  k+\alpha\right)  ^{3}}\right)  ^{2}+\frac{1}{8}\sum_{k=n}^{m}%
\frac{1}{\left(  k+\alpha\right)  ^{2}}\sum_{k=n}^{m}\frac{1}{\left(
k+\alpha\right)  ^{4}}-\frac{1}{6}\sum_{k=n}^{m}\frac{1}{\left(
k+\alpha\right)  ^{6}}%
\end{align*}
In passing by, our characteristic coefficients also enable to compute the
Stirling numbers $S_{m}^{\left(  k\right)  }$ via the generating function
(\ref{genfunc_stirling}) for large $m$ up to any order desired.

Let us proceed with $n=1$ (restricting $\alpha>-\frac{1}{2}$) and denote
$\gamma_{m}\left(  \alpha\right)  =\sum_{k=1}^{m}\frac{1}{\left(
k+\alpha\right)  }$, then
\[
\left(  m+1\right)  b_{m+1}(\frac{1}{2}+\alpha)=\prod_{k=1}^{m}\left(
k+\alpha\right)  \sum_{l=0}^{\infty}\phi_{l}\int_{-\frac{1}{2}}^{\frac{1}{2}%
}\left(  u+\alpha\right)  e^{u\gamma_{m}\left(  \alpha\right)  }u^{l+1}du
\]
The integral%
\[
\int_{-\frac{1}{2}}^{\frac{1}{2}}\left(  u+\alpha\right)  e^{u\gamma
_{m}\left(  \alpha\right)  }u^{l+1}du=\int_{-\frac{1}{2}}^{\frac{1}{2}%
}e^{u\gamma_{m}\left(  \alpha\right)  }u^{l+2}du+\alpha\int_{-\frac{1}{2}%
}^{\frac{1}{2}}e^{u\gamma_{m}\left(  \alpha\right)  }u^{l+1}du
\]
requires us to compute $I_{l}=\int_{-\frac{1}{2}}^{\frac{1}{2}}e^{au}u^{l}du$
for integer $l\geq0$. Partial integration leads to the recursion
\[
I_{l}=\int_{-\frac{1}{2}}^{\frac{1}{2}}e^{au}u^{l}du=\frac{1}{a2^{l}}\left(
e^{\frac{a}{2}}-\left(  -1\right)  ^{l}e^{-\frac{a}{2}}\right)  -\frac{l}%
{a}I_{l-1}%
\]
which, after iteration down to $I_{0}=\int_{-\frac{1}{2}}^{\frac{1}{2}}%
e^{au}du=\frac{e^{\frac{a}{2}}-e^{-\frac{a}{2}}}{a}$, leads to%
\[
I_{l}=\int_{-\frac{1}{2}}^{\frac{1}{2}}e^{au}u^{l}du=\frac{l!\left(
-1\right)  ^{l}}{a^{l+1}}\left(  e^{\frac{a}{2}}\sum_{j=0}^{l}\frac{\left(
-\frac{a}{2}\right)  ^{j}}{j!}-e^{-\frac{a}{2}}\sum_{j=0}^{l}\frac{\left(
\frac{a}{2}\right)  ^{j}}{j!}\right)
\]
Although the right-hand side seems to increase factorially with $l$, the
integral indicates that $\lim_{l\rightarrow\infty}I_{l}=0$. Thus, we obtain%
\[
\int_{-\frac{1}{2}}^{\frac{1}{2}}\left(  u+\alpha\right)  e^{u\gamma
_{m}\left(  \alpha\right)  }u^{l+1}du=\sum_{q=0}^{l+2}\left\{
\begin{array}
[c]{c}%
e^{\frac{\gamma_{m}\left(  \alpha\right)  }{2}}\left(  -1\right)  ^{q}\left(
\left(  \alpha+\frac{1}{2}\right)  l+1+\left(  2-q\right)  \alpha\right) \\
+\left(  -1\right)  ^{l}e^{-\frac{\gamma_{m}\left(  \alpha\right)  }{2}%
}\left(  \left(  \alpha-\frac{1}{2}\right)  l-1+\left(  2-q\right)
\alpha\right)
\end{array}
\right\}  \frac{2^{q-l-1}\left(  l+1\right)  !}{\left(  l+2-q\right)  !\left(
\gamma_{m}\left(  \alpha\right)  \right)  ^{q+1}}%
\]
Returning to $P=\frac{\left(  m+1\right)  b_{m+1}(\frac{1}{2}+\alpha)}%
{\prod_{k=1}^{m}\left(  k+\alpha\right)  }=\sum_{l=0}^{\infty}\phi_{l}%
\int_{-\frac{1}{2}}^{\frac{1}{2}}\left(  u+\alpha\right)  e^{u\gamma
_{m}\left(  \alpha\right)  }u^{l+1}du$ and after reversing of the $l$- and
$q$-sum, we obtain the expansion in inverse powers of $\gamma_{m}\left(
\alpha\right)  =\sum_{k=1}^{m}\frac{1}{\left(  k+\alpha\right)  }$,
\begin{align*}
P  &  =\frac{1}{2\gamma_{m}\left(  \alpha\right)  }\sum_{l=0}^{\infty}\phi
_{l}\left\{  e^{\frac{\gamma_{m}\left(  \alpha\right)  }{2}}\left(
\alpha+\frac{1}{2}\right)  +\left(  -1\right)  ^{l}e^{-\frac{\gamma_{m}\left(
\alpha\right)  }{2}}\left(  \alpha-\frac{1}{2}\right)  \right\}  \frac
{1}{2^{l}}\\
&  +\frac{1}{\left(  \gamma_{m}\left(  \alpha\right)  \right)  ^{2}}\sum
_{l=0}^{\infty}\phi_{l}\left\{  -e^{\frac{\gamma_{m}\left(  \alpha\right)
}{2}}\left(  \left(  \alpha+\frac{1}{2}\right)  l+1+\alpha\right)  +\left(
-1\right)  ^{l}e^{-\frac{\gamma_{m}\left(  \alpha\right)  }{2}}\left(  \left(
\alpha-\frac{1}{2}\right)  l-1+\alpha\right)  \right\}  \frac{1}{2^{l}}\\
&  +\sum_{q=0}^{\infty}\left(  \sum_{l=q}^{\infty}\phi_{l}\frac{\left(
l+1\right)  !\left\{  e^{\frac{\gamma_{m}\left(  \alpha\right)  }{2}}\left(
-1\right)  ^{q}\left(  \left(  \alpha+\frac{1}{2}\right)  l+1-q\alpha\right)
+\left(  -1\right)  ^{l}e^{-\frac{\gamma_{m}\left(  \alpha\right)  }{2}%
}\left(  \left(  \alpha-\frac{1}{2}\right)  l-1-q\alpha\right)  \right\}
}{2^{l-1}\left(  l-q\right)  !}\right)  \frac{2^{q}}{\left(  \gamma_{m}\left(
\alpha\right)  \right)  ^{q+3}}%
\end{align*}
With $\exp\left(  \sum_{j=2}^{\infty}\frac{\left(  -1\right)  ^{j-1}}{j}%
\sum_{k=n}^{m}\frac{1}{\left(  k+\alpha\right)  ^{j}}u^{j}\right)  =\sum
_{l=0}^{\infty}\phi_{l}u^{l}$ and $\prod_{k=n}^{m}\left(  1+\frac{u}{k+\alpha
}\right)  =\exp\left(  \sum_{j=1}^{\infty}\frac{\left(  -1\right)  ^{j-1}}%
{j}\sum_{k=n}^{m}\frac{1}{\left(  k+\alpha\right)  ^{j}}u^{j}\right)  $, we
observe that all $l$-sums in the last double sum are derivatives evaluated at
$u=\pm\frac{1}{2}$. For example, the first sum equals%
\begin{align*}
S_{1}  &  =\sum_{l=0}^{\infty}\phi_{l}\left\{  e^{\frac{\gamma_{m}\left(
\alpha\right)  }{2}}\left(  \alpha+\frac{1}{2}\right)  +\left(  -1\right)
^{l}e^{-\frac{\gamma_{m}\left(  \alpha\right)  }{2}}\left(  \alpha-\frac{1}%
{2}\right)  \right\}  \frac{1}{2^{l}}\\
&  =\left(  \alpha+\frac{1}{2}\right)  \prod_{k=1}^{m}\left(  1+\frac
{1}{2\left(  k+\alpha\right)  }\right)  +\left(  \alpha-\frac{1}{2}\right)
\prod_{k=1}^{m}\left(  1-\frac{1}{2\left(  k+\alpha\right)  }\right)
\end{align*}
We arrive at the expansion in powers of $\frac{1}{\gamma_{m}\left(
\alpha\right)  }$,
\begin{align}
P  &  =\frac{\left(  m+1\right)  b_{m+1}(\frac{1}{2}+\alpha)}{\prod_{k=1}%
^{m}\left(  k+\alpha\right)  }=\frac{\left(  \alpha+\frac{1}{2}\right)
\prod_{k=1}^{m}\left(  1+\frac{1}{2\left(  k+\alpha\right)  }\right)  +\left(
\alpha-\frac{1}{2}\right)  \prod_{k=1}^{m}\left(  1-\frac{1}{2\left(
k+\alpha\right)  }\right)  }{2\gamma_{m}\left(  \alpha\right)  }\nonumber\\
&  +\frac{1}{\left(  \gamma_{m}\left(  \alpha\right)  \right)  ^{2}}\sum
_{l=0}^{\infty}\phi_{l}\left\{  -e^{\frac{\gamma_{m}\left(  \alpha\right)
}{2}}\left(  \left(  \alpha+\frac{1}{2}\right)  l+1+\alpha\right)  +\left(
-1\right)  ^{l}e^{-\frac{\gamma_{m}\left(  \alpha\right)  }{2}}\left(  \left(
\alpha-\frac{1}{2}\right)  l-1+\alpha\right)  \right\}  \frac{1}{2^{l}%
}\nonumber\\
&  +\sum_{q=0}^{\infty}\sum_{l=q}^{\infty}\phi_{l}\left\{  e^{\frac{\gamma
_{m}\left(  \alpha\right)  }{2}}\left(  -1\right)  ^{q}\left(  \left(
\alpha+\frac{1}{2}\right)  l+1-q\alpha\right)  +\left(  -1\right)
^{l}e^{-\frac{\gamma_{m}\left(  \alpha\right)  }{2}}\left(  \left(
\alpha-\frac{1}{2}\right)  l-1-q\alpha\right)  \right\}  \frac{2^{q-l+1}%
\left(  l+1\right)  !}{\left(  l-q\right)  !\left(  \gamma_{m}\left(
\alpha\right)  \right)  ^{q+3}} \label{asymptotic_expansion_b_m(alpha)}%
\end{align}
In particular, for large $m$, where $\gamma_{m}\left(  \alpha\right)
=\sum_{k=1}^{m}\frac{1}{\left(  k+\alpha\right)  }=O\left(  \log\left(
m+a\right)  \right)  $, the expansion (\ref{asymptotic_expansion_b_m(alpha)})
shows that%
\begin{equation}
\frac{b_{m+1}(\frac{1}{2}+\alpha)}{m!}=\frac{\prod_{k=1}^{m}\left(
1+\frac{\alpha}{k}\right)  }{m+1}\frac{\left(  \alpha+\frac{1}{2}\right)
\prod_{k=1}^{m}\left(  1+\frac{1}{2\left(  k+\alpha\right)  }\right)  +\left(
\alpha-\frac{1}{2}\right)  \prod_{k=1}^{m}\left(  1-\frac{1}{2\left(
k+\alpha\right)  }\right)  }{2\sum_{k=1}^{m}\frac{1}{\left(  k+\alpha\right)
}}+O\left(  \frac{1}{\log^{2}m}\right)
\label{asymptotic_b(alpha_1/2)_large_m}%
\end{equation}
For $\alpha\rightarrow-\frac{1}{2}$ and $b_{m}(0)=c_{m}$, we find that
$c_{m}<0$ and that%
\[
\frac{c_{m+1}}{m!}=-\frac{\prod_{k=1}^{m}\left(  1-\frac{1}{2k}\right)  }%
{m+1}\frac{\prod_{k=1}^{m}\left(  1-\frac{1}{2\left(  k-\frac{1}{2}\right)
}\right)  }{2\sum_{k=1}^{m}\frac{1}{k-\frac{1}{2}}}+O\left(  \frac{1}{\log
^{2}m}\right)  =O\left(  \frac{1}{m\log m}\right)
\]
while for $\alpha=\frac{1}{2}$, $b_{m}\left(  1\right)  =\beta_{m}>0$ and
\[
\frac{b_{m+1}(1)}{m!}=\frac{\prod_{k=1}^{m}\left(  1+\frac{1}{2k}\right)
}{m+1}\frac{\prod_{k=1}^{m}\left(  1+\frac{1}{2\left(  k+\alpha\right)
}\right)  }{2\sum_{k=1}^{m}\frac{1}{k+\frac{1}{2}}}+O\left(  \frac{1}{\log
^{2}m}\right)  =O\left(  \frac{1}{m\log m}\right)
\]
Although $\frac{c_{m+1}}{m!}=O\left(  \frac{1}{m\log m}\right)  $ and
$\frac{b_{m+1}(1)}{m!}=O\left(  \frac{1}{m\log m}\right)  $, the products for
$c_{m}=b_{m}\left(  0\right)  $ are smaller than for $\beta_{m}=b_{m}\left(
1\right)  $, illustrating that the $\alpha=0$ case converges faster than the
$\alpha=1$ case (as in Fig. \ref{Fig_Stirlingoptimalaccuracy}).

\section{Verification of Theorem
\ref{theorem_Laplace_transform_infinite_convergent_factorial_expansions}}

\label{sec_inverse_Laplace_transform_factorial_series}Substituting the
factorial series (\ref{factorial_series_Laplace_transform_in_alpha_en_beta})
into the inverse Laplace transformation (\ref{inverse_Laplace}) and assuming
that summation and integration can be reversed, yields%
\[
f(t)=\sum_{m=0}^{\infty}m!\phi_{m}\left(  \alpha,\beta\right)  \frac{1}{2\pi
i}\int_{c-i\infty}^{c+i\infty}\frac{e^{zt}dz}{\prod_{k=0}^{m}(\beta
z+\alpha+k)}%
\]
For $\operatorname{Re}\left(  t\right)  \geq0$, $\operatorname{Re}\left(
\alpha\right)  >0$ and limiting ourselves to $\beta=1$, the contour can be
closed over the negative $\operatorname{Re}\left(  z\right)  $-plane, where
simple poles at $z=-\alpha-k$ are enclosed. Cauchy's residue theorem
\cite{Titchmarshfunctions} then indicates that%
\[
\frac{1}{2\pi i}\int_{c-i\infty}^{c+i\infty}\frac{e^{zt}dz}{\prod_{k=0}%
^{m}(z+\alpha+k)}=\sum_{j=0}^{m}\lim_{z\rightarrow-\left(  \alpha+j\right)
}\frac{\left(  z+\alpha+j\right)  e^{zt}}{\prod_{k=0}^{m}(z+\alpha+k)}%
=\sum_{j=0}^{m}\frac{e^{-\left(  \alpha+j\right)  t}}{\prod_{k=0;k\neq j}%
^{m}(k-j)}%
\]
With $\prod_{k=0;k\neq j}^{m}(k-j)=\left(  -1\right)  ^{j}j!\left(
m-j\right)  !$, we have%
\[
\frac{1}{2\pi i}\int_{c-i\infty}^{c+i\infty}\frac{e^{zt}dz}{\prod_{k=0}%
^{m}(z+\alpha+k)}=\frac{1}{m!}\sum_{j=0}^{m}\binom{m}{j}\left(  -1\right)
^{j}e^{-\left(  \alpha+j\right)  t}=\frac{e^{-\alpha t}}{m!}\left(
1-e^{-t}\right)  ^{m}%
\]
Thus, we obtain%
\[
f(t)=e^{-\alpha t}\sum_{m=0}^{\infty}\phi_{m}\left(  \alpha,1\right)  \left(
1-e^{-t}\right)  ^{m}%
\]
Introducing the form
(\ref{factorial_polynomial_in_alpha_reversed_for_Taylor_coefficients}) for
$\phi_{m}\left(  \alpha,1\right)  $ yields%
\[
f(t)=e^{-\alpha t}\sum_{m=0}^{\infty}\frac{\left(  1-e^{-t}\right)  ^{m}}%
{m!}\sum_{l=0}^{m}\frac{1}{l!}\left.  \frac{d^{l}f(t)}{dt^{l}}\right\vert
_{t=0}\sum_{k=0}^{m-l}\frac{\left(  k+l\right)  !}{k!}S_{m}^{(k+l)}\left(
-1\right)  ^{m-\left(  l+k\right)  }\alpha^{k}%
\]
After reversing the $m$- and $l$- sum,%
\[
f(t)=e^{-\alpha t}\sum_{l=0}^{\infty}\frac{1}{l!}\left.  \frac{d^{l}%
f(t)}{dt^{l}}\right\vert _{t=0}\sum_{m=l}^{\infty}\frac{\left(  1-e^{-t}%
\right)  ^{m}}{m!}\sum_{k=l}^{m}S_{m}^{(k)}\left(  -1\right)  ^{m-k}\frac
{k!}{\left(  k-l\right)  !}\alpha^{k-l}%
\]
we recognize that%
\[
\sum_{k=l}^{m}S_{m}^{(k)}\left(  -1\right)  ^{m-k}\frac{k!}{\left(
k-l\right)  !}\alpha^{k-l}=\frac{d^{l}}{d\alpha^{l}}\sum_{k=0}^{m}S_{m}%
^{(k)}\left(  -1\right)  ^{m-k}\alpha^{k}%
\]
The generating function (\ref{genfunc_stirling}) indicates that $\sum
_{k=0}^{m}S_{m}^{(k)}\left(  -1\right)  ^{m-k}\alpha^{k}=\prod_{k=0}%
^{m-1}(k+\alpha)=m!\left(  -1\right)  ^{m}\binom{-\alpha}{m}$ and we have%
\begin{align*}
f(t)  &  =e^{-\alpha t}\sum_{l=0}^{\infty}\frac{1}{l!}\left.  \frac{d^{l}%
f(t)}{dt^{l}}\right\vert _{t=0}\frac{d^{l}}{d\alpha^{l}}\sum_{m=l}^{\infty
}\binom{-\alpha}{m}\left(  e^{-t}-1\right)  ^{m}\\
&  =e^{-\alpha t}\sum_{l=0}^{\infty}\frac{1}{l!}\left.  \frac{d^{l}%
f(t)}{dt^{l}}\right\vert _{t=0}\left(  \frac{d^{l}}{d\alpha^{l}}\sum
_{m=0}^{\infty}\binom{-\alpha}{m}\left(  e^{-t}-1\right)  ^{m}-\frac{d^{l}%
}{d\alpha^{l}}\sum_{m=0}^{l-1}\binom{-\alpha}{m}\left(  e^{-t}-1\right)
^{m}\right)
\end{align*}
For any $\alpha$ and real $t\geq0$, the binomial sum $\sum_{m=0}^{\infty
}\binom{-\alpha}{m}\left(  e^{-t}-1\right)  ^{m}=\left(  1+e^{-t}-1\right)
^{-\alpha}=e^{\alpha t}$, while $\frac{d^{l}}{d\alpha^{l}}\sum_{m=0}%
^{l-1}\binom{-\alpha}{m}\left(  e^{-t}-1\right)  ^{m}=0$ because $\sum
_{m=0}^{l-1}\binom{-\alpha}{m}\left(  e^{-t}-1\right)  ^{m}$ is a polynomial
in $\alpha$ of degree $l-1$. Finally, with $\frac{d^{l}}{d\alpha^{l}}\left(
e^{\alpha t}\right)  =t^{l}e^{\alpha t}$, we return, indeed, to the Taylor
expansion of $f\left(  t\right)  $ around the point $t_{0}=0$.

\end{document}